\documentstyle[12pt,aps,floats,amssymb,emlines,bezier,eqsecnum]{revtex}
\tighten

\renewcommand{\thesection}{\arabic{section}}
\newcommand{\R}{{\if mm {\rm I}\mkern -3mu{\rm R}\else \leavevmode
     \hbox{I}\kern -.17em \hbox{R} \fi}}
\renewcommand{\R}{{\Bbb R}}
\newcommand{\C}{{\if mm {{\rm C}\mkern -15mu{\phantom{\rm t}\vrule}}
    \mkern +10mu \else \leavemode \hbox{I}\kern -.17em \hbox{C} \fi}}
\renewcommand{\C}{{\Bbb C}}
\newcommand{\N}{{\if mm {\rm I}\mkern -3mu{\rm N}\else \leavevmode
\hbox{I}\kern -.17em \hbox{N} \fi}}
\renewcommand{\N}{{\Bbb N}}
\newcommand{\Z}{{\if mm {\sf Z}\mkern -8mu{\sf Z}\else \leavevmode
\hbox{Z}\kern -.17em \hbox{Z} \fi}}
\renewcommand{\Z}{{\Bbb Z}}
\newcommand{\bB}{{\bf B}}

\newcommand{\bH}{{\bf H}}
\newcommand{\bI}{{\bf I}}
\newcommand{\bR}{{\bf R}}

\newcommand{\cD}{{\cal D}}
\newcommand{\cG}{{\cal G}}
\newcommand{\cH}{{\cal H}}
\newcommand{\cI}{{\cal I}}
\newcommand{\cM}{{\cal M}}
\newcommand{\cN}{{\cal N}}
\newcommand{\cO}{{\cal O}}
\newcommand{\cR}{{\cal R}}
\newcommand{\cQ}{{\cal Q}}

\newcommand{\cV}{{\cal V}}
\newcommand{\sD}{{\sf D}}
\newcommand{\sN}{{\sf N}}
\newcommand{\sP}{{\sf P}}
\newcommand{\sQ}{{\sf Q}}
\newcommand{\sm}{{m}}
\newcommand{\tH}{\tilde{H}}
\newcommand{\tX}{\tilde{X}}
\newcommand{\tsN}{\tilde{\sf N}}
\newcommand{\tsP}{\tilde{\sf P}}

\newcommand{\Sum}{\displaystyle\sum\limits}
\newcommand{\Int}{\displaystyle\int\limits}

\newcommand{\Min}{\mathop{{\rm min}}\limits}
\newcommand{\Max}{\mathop{{\rm max}}\limits}
\newcommand{\Inf}{\mathop{\rm inf}}
\newcommand{\Sup}{\mathop{\rm sup}}
\newcommand{\Supp}{\mathop{\rm supp}}

\newcommand{\opla}{\mathop{\oplus}\limits}
\newcommand{\plusp}{\mathop{\mathop{+}\limits^{\mbox{\bf.}}}\limits}
\newcommand{\reduction}[2]{#1 \biggr|_{#2}}
\newcommand{\Reduction}[2]{\left.#1\right|_{#2}}

\newcommand{\Lim}{\mathop{{\rm lim}}\limits}
\newcommand{\diag}{\mathop{\rm diag}}
\newcommand{\dist}{\mathop{\rm dist}}
\newcommand{\Img}{\mathop{\rm Im}}
\newcommand{\Real}{\mathop{\rm Re}}
\newcommand{\Inter}{\mathop{\rm Int}}
\newcommand{\ri}{{\rm i}}
\newcommand{\lal}{{\langle}}
\newcommand{\ral}{{\rangle}}
\newcommand{\D}{\displaystyle}
\newcommand{\be}{\begin{equation}}
\newcommand{\ee}{\end{equation}}

\newtheorem{theorem}{\sc Theorem}
\newtheorem{lemma}{\sc Lemma}

\newtheorem{corollary}{\sc Corollary}
\newtheorem{remark}{\sc Remark}
\catcode`\@=11
\def\draftlabel#1{{\@bsphack\if@filesw {\let\thepage\relax
   \xdef\@gtempa{\write\@auxout{\string
      \newlabel{#1}{{\@currentlabel}{\thepage}}}}}\@gtempa
   \if@nobreak \ifvmode\nobreak\fi\fi\fi\@esphack}
        \gdef\@eqnlabel{#1}}
\def\@eqnlabel{}
\def\@vacuum{}
\def\draftmarginnote#1{\marginpar{\raggedright\scriptsize\tt#1}}
\def\draft{\oddsidemargin -.5truein
        \def\@oddfoot{\sl preliminary draft \hfil
        \rm\thepage\hfil\sl\today\quad\militarytime}
        \let\@evenfoot\@oddfoot \overfullrule 3pt
        \let\label=\draftlabel
        \let\marginnote=\draftmarginnote
@
 \def\@eqnnum{(\theequation)\rlap{\kern\marginparsep\tt\@eqnlabel}%
\global\let\@eqnlabel\@vacuum}  }
\def\numberbysection{\@addtoreset{equation}{section}
        \def\theequation{\thesection.\arabic{equation}}
             \@addtoreset{definition}{section}
        \def\thedefinition{\thesection.\arabic{definition}}
             \@addtoreset{theorem}{section}
        \def\thetheorem{\thesection.\arabic{theorem}}
                     \@addtoreset{lemma}{section}
        \def\thelemma{\thesection.\arabic{lemma}}
             \@addtoreset{proposition}{section}
        \def\theproposition{\thesection.\arabic{proposition}}
             \@addtoreset{corollary}{section}
        \def\thecorollary{\thesection.\arabic{corollary}}
             \@addtoreset{remark}{section}
        \def\theremark{\thesection.\arabic{remark}}
}
\def\underline#1{\relax\ifmmode\@@underline#1\else
        $\@@underline{\hbox{#1}}$\relax\fi}
\catcode`@=12
\relax
\numberbysection
\begin{document}
\title{Operator interpretation of resonances
       arising in spectral problems for
       ${\bf 2}\times{\bf 2}$ operator matrices}
\author{ R. Mennicken, A. K. Motovilov%
\thanks{On leave of absence from the Laboratory of Theoretical Physics,
 Joint Institute  for Nuclear Research, Dubna, 141980, Russia}}
\address{Department of Mathematics, University of Regensburg,
           D-93040 Regensburg, Germany}
\date{March 17, 1998}
\maketitle

\bigskip

\bigskip

\small
\noindent {\bf Abstract.} We consider operator matrices
$
\bH=\left(\matrix{
  A_0        &   B_{01}  \cr
  B_{10}     &   A_{1}
}\right)
$
with self-adjoint entries $A_i$, $i=0,1,$ and bounded
$B_{01}=B_{10}^*$, acting in the orthogonal sum
\mbox{${\cal H}={\cal H}_0\oplus{\cal H}_1$}
of Hilbert spaces ${\cal H}_0$ and ${\cal H}_1$.
We are especially interested in the case where the
spectrum of, say, $A_1$ is partly or totally embedded into the
continuous spectrum of $A_0$ and the transfer function
$M_1(z)=A_1-z+V_1(z)$, where $V_1(z)=B_{10}(z-A_0)^{-1}B_{01}$,
admits analytic continuation (as an operator-valued
function) through the cuts along branches of the continuous
spectrum of the entry $A_0$ into the unphysical sheet(s) of the
spectral parameter plane. The values of $z$ in the unphysical
sheets where $M_1^{-1}(z)$ and consequently the resolvent
$(\bH-z)^{-1}$ have poles are usually called resonances.  A main
goal of the present work is to find non-selfadjoint operators
whose spectra include the resonances as well as to
study the completeness and basis properties of the resonance
eigenvectors of $M_1(z)$ in ${\cal H}_1$.  To this end we
first construct an operator-valued function $V_1(Y)$ on the space
of operators in ${\cal H}_1$ possessing the property:
$V_1(Y)\psi_1=V_1(z)\psi_1$ for any eigenvector $\psi_1$ of $Y$
corresponding to an eigenvalue $z$ and then study the equation
$
   H_1=A_1+V_1(H_1).
$
We prove the solvability of this equation even in the
case where the spectra of $A_0$ and $A_1$ overlap.
Using the fact that the root vectors of the solutions $H_1$
are at the same time such vectors for $M_1(z)$, we prove
completeness and even basis properties for the root vectors
(including those for the resonances).

\footnotesize


\bigskip

\medskip

\normalsize
\noindent {LANL E-print {\tt funct-an/9708001}.
Submitted to {\sl Mathematische Nachrichten.}}

\vspace*{0.7cm}

\section{Introduction}
\label{Intro}
In this paper we deal with $2\times2$ operator matrices
\begin{equation}
\label{twochannel}
\bH=\left(
\begin{array}{cc}
                     A_0              &          B_{01}          \\
                     B_{10}           &          A_1
\end{array}
\right),
\end{equation}
acting in an orthogonal sum  ${\cal H}={\cal H}_0\oplus{\cal H}_1$
of separable Hilbert spaces ${\cal H}_0$ and ${\cal H}_1$.  The
entries $A_0:{\cal H}_0\rightarrow{\cal H}_0$, and
$A_1:{\cal H}_1\rightarrow{\cal H}_1$, are assumed to be
self-adjoint operators with
domains ${\cal D}(A_0)$ and ${\cal D}(A_1)$, respectively.  It
is assumed that the couplings $B_{ij}:
\cH_j\rightarrow\cH_i$, $i,j=0,1$, $i{\neq}j$, are bounded
operators  (i.\,e., $B_{ij}\in\bB(\cH_j,\cH_i)$) and $B_{01}=B^*_{10}$.
Under these assumptions the matrix $\bH$ is a self-adjoint operator
in $\cH$ with domain $\cD(\bH)=\cD(A_0)\oplus\cD(A_1)$.

Operators of the form~(\ref{twochannel}) arise in many
of physical problems (see
e.\,g., \cite{Dashen,GlockleMuller,JaffeLow,KorchinShebeko,%
YaF88,KMMMP,MalyshevMinlos,MotJMPh91,Okubo,%
Pavlov1984,PavlovShushkov,Simonov}),
typically as a result of decomposing the Hilbert
space $\cH$ of a quantum system in two
``channel'' subspaces. The first one, say $\cH_0$,
may be interpreted as a space of ``external''
(for example hadronic as in \cite{Dashen,JaffeLow,Simonov})
degrees of freedom. The second one,
$\cH_1$, is associated with an ``internal''
(for example, quark \cite{Dashen,JaffeLow,Simonov})
structure of the system. We mention also that
spectral problems for a class
of $2\times2$ operator matrices arise in
magnetohydrodynamics \cite{Goedbloed,Lifschitz}.

In the spectral theory of operators of the form~(\ref{twochannel}) an
important role is played by the {\em transfer functions}
\begin{equation}
\label{trf}
M_i(z)=A_i-z+V_i(z), \quad i=0,1,
\end{equation}
where
\begin{equation}
\label{Poten}
V_i(z)=-B_{ij}R_j(z)B_{ji}, \quad j\neq i.
\end{equation}
Hereafter the notations $R_j(z)$ are used for the resolvents of
the operators $A_j$, \mbox{$R_j(z)=(A_j-zI_j)^{-1}$} where
$I_j$ stands for the identity operator in $\cH_j$. A particular
role of the transfer functions $M_i(z)$ can be understood
already from the fact that the resolvent $\bR(z)$ of the operator
$\bH$, \mbox{$\bR(z)=(\bH-z\bI)^{-1}$} where $\bI$ is the identity
operator in $\cH$, can be expressed explicitly in terms of
$M_0(z)$ or $M_1(z)$:
\begin{equation}
\label{ReprRM}
\begin{array}{rcl}
\bR(z) & = & \left(\begin{array}{ccc}
R_{00}(z)             &\quad&   R_{01}(z)\\
R_{10}(z)             &\quad&   R_{11}(z)
\end{array}\right)\\
 &\phantom{\cdot}& \\
& = & \left(\begin{array}{ccc}
M_0^{-1}(z)             &\quad&-M_0^{-1}(z)B_{01}R_1(z)\\
-R_1(z)B_{10}M_0^{-1}(z)&\quad&R_1(z)+R_1(z)B_{10}M_0^{-1}(z)B_{01}R_1(z)
\end{array}\right)\\
  &\phantom{\cdot}& \\
     & = & \left(\begin{array}{ccc}
R_0(z)+R_0(z)B_{01}M_1^{-1}(z)B_{10}R_0(z)&\quad&-R_0(z)B_{01}M_1^{-1}(z)\\
-M_1^{-1}(z)B_{10}R_0(z)                  &\quad& M_1^{-1}(z)
\end{array}\right).
\end{array}
\end{equation}
It follows from the representations~(\ref{ReprRM}) that $\bR(z)$
and, hence, its components $R_{ij}$, $i,j=0,1,$ may partly inherit
the singularities of the channel resolvents $R_0(z)$ and $R_1(z)$.
However, all the nontrivial singularities of $\bR(z)$, differing from
those of $R_0(z)$ and $R_1(z)$, are singularities of the
inverse transfer functions $R_{00}(z)=M_0^{-1}(z)$
and/or $R_{11}(z)=M_1^{-1}(z)$. Therefore, in studying the spectral
properties of the transfer functions one studies at the same
time the spectral properties of the initial operator matrix $\bH$.

Often the study of the spectral properties of the transfer
functions indeed turns out to be a simpler task than an
immediate study of the spectral problem for the
total matrix~(\ref{twochannel}). In particular, the described
reduction of the spectral problem $\bH\Psi=z\Psi$ for
an initial two-channel Hamiltonian $\bH$
of the form~(\ref{twochannel}) to
the channel spectral problems
\begin{equation}
\label{Schrodinger}
\biggl(A_i+V_i(z)\biggr)\psi^{(i)}=z\psi^{(i)}
\end{equation}
is common place in quantum physics where the perturbations
$V_i(z)$ are called energy-dependent potentials, energy--dependent
interactions {\em etc.} Regarding this subject see, e.\,g.,
the papers~\cite{McKellarMcCay,MotRemTMF,SchmidEW} discussing some
problems related to the use of the energy--dependent potentials in
the physics of few-body systems.

In the case where one of the spaces $\cH_0$ and $\cH_1$ is
finite-dimensional, say, the space $\cH_i$, the respective
transfer function $M_i(z)$ is also known as the
Liv\v{s}ic matrix~\cite{Livsic} (see Ref.~\cite{Howland}
for applications of the Liv\v{s}ic matrices to perturbation
theory and for further references).

In the papers~\cite{BraunMA,McKellarMcCay} the following
question was raised: {\em Is it possible to introduce an operator
$H_i$, $i=0,1,$ independent of the spectral parameter $z$, such
that its spectrum coincides with the spectrum of
Eq.~{\rm(\ref{Schrodinger})} while the eigenvector of $H_i$ is at the
same time an eigenvector of the transfer function $M_i$
{\rm(}i.\,e., $H_i\psi^{(i)}=z\psi^{(i)}$ implies
that~{\rm(\ref{Schrodinger})} holds{\rm)}?} Obviously,
having found such an operator
one would reduce the spectral problem for
the transfer-function $M_i(z)$ to the standard spectral problem
for the operator $H_i$ and, thus, the questions regarding
completeness and basis properties for the eigenvectors of $M_i$
could be answered in terms of the operator $H_i$ referring to
well known facts from operator theory.

A rigorous answer to the above question was found by
{\sc A.\,K.\,Motovilov}~\cite{MotRemTMF,MotSPbWorkshop,MotRem,MotJMPh91}
in the case where the spectra $\sigma(A_0)$ and $\sigma(A_1)$
of the entries $A_0$ and $A_1$ can be interwoven with each other but
must be strictly separated,
\begin{equation}
\label{SeparCondSigma}
\dist\{\sigma(A_0),\sigma(A_1)\}>0.
\end{equation}
To this end an operator-valued function $V_i(Y_i)$ on the
space of linear operators in $\cH_i$ was constructed
in~\cite{MotRemTMF,MotSPbWorkshop,MotRem,MotJMPh91} such that
\begin{equation}
\label{VzPsi}
V_i(Y_i)\psi^{(i)}=V_i(z)\psi^{(i)}
\end{equation}
for any eigenvector $\psi^{(i)}$ corresponding
to an eigenvalue $z$ of the operator $Y_i$.
The desired operator $H_i$ was
searched for as a solution of the operator equation
\begin{equation}
\label{MainOld}
   H_i=A_i+V_i(H_i), \qquad i=0,1.
\end{equation}
Notice that an equation of the form~(\ref{MainOld})
first appeared explicitly in the paper~\cite{BraunMA}
by {\sc M.\,A.\,Braun}.
Obviously, if $H_i$ is a solution of Eq.~(\ref{MainOld}) and
$H_i\psi^{(i)}=z\psi^{(i)}$ then, due to~(\ref{VzPsi}),
automatically
\mbox{$z\psi^{(i)}=\biggl(A_i+V_i(H_i)\biggr)\psi^{(i)}=
\biggl(A_i+V_i(z)\biggr)\psi^{(i)}$}
and, thus, for these $z$ and $\psi^{(i)}$ the
equality~(\ref{Schrodinger}) holds.
The solvability of the equation~(\ref{MainOld})
was announced in~\cite{MotSPbWorkshop}
and proved in~\cite{MotRemTMF,MotRem} under the assumption
\begin{equation}
\label{SeparCondOld}
\|B_{ij}\|_2<\D\frac{1}{2}\dist\{\sigma(A_0),\sigma(A_1)\}
\end{equation}
where $\|B_{ij}\|_2$ stands for the Hilbert-Schmidt norm
of the couplings $B_{ij}$. It was found~\cite{MotRemTMF,MotRem}
that the problem of constructing the operators $H_i$
is closely related to the problem of searching for
the invariant subspaces $\cG_i$, $i=0,1,$ of the
matrix $\bH$ which admit the graph representations,
\begin{equation}
\label{GraphRepr}
\begin{array}{rcl}
\cG_0 &=& \left\{
u\in\cH:\,\, u=\left(\begin{array}{c}
                         u_0   \\ Q_{10}u_0
               \end{array}\right),\,
u_0\in\cH_0
\right\}\,,   \\ && \phantom{.}\\
\cG_1 &=& \left\{
u\in\cH:\,\, u=\left(\begin{array}{c}
                        Q_{01} u_1   \\ u_1
               \end{array}\right),\,
u_1\in\cH_1
\right\}\,,
\end{array}
\end{equation}
with bounded $Q_{ji}:\cH_i\to\cH_j$ such that $Q_{ij}=-Q_{ji}^*$ and
$\cH=\cG_0\oplus\cG_1$.
The point is that under the assumption~(\ref{SeparCondOld}) the solutions
$H_i$, $i=0,1,$ of Eqs.~(\ref{MainOld}) read
\begin{equation}
\label{HiOld}
  H_i=A_i+B_{ij}Q_{ji}
\end{equation}
where $Q_{ji}$ are contractions just realizing
the above representations~(\ref{GraphRepr}).
The operators $Q_{ji}$ satisfy
the stationary Riccati equations
\begin{equation}
\label{Riccati}
 Q_{ji}A_i-A_j Q_{ji}+Q_{ji}B_{ij}Q_{ji}=B_{ji},
\qquad i,j=0,1,\quad j\neq i\,
\end{equation}
while the similarity transform $\bH'=\cQ^{-1}\bH\cQ$
with
$
\cQ=\left(\begin{array}{cc}
     I_0   &   Q_{01}   \\
     Q_{10}  &    I_1
\end{array}\right)
$
reduces the operator $\bH$ to the block-diagonal form
$\bH'=\diag\{H_0,H_1\}$. Under the
condition~(\ref{SeparCondOld}) the spectra $\sigma(H_0)$ and
$\sigma(H_1)$ do not intersect each other, i.\,e.,
\begin{equation}
\label{SigNonInt}
\sigma(H_0)\cap\sigma(H_1)=\emptyset,
\end{equation}
and
\begin{equation}
\label{SigUn}
\sigma(\bH)=\sigma(H_0)\cup\sigma(H_1)
\end{equation}
while the $H_i$, $i=0,1,$ represent  parts (spectral components)
of the operator $\bH$ in the corresponding invariant subspaces
$\cG_i$.

The idea of the block diagonalization of the $2\times2$ operator
matrices in terms of the invariant subspaces~(\ref{GraphRepr})
is a rather old one going back to the paper~\cite{Okubo} by {\sc
S.\,Okubo} (regarding applications of Okubo's approach to
particle physics see, e.\,g.,
\cite{GlockleMuller,KorchinShebeko} and Refs. cited therein). In
a mathematically rigorous way this idea was used by {\sc
V.\,A.\,Malyshev} and {\sc R.\,A.\,Minlos} in their
method~\cite{MalyshevMinlos} for the construction of invariant
subspaces for a class of selfadjoint operators in statistical
physics. Regarding a proof of solvability of the Riccati
equations~(\ref{Riccati}), the techniques of
Ref.~\cite{MalyshevMinlos} are restricted to the case where the
norms of the entries $B_{ij}$ are sufficiently small and the
separation condition~(\ref{SeparCondSigma}) holds, too.
Recently, the existence of invariant subspaces of the
form~(\ref{GraphRepr})  was proved by {\sc V.\,M.\,Adamyan} and
{\sc H.\,Langer}~\cite{AdL} for arbitrary bounded entries
$B_{ij}$ however assuming, instead of the
condition~(\ref{SeparCondSigma}), the essentially different
assumption that the spectrum of one of the entries $A_i$,
$i=0,1$ is situated strictly below the spectrum of the other
one, say
\begin{equation}
\label{CondAdL}
\Max\sigma(A_1)<\Min\sigma(A_0)\,.
\end{equation}
Soon, the result of~\cite{AdL} was extended by
{\sc V.\,M.\,Adamyan, H.\,Langer, R.\,Mennicken} and
{\sc J.\,Saurer}~\cite{AdLMSr} to the case where
\begin{equation}
\label{CondAdLMSr}
\Max\sigma(A_1)\leq\Min\sigma(A_0)\,
\end{equation}
and where the couplings $B_{ij}$ were allowed to be
unbounded operators such that,
for \mbox{$\alpha_0<\Min\sigma(A_0)$},  the product
\mbox{$(A_0-\alpha_0)^{-1/2}B_{01}$} makes sense as a
bounded operator. The conditions~(\ref{CondAdL}),
(\ref{CondAdLMSr}) were then somewhat weakened by {\sc
R.\,Mennicken} and {\sc A.\,A.\,Shkalikov}~\cite{MenShk} in the
case of a bounded entry $A_1$ and the same type of entries
$B_{ij}$ as in~\cite{AdLMSr}. Instead of the explicit
conditions~(\ref{CondAdL}), (\ref{CondAdLMSr}) on the spectra of
$A_i$, the paper~\cite{MenShk} uses a condition on the spectrum
of the transfer function $M_1(z)$ requiring the existence of a
regular point $\beta>\Min\sigma(A_1)$ for $M_1$ such that
\begin{equation}
\label{CondMenShk}
   M_1(\beta)\leq\alpha_1<0.
\end{equation}
This condition still allows one to prove the existence of the
invariant subspaces of $\bH$ in the
form~(\ref{GraphRepr})~\cite{MenShk}. It should be noted that
the condition~(\ref{CondMenShk}) may hold even in the case
where the spectra $\sigma(A_0)$ and $\sigma(A_1)$ weakly
overlap but that the requirement above regarding the unboundedness
of the entries $B_{ij}$ is strictly necessary in this case
(see~\cite{MenShk}). For all the cases considered
in~\cite{AdL,AdLMSr,MenShk} the relations~(\ref{SigNonInt})
and~(\ref{SigUn}) also hold true. However, if the coupling
$B_{ij}$ is unbounded, the Riccati equation for the operator
$Q_{ji}$ determining the representations~(\ref{GraphRepr})
must be written in general in a more complicated
form as compared to Eq.~(\ref{Riccati})
(see details in~\cite{AdLMSr,MenShk}). One can check nevertheless
that the spectral component $H_i$ of the matrix $\bH$ constructed
in~\cite{AdL,AdLMSr,MenShk} satisfies the equation~(\ref{MainOld}),
at least in the case where for $j\neq i$ the entry $A_j$
is bounded.

In the present work we study the equation~(\ref{MainOld}) in a
case which is totally different from the spectral situations
considered in~\cite{AdL,AdLMSr,MalyshevMinlos,MenShk,%
MotRemTMF,MotSPbWorkshop,MotRem};
namely, we suppose from the beginning that
\mbox{$\sigma(A_0)\cap\sigma(A_1)\neq\emptyset$.} In fact, we are
especially interested in the case where the spectrum of, say
$A_1$, is partly or totally embedded into the continuous spectrum
of $A_0$. Some remarks concerning the solvability of the
equation~(\ref{MainOld}) in this case may be found only in
Ref.~\cite{MotJMPh91}.

We work under the assumption that the coupling
operators $B_{ij}$ are such that the transfer function $M_1(z)$
admits analytic continuation, as an operator-valued
function, under the cuts along the branches of the absolutely
continuous spectrum $\sigma_{ac}(A_0)$ of the entry $A_0$.
Among other things, Sect.~\ref{Transfer_functions} includes
a detailed description of the conditions on $B_{ij}$ making such a
continuation of $M_1(z)$ possible.

The problem considered is closely related to the resonances
generated by the matrix $\bH$. Regarding a definition of the
resonance and a history of the subject see, e.\,g., the
books~\cite{AlbeverioBook,Exner,Baz,ReedSimonIII-IV}.  A recent
survey of the literature on resonances can be found
in~\cite{MotMathNach}. Throughout the paper we treat resonances
as the discrete spectrum of the transfer function $M_1(z)$ situated
in the so-called unphysical sheets of its Riemann surface. One
can find some definitions regarding the
unphysical sheets and the resonances in
Sect.~\ref{Transfer_functions}.

Sect.~\ref{SmainEq} starts with adjusting the definition of the
function $V_1(Y)$ of
Refs.~\cite{MotRemTMF,MotSPbWorkshop,MotRem} to the spectral
situation considered here. Since we deal with $\sm$
($1\leq\sm<\infty$) distinct intervals of the absolutely
continuous spectrum of the entry $A_0$, we get as a result
$2^\sm$ variants of the function $V_1(Y)$ and, consequently,
$2^\sm$ different variants of the equation~(\ref{MainOld}) which
read now as Eq.~(\ref{MainEq}). This circumstance
corresponds to the $2^\sm$ possible ways of realizing the analytic
continuation of the transfer function $M_1(z)$ under $m$
distinct cuts into the unphysical sheet(s) neighboring the
physical one. It should be stressed that in this paper we deal
only with the neighboring unphysical sheets. For convenience,
Eq.~(\ref{MainEq}) is referred to as the basic equation in
the following.
The solvability of this equation is proved%
\footnote{It should be noted that having solved the basic
equation~(\ref{MainEq}) one can find as well some formal
solutions for the Riccati equations~(\ref{Riccati}).  However, in
this case the formal solutions $Q_{ij}$ of~(\ref{Riccati}) can
not be treated in the conventional operator sense. A generalized
interpretation of these solutions as well as the construction of
the generalized invariant subspaces~(\ref{GraphRepr}) are beyond
the scope of the present work and will be a subject of another
paper.}
under the assumption~(\ref{Best}) recalling the
condition~(\ref{SeparCondOld}) but without already requiring the
entries $B_{ij}$ to be of the Hilbert-Schmidt class.  The
solutions of~(\ref{MainEq}) represent non-selfadjoint operators
the spectrum of which includes the resonances. In general, these
operators are not even dissipative.

In Sect.~\ref{SecFactor} we first prove the factorization
theorem (Theorem~\ref{factorization}) for the transfer function
$M_1(z)$.  It follows from this theorem that there exist certain
domains surrounding the set $\sigma(A_1)$ and lying partly in
the unphysical sheet(s) where the spectrum of $M_1$ is
represented only by the spectrum of the respective solutions of
the basic equation~(\ref{MainEq}). Since the root vectors of
these solutions are also root vectors for $M_1$,
this fact allows us to talk further, in
Sect.~\ref{RealEigen}--\ref{BasisnessGeneral}, about
completeness and basis properties%
\footnote{Note that these results recall those of the Lax--Phillips
scattering theory~\cite{LaxPhillips}
(see also~\cite{AdamyanLP,PavlovLP,PavlovFedorov} and Refs. therein)
where resonances appear as the spectrum of a dissipative
operator representing the generator of the compressed
evolution (semi)group and this implies completeness and
basis properties of the resonance states in a translationally
invariant subspace.}
of the root vectors of the transfer function $M_1$ corresponding
to its spectrum in the above domains, including the resonance
spectrum. In Sect.~\ref{SecFactor} we describe as well some
relations between the different solutions of~(\ref{MainEq}) and
some relations between their spectra which reflect the symmetry
of the resonance sets with respect to the real axis.

In Sect.~\ref{RealEigen} we pay special attention to the real
point spectrum of the solutions of the basic
equation~(\ref{MainEq}) and, thereby, to this part of the
spectrum of the transfer function $M_1$ as well. It is found in
this section that the real isolated eigenvalues
of all the considered solutions of~(\ref{MainEq}) are the same
and, moreover, the real
eigenvalues correspond to the same algebraic eigenspaces which
consist in this case only of eigenvectors. We prove the basis
property of these eigenvectors with respect to their closed
linear span.

In contrast to Sect.~\ref{Transfer_functions}--\ref{RealEigen}
we suppose in Sect.~\ref{Sfinite-dim} and~\ref{BasisnessGeneral}
that the entry $A_1$ only has discrete spectrum.

Sect.~\ref{Sfinite-dim} is devoted to a detailed consideration of
the case where the space $\cH_1$ is finite-dimensional.  In
particular, we describe in this section the relations between the
eigenprojections and eigennilpotents corresponding to the
resonances.

The results of Sect.~\ref{BasisnessGeneral} are
obtained under the assumption that the operator $A_1$ has a
compact resolvent. Here, to prove the completeness and basis
properties for the root vector systems of the solutions of
Eq.~(\ref{MainEq}), we rely mainly on the respective statements
regarding non-selfadjoint operators from the books by {\sc
I.\,C.\,Gohberg} and {\sc M.\,G.\,Krein}~\cite{GK} and by {\sc
T.\,Kato}~\cite{Kato}.

In Sect.~\ref{Simple-Example} we present an illustration
of the results obtained for a simple example going back
to one of the Friedrichs models~\cite{Fried} while
Appendices~\ref{app1} and~\ref{IntOpMer} contain
some auxiliary material used throughout the paper.

\bigskip

\bigskip

The authors thank Prof.\,\,M.\,M\"oller who had carefully read the 
manuscript and made a number of important remarks.  The authors are 
indebted to Dr.\,\,J.\,Lutgen for proofreading the manuscript, too.  
The support of this work by the Deutsche Forschungsgemeinschaft, 
the Russian Foundation for Basic Research and the Joint 
Institute for Nuclear Research (Dubna) is gratefully 
acknowledged.  One of the authors (A.\,K.\,M.) is much indebted 
to the Department of Mathematics, University of Regensburg, for 
the kind hospitality.

\section{Transfer function:
        Analytic continuation
        through the continuous spectrum}
\label{Transfer_functions}

The transfer function $M_i(z)$, $i=0,1$, considered on
the resolvent set $\rho(A_j)$ of the entry $A_j$, $j\neq i$,
represents a particular case of a holomorphic
operator-valued functions. In the present work we use the standard
definition of holomorphy of an operator-valued function
with respect to the operator norm topology. Namely:

Let $D$ be a domain in $\C$ and $\bB(\cG',\cG'')$
the Banach space of bounded operators between the
Hilbert spaces $\cG'$ and $\cG''$. A mapping
$T:\,D\to\bB(\cG',\cG'')$ is said to be a holomorphic
operator-valued function on $D$ if it is differentiable
at every $z\in D$ with respect to the operator norm
topology.

Let $T(z)=A+F(z)$ with $A:\cG'\to\cG''$ a closed
(in particular self-adjoint if $\cG'=\cG''\equiv\cG$)
operator on a domain $\cD(A)$
and $F:\,D\to\bB(\cG',\cG'')$. Such a function $T$ is called
holomorphic on $D$ if $F$ is holomorphic on $D$.

Each transfer function $M_i(z)$, $i=0,1$,
is holomorphic at least in the resolvent set $\rho(A_j)$
of the entry $A_j$, $j\neq i$. Since the inverse transfer
functions $M_i^{-1}(z)$ coincide with the respective
block components $R_{ii}(z)$ of the resolvent $\bR(z)$,
they are both holomorphic at least in the set $\rho(\bH)$.

One can extend to operator-valued functions the usual
definitions of the spectrum and its components. We recall these
definitions here, restricting ourselves to the cases above of
holomorphic functions $T(\cdot)$ on a domain $D\subset\C$ where
either the function $T(\cdot)$ is bounded itself,
$T:\,D\to\bB(\cG,\cG)$, or $T(\cdot)$ is a sum,
$T(\cdot)=A+F(\cdot)$, of a fixed closed operator $A$ in a
Hilbert space $\cG$ and a bounded function
$F:\,D\to\bB(\cG,\cG)$.  A point
$\lambda\in D$ is called a regular point of the function $T$
if $T^{-1}(\lambda)\in\bB(\cG,\cG)$ exists. The set $\rho(T)$
consisting of all the regular points $\lambda\in D$ is called
the resolvent set of the function $T$ in the domain $D$. The set
$\sigma(T)=D\setminus\rho(T)$ is called the spectrum of $T$ in
$D$. If $\mathop{\rm Ker}T(\lambda)\neq\{0\}$, $\lambda\in D$,
then one says the $\lambda$ is an eigenvalue of $T$ in $D$,
$\lambda\in\sigma_p(T)\cap D$. If
$x\in\mathop{\rm Ker}T(\lambda)$, $x\neq0,$ and, thus,
$T(\lambda)x=0,$ then such an $x$ is called an eigenvector
corresponding to the eigenvalue $\lambda$. The continuous spectrum
$\sigma_c(T)$ of the function $T$  in $D$ is introduced
as the set of all those points $\lambda\in D$ for which the
image $\cR(T(\lambda))$ does not coincide with its closure,
$\cR(T(\lambda))\neq\overline{\cR(T(\lambda))}$.
Obviously, the standard definitions for the spectrum of a closed
operator $A$ coincide with the definitions above
if one takes $T(z)=A-z$ and $D=\C$.

Let  $E_j$ be the spectral measure for the entry
$A_j$, $A_j=\D\int_{\sigma(A_j)}\lambda\,dE_j(\lambda)$\,, $j=0,1,$
\mbox{$\sigma(A_j)\subset\R$}. Then the functions $V_i(z)$ given
in~(\ref{Poten}) can be written
$$
    V_i(z)=B_{ij}\Int_{\sigma(A_j)}dE_j(\mu)\D\frac{1}{z-\mu}B_{ji}.
$$
Thus, it is convenient to introduce the quantities
$$
\cV_j(B)=\|B_{ij}\|_{E_j}^2
$$
where, by definition,
$$
\|B_{ij}\|_{E_j}^2=\|B_{ji}\|_{E_j}^2=\Sup_{\left\{\delta_k\right\}}
\Sum_k \|B_{ij}E_j(\delta_k)B_{ji}\|,
$$
with $\left\{\delta_k\right\}$ being a finite or countable
complete system of pairwise nonintersecting subsets of the
spectrum $\sigma(A_j)$ being measurable with respect to $E_j$
(i.\,e., $\delta_k$ are Borel subsets of $\sigma(A_j)$ such
that $\delta_k\cap\delta_l=\emptyset$, if $k\neq l$ and
$\mathop{\bigcup}\limits_k\delta_k=\sigma(A_j)$).  The number
$\cV_j(B)$ is called the variation of the operators $B_{ij}$
with respect to the spectral measure $E_j$. At the same time the
quantity  $\|B_{ij}\|_{E_j}=\|B_{ji}\|_{E_j}$ will be called the
norm of the operators $B_{01}$ and $B_{10}$ with respect to this
measure.  Some properties of such a norm are described below in
Appendix~\ref{app1}.  It follows from the results of this
appendix that $\cV_j(B)$ satisfies the estimates
$$
\|B_{ij}\|^2\leq \cV_j(B)\leq\|B_{ij}\|_2^2
$$
where $\|B_{ij}\|_2=\|B_{ji}\|_2$ is the Hilbert-Schmidt norm of
the couplings $B_{ij}$.  The equality $\cV_j(B)=\|B_{ij}\|^2$ is
attained in the case where $A_j$ is a multiple of the identity
operator. The equality $\cV_j(B)=\|B_{ij}\|_2^2$ holds
if $A_j$ possesses only a pure discrete spectrum which is at the same
time simple.

Note that along with the ``total'' variation
$\cV_j(B)$ we shall use the ``truncated'' variations
$$
    \reduction{\cV_j(B)}{\Delta}=\Sup_{\left\{\delta_k\right\}}
 \Sum_k \|B_{ij}E_j(\delta_k\cap\Delta)B_{ji}\|
$$
where $\Delta$ is a certain Borel subset of
$\sigma(A_j)$. Obviously, for any such $\Delta\subset\sigma(A_j)$
one has $\reduction{\cV_j(B)}{\Delta}\leq \cV_j(B)$.

Further, we shall suppose that at least one of the variations
$\cV_j(B)$, $j=0,1$, is finite; say, for example, the
variation $\cV_0(B)$:
\begin{equation}
\label{N0est}
   \cV_0(B)<\infty.
\end{equation}

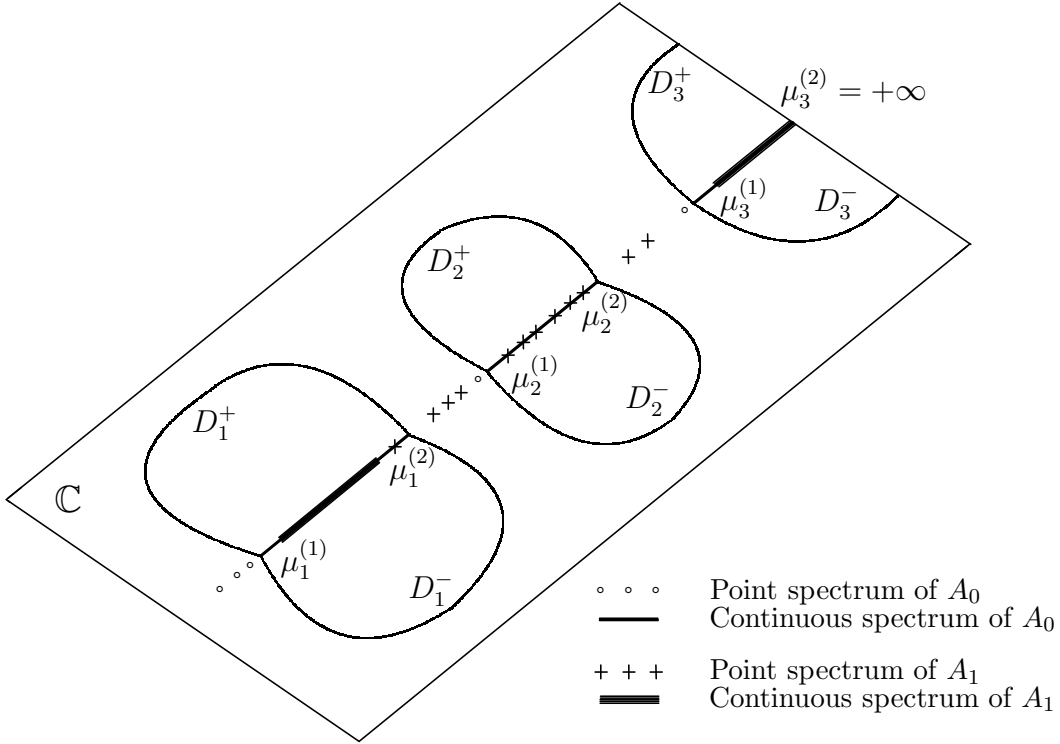
\begin{figure}
\centering
\unitlength=0.80mm
\special{em:linewidth .6pt}
\linethickness{.6pt}
\begin{picture}(160.00,122.83)
\emline{0.00}{40.00}{1}{58.33}{0.00}{2}
\emline{101.67}{122.67}{3}{160.00}{82.67}{4}
\emline{0.00}{40.49}{5}{101.67}{122.83}{6}
\emline{58.33}{0.00}{7}{159.33}{82.33}{8}
\emline{41.96}{30.83}{9}{66.53}{50.98}{10}
\emline{66.53}{50.98}{11}{66.53}{50.98}{12}
\emline{79.56}{61.49}{13}{97.85}{76.40}{14}
\emline{113.84}{89.44}{15}{130.46}{103.04}{16}
\emline{41.68}{31.04}{17}{42.22}{30.62}{18}
\emline{66.26}{51.20}{19}{66.80}{50.78}{20}
\emline{79.29}{61.69}{21}{79.83}{61.27}{22}
\emline{97.59}{76.61}{23}{98.13}{76.19}{24}
\emline{113.59}{89.63}{25}{114.13}{89.21}{26}
\put(40.21,29.31){\circle{0.83}}
\put(38.06,27.66){\circle{0.83}}
\put(35.19,25.31){\circle{0.83}}
\put(78.15,60.23){\circle{0.83}}
\put(112.47,88.31){\circle{0.83}}
\emline{157.40}{80.74}{27}{159.90}{82.78}{28}
\emline{2.17}{42.22}{29}{-0.31}{40.23}{30}
\emline{-0.31}{40.23}{31}{0.29}{39.79}{32}
\emline{45.14}{33.80}{33}{61.28}{47.02}{34}
\emline{45.47}{33.33}{35}{61.61}{46.55}{36}
\emline{117.45}{92.71}{37}{130.21}{103.07}{38}
\emline{117.18}{92.87}{39}{129.95}{103.23}{40}
\emline{45.01}{33.95}{41}{61.15}{47.17}{42}
\emline{45.60}{33.14}{43}{61.74}{46.37}{44}
\emline{118.14}{92.25}{45}{130.90}{102.60}{46}
\emline{117.87}{92.41}{47}{130.64}{102.76}{48}
\bezier{176}(114.00,89.33)(133.00,76.33)(147.67,90.67)
\bezier{180}(114.00,89.33)(94.67,105.00)(111.33,116.00)
\bezier{148}(98.00,76.33)(88.33,92.33)(71.67,85.00)
\put(34.33,53.00){\makebox(0,0)[cc]{$D_1^+$}}
\put(70.33,25.33){\makebox(0,0)[cc]{$D_1^-$}}
\put(73.33,79.66){\makebox(0,0)[cc]{$D_2^+$}}
\put(106.33,57.00){\makebox(0,0)[cc]{$D_2^-$}}
\put(110.00,109.33){\makebox(0,0)[cc]{$D_3^+$}}
\put(137.67,89.67){\makebox(0,0)[cc]{$D_3^-$}}
\put(10.00,40.33){\makebox(0,0)[cc]{\large$\C$}}
\put(49.33,30.33){\makebox(0,0)[cc]{$\mu_1^{(1)}$}}
\put(67.33,45.67){\makebox(0,0)[cc]{$\mu_1^{(2)}$}}
\put(87.66,61.00){\makebox(0,0)[cc]{$\mu_2^{(1)}$}}
\put(99.33,71.33){\makebox(0,0)[cc]{$\mu_2^{(2)}$}}
\put(122.33,90.33){\makebox(0,0)[cc]{$\mu_3^{(1)}$}}
\put(152.66,108.67){\makebox(0,0)[rc]{$\mu_3^{(2)}=+\infty$}}
\emline{41.85}{30.90}{49}{66.42}{51.05}{50}
\emline{42.06}{30.75}{51}{66.63}{50.90}{52}
\emline{113.72}{89.51}{53}{130.34}{103.11}{54}
\emline{113.96}{89.35}{55}{130.58}{102.95}{56}
\emline{70.65}{55.48}{57}{70.68}{53.20}{58}
\emline{73.13}{57.51}{59}{73.16}{55.23}{60}
\emline{75.24}{59.17}{61}{75.27}{56.89}{62}
\emline{74.29}{58.01}{63}{76.29}{58.01}{64}
\emline{72.17}{56.30}{65}{74.17}{56.30}{66}
\emline{69.71}{54.47}{67}{71.72}{54.47}{68}
\emline{64.26}{50.19}{69}{64.29}{47.91}{70}
\emline{63.31}{49.03}{71}{65.31}{49.03}{72}
\emline{79.45}{61.56}{73}{97.74}{76.47}{74}
\emline{79.67}{61.40}{75}{97.96}{76.31}{76}
\emline{90.86}{71.82}{77}{90.89}{69.54}{78}
\emline{93.41}{74.14}{79}{93.44}{71.86}{80}
\emline{95.53}{75.79}{81}{95.56}{73.51}{82}
\emline{94.57}{74.63}{83}{96.57}{74.63}{84}
\emline{92.46}{72.92}{85}{94.46}{72.92}{86}
\emline{89.93}{70.81}{87}{91.93}{70.81}{88}
\emline{83.05}{65.27}{89}{83.08}{62.99}{90}
\emline{85.59}{67.59}{91}{85.62}{65.31}{92}
\emline{87.71}{69.24}{93}{87.74}{66.96}{94}
\emline{86.76}{68.08}{95}{88.76}{68.08}{96}
\emline{84.64}{66.37}{97}{86.64}{66.37}{98}
\emline{82.11}{64.26}{99}{84.11}{64.26}{100}
\emline{103.07}{81.78}{101}{103.10}{79.50}{102}
\emline{106.29}{84.39}{103}{106.32}{82.11}{104}
\emline{105.34}{83.23}{105}{107.34}{83.23}{106}
\emline{102.11}{80.57}{107}{104.11}{80.57}{108}
\put(98.47,25.22){\circle{0.83}}
\emline{98.49}{12.72}{109}{98.52}{10.44}{110}
\emline{97.54}{11.56}{111}{99.54}{11.56}{112}
\emline{98.45}{7.23}{113}{108.15}{7.23}{114}
\emline{98.45}{7.54}{115}{108.15}{7.54}{116}
\emline{98.41}{6.36}{117}{108.11}{6.36}{118}
\emline{98.41}{6.67}{119}{108.11}{6.67}{120}
\emline{98.41}{6.94}{121}{108.11}{6.94}{122}
\emline{98.27}{20.48}{123}{107.98}{20.48}{124}
\emline{98.30}{20.12}{125}{108.01}{20.12}{126}
\emline{98.30}{20.28}{127}{108.01}{20.28}{128}
\put(116.54,20.34){\makebox(0,0)[lc]{{\small Continuous spectrum of $A_0$}}}
\put(116.50,6.96){\makebox(0,0)[lc]{{\small Continuous spectrum of $A_1$}}}
\bezier{184}(79.56,61.53)(55.79,73.75)(71.69,84.88)
\bezier{172}(79.51,61.43)(95.46,42.65)(109.87,53.29)
\bezier{188}(97.91,76.39)(123.89,67.71)(110.02,53.31)
\bezier{200}(41.91,30.78)(53.57,8.57)(74.24,22.70)
\bezier{272}(41.95,30.80)(7.56,41.34)(34.42,59.05)
\bezier{172}(66.53,51.05)(49.55,69.25)(34.38,59.18)
\bezier{228}(73.50,22.17)(93.78,40.94)(66.54,50.92)
\put(103.11,25.22){\circle{0.83}}
\put(107.63,25.22){\circle{0.83}}
\emline{103.13}{12.72}{129}{103.16}{10.44}{130}
\emline{102.18}{11.56}{131}{104.18}{11.56}{132}
\emline{107.65}{12.72}{133}{107.68}{10.44}{134}
\emline{106.70}{11.56}{135}{108.70}{11.56}{136}
\put(116.67,24.98){\makebox(0,0)[lc]{{\small Point spectrum of $A_0$}}}
\put(116.57,11.58){\makebox(0,0)[lc]{{\small Point spectrum of $A_1$}}}
\end{picture}
\bigskip
\caption{
An example of the spectral situation considered in the paper for
the case where $\sm=3$ and $\mu_1^{(1)}>-\infty$,
$\mu_3^{(2)}=+\infty$.
}
\label{figSpPic}
\end{figure}

We assume that the spectrum of the operator $A_1$ intersects
only the continuous spectrum of the operator $A_0$ and this
intersection is only realized on (every of) the pairwise
nonintersecting open intervals (see Fig.~\ref{figSpPic})
$\Delta_k^0=(\mu_k^{(1)},\mu_k^{(2)})\subset\sigma_c(A_0)$,
$\mu_k^{(1)}<\mu_k^{(2)}$,
\mbox{$k=1,2,\ldots,\sm$,}
\mbox{$\sm<\infty$,}
and \mbox{$-\infty\leq\mu_1^{(1)}$,}
\mbox{$\mu_\sm^{(2)}\leq+\infty$.}
Therefore, we assume that
$\Delta_k^0\cap\sigma(A_1)\neq\emptyset$ for all
$k=1,2,\ldots,\sm$ and $\sigma(A_1)\cap\sigma'(A_0)=\emptyset$
where $\sigma'(A_0)=\sigma(A_0)\setminus
\displaystyle\mathop{\bigcup}_{k=1}^\sm\Delta_k^0$ denotes the
remaining spectrum of $A_0$.


We shall suppose that the product
\begin{equation}
\label{KBproduct}
  K_B(\mu)\mathop{\mbox{\large$=$}}\limits^{\rm def}
  B_{10}E^0(\mu)B_{01}
\end{equation}
where $E^0(\mu)$ stands for the spectral function of $A_0$,
$E^0(\mu)=E_0\biggl((-\infty,\mu)\biggr)$,
is differentiable in $\mu$ for all $\mu\in\Delta_k^0$,
$k=1,2,\ldots,\sm,$ in the operator norm topology, i.\,e.
the limits
$$
\Lim_{\lambda\to\mu}\left\|\D\frac{K_B(\lambda)-K_B(\mu)}
{\lambda-\mu}-K'_B(\mu)\right\|=0,
\qquad \lambda,\mu\in\Delta_k^0,
$$
exist with $K'_B(\mu)\in\bB(\cH_1,\cH_1)$.
Obviously, the derivative $K'_B(\mu)$
is non-negative,
$$
 K'_B(\mu)\geq0,
$$
since $K_B(\mu)$ is a non-decreasing function.
Differentiability of $K_B(\mu)$ means that the continuous
spectrum of the entry $A_0$ includes, in each $\Delta_k^0$,
$k=1,2,\ldots,\sm,$ a branch of the absolutely continuous%
\footnote{Recall, for convenience of the reader, definition of
the absolutely continuous spectrum of a
selfadjoint operator $A$ acting in a Hilbert space $\cG$.  Let
$E^A(\mu)$ be the spectral function of $A$,
$E^A(\mu)=E_A\biggl((-\infty,\mu)\biggr)$, $\mu\in\R$, where
$E_A$ stands for the spectral measure of $A$. Denote by
$\cG_{ac}$ the invariant subspace of $A$ consisting of all the
vectors $f\in\cG$ for which the function $\mu\rightarrow\lal
E^A(\mu)f,f\ral$ is absolutely continuous in $\mu\in\R$. Then
the spectrum of the restriction $\reduction{A}{\cG_{ac}}$ of $A$
on the subspace $\cG_{ac}$ is called absolutely continuous
spectrum $\sigma_{ac}(A)$ of the operator $A$.  Also, one says
the subspace $\cG_{ac}$ corresponds to the spectrum
$\sigma_{ac}(A)$. For more details see, e.\,g.,
Ref.~\cite{Kato}, \S\,X.1.2, and Ref.~\cite{ReedSimonI},
Section~VII.2.  It should be noted that in most of physical
applications all the continuos spectrum of selfadjoint operators
involved is typically absolutely continuous.}
spectrum $\sigma_{ac}(A_0)$:  For any Borel subset
$\delta\subset\Delta^0_k$ and for any $u_1\in\cH_1$ the vector
$E_0(\delta)B_{01}u_1$ belongs to the invariant subspace
$\cH_0^{ac}\subset\cH_0$ of $A_0$ corresponding to the
absolutely continuous spectrum $\sigma_{ac}(A_0)$ and
$E_0(\delta)B_{01}u_1=E^{ac}_0(\delta)B_{01}u_1$ where
$E^{ac}_0$ denotes the part of the spectral measure $E_0$
corresponding to $\sigma_{ac}(A_0)$.  Obviously,
$$
   \reduction{\cV_0(B)}{\Delta_k^0}=\Int_{\Delta_k^0}d\mu\,\|K'_B(\mu)\|.
$$

Further, we suppose that the function $K'_B(\mu)$ is continuous
within the closed intervals $\overline{\Delta_k^0}$ and,
moreover, that it admits analytic continuation from each of
these intervals to a simply connected domain situated, say, in
$\C^-$.
For the interval $\Delta_k^0$, let this domain be
called $D_k^-$ (see Fig~\ref{figSpPic}). We assume
that the boundary of each domain $D_k^-$ includes the
entire spectral interval $\Delta_k^0$.
Let $K'^{(k)}_B(\mu)$ denote
the continuation of $K'_B$ from $\Delta_k^0$ into $D_k^-$.
The presence of the index $k$ in this notation is related to the
fact that in general the continuation of $K'_B$ to the same
domain of $\C$ can be different if one starts from different
intervals $\Delta_k^{0}$. Thus, the notation $K'^{(k)}_B(\mu)$
relates to the distinct branches of the function $K'_B$.  In the
case where $D_j^{-}\cap D_k^{-}\neq\emptyset$ for $j\neq k$
and $K'^{(j)}_B(\mu)\neq K'^{(k)}_B(\mu)$ for $\mu\in
D_j^{-}\cap D_k^{-}$ the points $z\in D_j^{-}$
and $z\in D_k^{-}$ must be considered as distinct (namely,
one must assume that these points belong to different sheets
of the Riemann surface of the function $K'_B$). To avoid
unnecessary complications regarding such a treatment of the
different branches of $K'_B$, we shall further assume that {\em the
domains $D_k^-$ for the different $k$ do not intersect each
other,} i.\,e.,
\begin{equation}
\label{DnotIntersect}
   D_j^-\cap D_k^-=\emptyset, \quad j\neq k\,.
\end{equation}
It is implied that one can always consider more narrow initial
holomorphy domains of the function $K'_B$ if it is necessary.
The assumption~(\ref{DnotIntersect}) allows us to drop the
branch identification (the index $k$) in the notation above for
the analytic continuation of $K'_B$ and this is done throughout
the paper.

Since $K'_B(\mu)$ represents a self-adjoint operator for
$\mu\in\Delta_k^0$ and $\Delta_k^0\subset\R$, the function
$K'_B(\mu)$ also automatically admits analytic continuation from
$\Delta_k^0$ into the domain $D_k^+$, symmetric to $D_k^-$ with
respect to the real axis, $D_k^+=\{z:\,\overline{z}\in D_k^-\}$.
For the continuation into $D_k^+$ we retain the same notation
$K'_B(\mu)$. The selfadjoitness of $K'_B(\mu)$ for
$\mu\in\Delta_k^{0}$ implies \begin{equation} \label{KBsym}
[K'_B(\mu)]^*=K'_B(\bar{\mu}), \quad \mu\in D_k^\pm\,.
\end{equation}
Also, we shall always suppose the $K'_B(\mu)$
satisfies the following H\"older condition at the end points
$\mu_k^{(1)}$, $\mu_k^{(2)}$ of the spectral intervals
$\Delta_k^0$, $$ \|K'_B(\mu)-K'_B(\mu_k^{(i)})\|
\leq C|\mu-\mu_k^{(i)}|^\gamma,
\qquad i=1,2,\quad\mu\in D_k^\pm,
$$ with some positive $C$ and $\gamma$.

Let $l=(l_1,l_2,\ldots,l_\sm)$ be a multi-index having the
components $l_k=+1$ or $l_k=-1$, $k=1,2,\ldots,\sm$. In what
follows we consider the domains
$D_l=\mathop{\bigcup}\limits_{k=1}^\sm D_k^{l_k}$ where
$D_k^{l_k}$ are the holomorphy domains of $K'_B$ described
above.  Let $\Gamma_k^{l_k}$ be a rectifiable Jordan curve in
$D_k^{l_k}$ resulting from continuous deformation of the
interval $\Delta^0_k$, the end points of this interval being
fixed (except
$\mu_1^{(1)}=-\infty$ and $\mu_\sm^{(2)}=+\infty$ which are
allowed, if this is possible, to be shifted respectively to
$\tilde{\mu}_1^{(1)}=-\infty+\ri y_{-\infty}\in D_1^{l_1}$ and
$\tilde{\mu}_\sm^{(2)}=+\infty+\ri y_{+\infty}\in D_m^{l_m}$
with some real $y_{-\infty}=l_1|y_{-\infty}|$ and
$y_{+\infty}=l_m|y_{+\infty}|$).  With the exception of the end
points, the closure $\overline{\Gamma}_k^{l_k}$ 
of the contour $\Gamma_k^{l_k}$ has no other common points
with the set $\sigma_c(A_0)$.  Note that under the
condition~(\ref{DnotIntersect})
$\Gamma_j^{\pm}\cap\Gamma_k^{\pm}=\emptyset$ for any
$j,k\in\{1,2,\ldots,m\}$ such that $j\neq k$.  By $\Gamma_l$,
$l=(l_1,l_2,\ldots,l_\sm)$, we shall denote the union of the
contours $\Gamma_k^{l_k}$,
$
\Gamma_l=\mathop{\bigcup}\limits_{k=1}^\sm\Gamma_k^{l_k}.
$

Also, we extend the definition of the variation $\cV_0(B)$
to the set $\sigma'(A_0)\cup\Gamma_l$ by introducing
the modified variation
\begin{equation}
\label{NBnorm}
\cV_0(B,\Gamma_l)=\reduction{\cV_0(B)}{\sigma'(A_0)}+
\displaystyle\int\limits_{\Gamma_l}|d\mu|\,\|K'_B(\mu)\|
\end{equation}
with $|d\mu|$ Lebesgue measure on $\Gamma_l$.
It is clear that if the length $\ell_{\Gamma_l}$ of the curve
$\Gamma_l$ is finite (in the case where the set
$\displaystyle\mathop{\cup}_{k=1}^\sm\Delta^0_k$ is bounded) the
value $\cV_0(B,\Gamma_l)$ is also finite,
$$
\cV_0(B,\Gamma_l)\leq\reduction{\cV_0(B)}{\sigma'(A_0)}+
{\ell_{\Gamma_l}}\cdot\Max_{\mu\in\Gamma_l}\|K'_B(\mu)\|.
$$
We suppose that the operators $B_{ij}$ are such that
there exists a contour (contours) $\Gamma_l$
where the value $\cV_0(B,\Gamma_l)$ is finite,
\begin{equation}
\label{Bbound}
\cV_0(B,\Gamma_l)<\infty,
\end{equation}
including also the case of the unbounded set
$\displaystyle\mathop{\cup}_{k=1}^\sm\Delta^0_k$.  It is
assumed that the inequality~(\ref{Bbound}) must hold during
the reverse continuous deformation of the contour $\Gamma_l$
back to the set
$\displaystyle\mathop{\cup}_{k=1}^\sm\Delta^0_k$.  The contours
$\Gamma_l$ satifying the condition~(\ref{Bbound}) are said to be
{\em $K_B$-bounded} contours.

In the following we deal mainly with the analytic continuation of the
transfer function $M_1(z)$ and its inverse, $[M_1(z)]^{-1}$,
through the spectral intervals $\Delta^0_k$ into the domains
$D_l$. Under the assumed conditions the function $M_1(z)$
admits such a continuation in the conventional sense, i.\,e., as
an operator-valued function. Namely, the following statement
holds.

\begin{lemma}\label{M1-Continuation}
The analytic continuation of the transfer function $M_1(z)$,
$z\in\C\setminus{\sigma(A_0)}$,
through the spectral intervals $\Delta^0_k$ into the subdomain
$D(\Gamma_l)\subset D_l$ bounded by the set
$\mathop{\bigcup}\limits_{k=1}^{\sm}\Delta^0_k$ and
a $K_B$-bounded contour $\Gamma_l$ is given by
\begin{equation}
\label{Mcmpl}
    M_1(z,\Gamma_l)=A_1-z+V_1(z,\Gamma_l)
\end{equation}
where
\begin{eqnarray}
\label{MGamma}
V_1(z,\Gamma_l) &=&\Int_{\sigma'(A_0)\cup\Gamma_l}
K_B(d\mu)\,\,\D\frac{1}{z-\mu}   \\
\nonumber
 &\mathop{\mbox{\large$=$}}\limits^{\rm def} &
\Int_{\sigma'(A_0)}B_{10}E_0(d\mu)B_{01}\,\D\frac{1}{z-\mu} +
\Int_{\Gamma_l}
d\mu\,K'_B(\mu)\,\D\frac{1}{z-\mu}.
\end{eqnarray}
For $z\in D_k^{l_k}\cap D(\Gamma_l)$ the function
$M_1(z,\Gamma_l)$ may be written as
\begin{equation}
\label{M1Gresidue}
  M_1(z,\Gamma_l)=M_1(z)+2\pi\ri\,l_k K'_B(z).
\end{equation}
\end{lemma}

\noindent P~r~o~o~f~.~~The proof is reduced to the observation
that the function $M_1(z,\Gamma_l)$ is holomorphic
for $z\in\C\setminus[\sigma'(A_0)\cup\Gamma_l]$
and coincides with $M_1(z)$ for
$z\in\C\setminus[\sigma'(A_0)\cup\overline{D(\Gamma_l)}]$.
Eq.~(\ref{M1Gresidue}) is obtained from~(\ref{MGamma})
using the Residue Theorem. %
{\nopagebreak\mbox{\phantom{MMMM}}\hfill $\Box$\par\addvspace{0.25cm}}

\begin{remark}
\label{M1M1Gamma}
The definition Eq.~{\rm(\ref{MGamma})} defines the function
$V_1(z,\Gamma_l)$ and, hence, via Eq.~{\rm(\ref{Mcmpl})} the
function $M_1(z,\Gamma_l)$  for
$z\in\C\setminus(\sigma'(A_0)\cup\Gamma_l)$, the values of
$V_1(z,\Gamma_l)$ for such $z$ being bounded operators in
$\cH_1$.  As mentioned above the inverse transfer function
$M_1^{-1}(z)$ coincides with the block component $R_{11}(z)$ of
the resolvent ${\bf R}(z)$ and, thus, it is bounded and
holomorphic in $z\in\C\setminus\sigma({\bf H})$.  Since
$M_1(z,\Gamma_l)=M_1(z)$ for
$z\in\C\setminus[\sigma'(A_0)\cup\overline{D(\Gamma_l)}]$, one
concludes that $[M_1(z,\Gamma_l)]^{-1}$ exists and is bounded
and holomorphic in $z$ at least for $z\in\C\setminus[\sigma({\bf
H})\cup\overline{D(\Gamma_l)}]$.
\end{remark}

It is clear that by varying the contour $\Gamma_l$ one can represent
the continuation of the transfer function $M_1(z)$ in the
form~(\ref{MGamma}) for any subset of the domain $D_l$.
One notes also that if the subdomains $D(\Gamma^{(1)}_l)$ and
$D(\Gamma^{(2)}_l)$  correspond to two different
$K_B$-bounded contours $\Gamma^{(1)}_l$ and $\Gamma^{(2)}_l$,
then automatically
$V_1(z,\Gamma^{(1)}_l)=V_1(z,\Gamma^{(2)}_l)$ and, hence,
$M_1(z,\Gamma^{(1)}_l)=M_1(z,\Gamma^{(2)}_l)$ for
\mbox{$z\in{D}(\Gamma^{(1)}_l){\cap}D(\Gamma^{(2)}_l)$}
because of the uniqueness of the analytic continuation.

The formula~(\ref{M1Gresidue}) shows that in general the
transfer function $M_1$ has a multi-sheeted Riemann surface.
Properties of this surface such as the number of sheets, the
presence of branching points in addition to the points
$\mu_k^{(1)},\mu_k^{(2)}$, $k=1,2,\ldots,\sm$ (if $K'_B$ is
considered in a larger domain than
$\mathop{\cup}\limits_l D_l$) {\em etc.} are determined by the
analytic properties of the function $K'_B(\mu)$ itself.  The
sheet of the complex plane where the transfer function $M_1(z)$
is considered together with the resolvent $\bR(z)$ initially
is said to be the {\em physical sheet}.  The remaining sheets of the
Riemann surface are said to be {\em unphysical sheets.}

In the present work we only deal with the unphysical
sheets neighbouring the physical one, i.\,e.,
with the sheets connected through the intervals
$\Delta_k^0$ for some $k\in\{1,2,\ldots,\sm\}$
immediately to the physical sheet.
The index $l=(l_1,l_2,\ldots,l_m)$ can be considered
as an identifier of the neighboring sheet. It should be
noted however that some of these sheets
can turn out to be identical to each other if
one is able to consider a wider domain than
$\mathop{\cup}\limits_l D_l$ but all this
depends on a concrete form for the function $K'_B$
and we do not touch on this subject.

Regarding the total resolvent $\bR(z)$, it may admit
continuation only in a generalized sense.  First, one
can suppose that there exists a dense subset
$\tilde{\cH}_0$ of ${\cH}_0$ such that for any
$u_0,v_0\in\tilde{\cH}_0$ the bilinear form
$\langle{E^0}(\mu)u_0,v_0\rangle$ admits analytic continuation
in variable $\mu$ as a holomorphic function from each interval
$\Delta^0_k$ into the respective domain $D^{l_k}$.  So that the
previous assumption concerning the holomorphy of the function
$K'_B(\mu)$ implies that $B_{01}u_1$ is an element of the
subset $\tilde{\cH}_0$ for any $u_1\in\cH_1$. Obviously,
the analytic continuation of the form
$\langle{R_0}(z)u_0,v_0\rangle$ into $D_l$ reads as
(cf.~Lemma~\ref{M1-Continuation})
$$
\Reduction{\langle{R_0}(z)u_0,v_0\rangle}{D_l}=
\Int_{\sigma'(A_0)}\D\frac{d\langle{E^0}(\mu)u_0,v_0\rangle}{\mu-z}+
\Int_{\Gamma_l}d\mu\, \D\frac{K'_{u_0,v_0}(\mu)}{\mu-z}
$$
where $K'_{u_0,v_0}$ denotes the derivative of
the analytic continuation $K_{u_0,v_0}(\mu)$ of
the form $\langle{E^0}(\mu)u_0,v_0\rangle$.
Using the Residue Theorem one can verify
$$
\reduction{\langle{R_0}(z)u_0,v_0\rangle}{D_k^{l_k}\cap D_l}=
\langle{R_0}(z)u_0,v_0\rangle
+2\pi\ri\,l_k\,K'_{u_0,v_0}(z)
$$
where $\langle{R_0}(z)u_0,v_0\rangle$ stands for the
conventional bilinear form of the resolvent $R_0(z)$, i.\,e., this
form is taken for $R_0(z)$ in the phyisical sheet.

The analytic continuation of the resolvent $\bR(z)$ is
understood in terms of the continuation of the bilinear form
$\lal \bR(z)u,v\ral$ where $u=(u_0,u_1)$, $v=(v_0,v_1)$ are
elements of $\cH$ with $u_0,v_0\in\tilde{\cH}_0$ and
$u_1,v_1\in\cH_1$.  It follows from the
representation~(\ref{ReprRM}) that such a generalized
continuation is indeed possible if the function
$K'_{u_0,v_0}(\mu)$ is holomorphic in $D_l$ for any
$u_0,v_0\in\tilde{\cH}_0$. The holomorphy domain of the continuation
of $\bR(z)$ into $D_l$ thus has to coincide with just such
a domain for the continuation of the inverse transfer function
$R_{11}(z)=[M_1(z)]^{-1}$.

The spectral problem for the continued transfer function
$M_1(z,\Gamma)$, that is, the problem
\begin{equation}
\label{init}
[A_1+V_1(z,\Gamma)]\,u_1=z\,u_1\,, \quad u_1\in\cH_1,
\end{equation}
will be referred to in the following as the {\em initial spectral problem}.

\section{The basic equation. Solutions
$H_1^{(\lowercase{l})}$}
\label{SmainEq}

If an operator-valued function
$T:\,\sigma'(A_0)\cup\Gamma\rightarrow\bB(\cH_1,\cH_1)$
satisfies the Lipschitz condition [the inequality
(\ref{Lipschitz}) of Appendix~\ref{IntOpMer}] on $\sigma'(A_0)$
and is continuous and bounded on a $K_B$-bounded contour
$\Gamma$,
\begin{equation}
\label{supT}
\|T\|_{\infty,\Gamma}=
\Sup_{\mu\in\sigma'(A_0)\cup\Gamma}\|T(\mu)\|<\infty,
\end{equation}
then the integral
\begin{equation}
\label{Xi}
\Int_{\sigma'(A_0)\cup\Gamma}K_B(d\mu)\,T(\mu)
\mathop{\mbox{\large$=$}}\limits^{\rm def}
\Int_{\sigma'(A_0)}B_{10}E_0(d\mu)B_{01}T(\mu)+
\Int_{\Gamma}
d\mu\,K'_B(\mu)\,T(\mu),
\end{equation}
exists in the sense of the operator norm topology
(see Appendix~\ref{IntOpMer})
and
\begin{equation}
\label{Lem1}
\biggl\|\,\,\Int_{\sigma'(A_0)\cup\Gamma}K_B(d\mu)T(\mu)\biggr\|
\leq \cV_0(B,\Gamma)\cdot\|T\|_{\infty,\Gamma}.
\end{equation}
In particular, if $T(z)$ is the resolvent of an operator
$Y$, \mbox{$T(z)=(Y-z I_1)^{-1}$,} the spectrum of which
has no common points with
$\sigma'(A_0)\cup\Gamma$, then one can define the operator
\begin{equation}
\label{V1}
V_1(Y,\Gamma)=
\Int\limits_{\sigma'(A_0)\cup\Gamma}
K_B(d\mu)(Y-\mu{I_1})^{-1}.
\end{equation}
This operator is bounded, $V_1(Y,\Gamma)\in\bB(\cH_1,\cH_1)$,
and, because of~(\ref{Lem1}), its norm admits the estimate
\begin{equation}
\label{V1Yest}
\|V_1(Y,\Gamma)\|\leq \cV_0(B,\Gamma_l)\cdot
\Sup_{\mu\in\sigma'(A_0)\cup\Gamma}\|(Y-\mu{I_1})^{-1}\|.
\end{equation}
According to the definition~(\ref{V1}), the operator-valued
function $V_1(Y,\Gamma)$ of the operator variable
$Y:\cH_1\rightarrow\cH_1$ possesses the following
important property: If $u_1\in\cH_1$ is an
eigenvector of $Y$ corresponding to an eigenvalue  $z$,
$Yu_1=zu_1$, then
\begin{equation}
\label{V1Action}
V_1(Y,\Gamma)u_1=
\Int\limits_{\sigma'(A_0)\cup\Gamma}
K_B(d\mu)\frac{1}{z-\mu}u_1
\equiv V_1(z,\Gamma)u_1.
\end{equation}

In what follows we consider the equation%
\footnote{
In the case where the spectra of $A_0$ and $A_1$ have no
intersection, the equation~(\ref{MainEq}) is reduced to the
operator Riccati equation~(\ref{Riccati}) (see Sect.~\ref{Intro};
 for details see Refs.~\cite{MotRemTMF,MotSPbWorkshop,MotRem}).}
\begin{equation}
\label{MainEq}
Y=A_1+V_1(Y,\Gamma).
\end{equation}
We deal with this equation, since it possesses the following
characteristic property: If an operator $H_1$ is a
solution of~(\ref{MainEq}) and $u_1$ is an eigenvector of $H_1$
corresponding to an eigenvalue $z$, $H_1 u_1=zu_1$, then
automatically (cf.~Sect.~\ref{Intro})
\begin{equation}
\label{Return}
zu_1=A_1 u_1+V_1(H_1,\Gamma)u_1=A_1 u_1+V_1(z,\Gamma)u_1.
\end{equation}
This implies that any eigenvalue $z$ of such an operator $H_1$
is automatically an eigenvalue for the initial spectral
problem~(\ref{init}) and $u_1$ a corresponding eigenvector.  Thus, having
found the solution(s) of the equation~(\ref{MainEq}) one obtains an
effective means of studying the spectral properties of the transfer
function $M_1(z,\Gamma)$ itself.  This is why the
equation~(\ref{MainEq}) and its solutions represent one of the main
subjects of the present work.

Often it turns out to be convenient to rewrite Eq.~(\ref{MainEq}) in
the form
\begin{equation}
\label{MainEqC}
X=V_1(A_1+X,\Gamma)
\end{equation}
where $X=Y-A_1$. Both equation~(\ref{MainEq}) and its
variant~(\ref{MainEqC}) will be referred to in the following as the
{\em basic equations}. Sufficient conditions for solvability of
these equations are described in the following statement.

\begin{theorem}\label{Solvability}
Let:

a{\rm)} a contour $\Gamma$ be $K_B$-bounded;

b{\rm)} the spectrum of the operator $A_1$ be strictly separated
from
the set $\sigma'(A_0)\cup\Gamma$,
\begin{equation}
\label{Separ}
d_0(\Gamma)=\dist\{\sigma(A_1),\sigma'(A_0)\cup\Gamma\}>0;
\end{equation}

c{\rm)} the inequality
\begin{equation}
\label{Best}
\cV_0(B,\Gamma)< \displaystyle\frac{1}{4}\,d_0^2(\Gamma)
\end{equation}
be valid.
Then Eq.~{\rm(\ref{MainEqC})} is uniquely solvable
in any ball
${\cal S}_1(r)\subset\bB(\cH_1,\cH_1)$ including operators
$X:\cH_1\rightarrow\cH_1$ the norms of which are bounded as
$\|X\|\leq r$ with $r$ such that
\begin{equation}
\label{Br}
  r_{\rm min}(\Gamma)\leq r < r_{\rm max}(\Gamma)
\end{equation}
where
\begin{equation}
\label{rmin}
r_{\rm min}(\Gamma)=\D\frac{d_0(\Gamma)}{2}-
\sqrt{\D\frac{d_0^2(\Gamma)}{4}
-\cV_0(B,\Gamma)} \qquad (0<r_{\rm min}(\Gamma)<d_0(\Gamma)/2)
\end{equation}
and
\begin{equation}
\label{rmax}
r_{\rm max}(\Gamma)=d_0(\Gamma)-\sqrt{\cV_0(B,\Gamma)}
\qquad (d_0(\Gamma)/2 < r_{\rm max}(\Gamma) < d_0(\Gamma)).
\end{equation}
The solution $X$ of Eq.~{\rm(\ref{MainEqC})} is the same for any
$r$ satisfying~{\rm(\ref{Br})} and in fact it belongs to the smallest
ball ${\cal S}_1(r_{\rm min})$, $\|X\|\leq r_{\rm min}(\Gamma)$.
\end{theorem}

\noindent P~r~o~o~f~.~~The proof
will be based on the Banach's Fixed Point Theorem.

Let $F(X)=V_1(A_1+X,\Gamma)$ with $X{\in}{\cal S}_1(r)$.
To begin with we search for
a condition under which the function $F$ maps
the ball ${\cal S}_1(r)$ into itself.  Since, in view of~(\ref{Br}) and
(\ref{rmax}) the condition $0<r<d_0$, $d_0=d_0(\Gamma)$
automatically holds, the spectrum of the operator $A_1+X$ does
not intersect the set $\sigma'(A_0)\cup\Gamma$ because
of condition~(\ref{Separ}).  This means that for all
$\mu\in\sigma'(A_0)\cup\Gamma$ the resolvent
$(A_1+X-\mu{I_1})^{-1}$ exists as a bounded operator in
$\cH_1$. It follows from the estimate~(\ref{V1Yest}) that
$$
\|F(X)\|\leq \cV_0(B,\Gamma)
\Sup_{\mu\in\sigma'(A_0)\cup\Gamma}\|(A_1+X-\mu{I_1})^{-1}\|.
$$
Using the identity
\begin{equation}
\label{ResIdentity}
(A_1+X-\mu{I_1})^{-1}=
\left(I_1+(A_1-\mu{I_1})^{-1}X\right)^{-1}(A_1-\mu{I_1})^{-1},
\end{equation}
one obtains the estimate
$$
\begin{array}{rcl}
\|(A_1+X-\mu{I_1})^{-1}\| &\leq& \displaystyle
\frac{1}{1-\|(A_1-\mu{I_1})^{-1}\|\cdot\|X\|}\cdot
\|(A_1-\mu{I_1})^{-1}\|\\ &\phantom{.}& \\
 &\leq& \displaystyle\frac{1}{1-\displaystyle\frac{r}{d_0}}\cdot
\frac{1}{d_0}=\displaystyle\frac{1}{d_0-r}.
\end{array}
$$
It follows from this estimate
that the ball ${\cal S}_1(r)$ is necessarily mapped
by the function $F$ into itself if the radius $r$ and the value
$\cV_0(B,\Gamma)$ are such that
\begin{equation}
\label{Est1}
\cV_0(B,\Gamma)\cdot\displaystyle\frac{1}{d_0-r}\leq r.
\end{equation}

Now, we clarify the conditions for $F$ to be a contraction.
To this end we estimate the difference
$$
F(X)-F(Y)=\Int_{\sigma'(A_0)\cup\Gamma} K_B(d\mu)\,T(\mu)
$$
where
$$
\begin{array}{rcl}
T(\mu)&=&(A_1+X-\mu{I_1})^{-1}-(A_1+Y-\mu{I_1})^{-1}\\
\phantom{\cdot}&&\\
 &=&(A_1+X-\mu{I_1})^{-1}\,(Y-X)\,(A_1+Y-\mu{I_1})^{-1}\,.
\end{array}
$$
Using again the inequality~(\ref{V1Yest}) we find
$$
\|F(X)-F(Y)\|\leq
$$
$$
\leq \cV_0(B,\Gamma)\cdot
\Sup_{\mu\in\sigma'(A_0)\cup\Gamma}\|(A_1+X-\mu{I_1})^{-1}\|\cdot
\Sup_{\mu\in\sigma'(A_0)\cup\Gamma}\|(A_1+Y-\mu{I_1})^{-1}\|
\cdot \|Y-X\|
$$
$$
\leq \cV_0(B,\Gamma)\cdot\displaystyle\frac{1}{(d_0-r)^2}\|Y-X\|.
$$
Clearly, $F$ is a contraction if
\begin{equation}
\label{Est2}
\displaystyle\frac{\cV_0(B,\Gamma)}{(d_0-r)^2}<1\,.
\end{equation}
Under the condition~(c) the inequalities~(\ref{Est1})
and~(\ref{Est2}) considered together are just equivalent
to the condition~(\ref{Br}).
Thus if the condition~(c) is valid, then
$F$ is indeed a contraction of the ball ${\cal S}_1(r)$ into
itself for any radius $r$ satisfying~(\ref{Br}).  This implies
that Eq.~(\ref{MainEqC}) has a solution in any such ball and
this solution is unique.  Consequently, the solution is the same
for all the radii satisfying~(\ref{Br}). Moreover, it belongs to the
ball ${\cal S}_1(r_{\rm min})$ with the radius $r_{\rm min}$
given by~(\ref{rmin}).

The proof is complete.%
{\nopagebreak\mbox{\phantom{MMMM}}\hfill $\Box$\par\addvspace{0.25cm}}

\begin{remark}
\label{PartA1Distances}
It should be noted that the distance
$d_0(\Gamma)=\dist\{\sigma(A_1),\sigma'(A_0)\cup\Gamma\}$ always
satisfies, of course, the inequality
$d_0(\Gamma)\leq\dist\{\sigma(A_1),{\cal E}_c(A_0)\}$ where
${\cal E}_c(A_0)$ denotes the set of the end points
$\mu_k^{(1)}$, $\mu_k^{(2)}$ of the intervals $\Delta_k^0$,
$k=1,2,\ldots,\sm$.  Thus, the condition~{\rm(\ref{Separ})} assumes
that the distance between any two of the parts
$\sigma(A_1)\cap\Delta^0_k$ of the spectrum $\sigma(A_1)$ lying
inside the different intervals $\Delta^0_k$, $k=1,2,\ldots,\sm$
is greater than $2d_0(\Gamma)$.  The same is true for the rest
part $\sigma(A_1)\setminus\sigma(A_0)$ of the spectrum of $A_1$
if it is nonempty:  \mbox{$\dist\{\sigma(A_1)\setminus\sigma(A_0),
\sigma(A_1)\cap\Delta^0_k\}>2d_0(\Gamma)$,} $k=1,2,\ldots,\sm$.
\end{remark}

\begin{theorem}\label{Hunique}
Let the conditions of {\rm Theorem~\ref{Solvability}} be valid
for a $K_B$-bounded contour $\Gamma\subset D_l$
and let $X$ be the solution of Eq.{~\rm(\ref{MainEqC})}
referred to there. Then the analogous solution
$\tilde{X}$ for any other $K_B$-bounded
contour $\tilde{\Gamma}\subset D_l$ satisfying the estimate
$\cV_0(B,\tilde{\Gamma})<\tilde{d}_0^2/4$ with
$0 < \tilde{d}_0=\dist\{\sigma(A_1),
\sigma'(A_0)\cup\tilde{\Gamma}\}\leq d_0(\Gamma)$
coincides with $X$.
\end{theorem}

\noindent P~r~o~o~f.~  The solution $X$ satisfies the inequality
$\|X\|\leq r_{\rm min}(\Gamma)$ with $r_{\rm min}(\Gamma)$ given
by~(\ref{rmin}).  This means $\|X\|<d_0(\Gamma)/2$. Similarly,
$\|\tilde{X}\|<\tilde{d}_0/2$.  The resolvent
\mbox{$(A_1-\tilde{X}-\mu)^{-1}$} is, therefore, a holomorphic
operator-valued function with its values belonging to
$\bB(\cH_1,\cH_1)$ for any $\mu\in\C$ such that
\mbox{$\dist\{\mu,\sigma(A_1)\}>\tilde{d}_0/2$.}  Recall that we
consider only the contours which result from a continuous
deformation of respective spectral intervals $\Delta^0_k$,
$k=1,2,\ldots,\sm$. The paths $\Gamma,\tilde{\Gamma}\subset D_l$
are supposed to belong to this class of contours. So that the
contour $\tilde{\Gamma}$ may be continuously transformed to the
path $\Gamma$ in such a way that for any intermediate paths
$\Gamma'$ one still has
\mbox{$\dist\{\sigma(A_1),\sigma'(A_0)\cup\Gamma'\}
\geq\tilde{d}_0$, since $d_0(\Gamma)\geq\tilde{d}_0$.}  In view of
holomorphy of the resolvent $(A_1-\tilde{X}-\mu)^{-1}$ for $\mu$
such that $\dist\{\mu,\sigma(A_1)\}>\tilde{d}_0/2$ one finds
immediately that the r.\,h. side of the equality~(\ref{MainEqC}) for
$\tilde{X}$,
$
  \tilde{X}=V_1(A_1+\tilde{X},\Gamma'),
$
remains fixed when one transforms $\Gamma'$ from $\tilde{\Gamma}$
to $\Gamma$, keeping
$\dist\{\sigma(A_1),\sigma'(A_0)\cup\Gamma'\}
\geq\tilde{d}_0$ (or even
$\dist\{\sigma(A_1),\sigma'(A_0)\cup\Gamma'\}>\tilde{d}_0/2$).
This means that $\tilde{X}$ satisfies exactly the same
equation
$
  \tilde{X}=V_1(A_1+\tilde{X},\Gamma)
$
as $X$. According to Theorem~\ref{Solvability},
the solution of Eq.~(\ref{MainEqC}) is unique and the same
in any ball ${\cal S}_1(r)\subset\bB(\cH_1,\cH_1)$
with $r$ satisfying~(\ref{Br}).
In particular, the value $r=d_0(\Gamma)/2$ can also be
substituted into~(\ref{Br}).
Meanwhile, $\tilde{X}\in{\cal S}_1(\tilde{r}_{\rm min})$
with
$
\tilde{r}_{\rm min}={\tilde{d}_0}/{2}-
\sqrt{{\tilde{d}_0^2}/{4}
-\cV_0(B,\tilde{\Gamma})}.
$
Obviously, $\tilde{r}_{\rm min}<\tilde{d}_0/2\leq d_0(\Gamma)/2$,
since $\tilde{d}_0\leq d_0(\Gamma)$
according to our assumption. Hence, $\tilde{X}$
must coincide with $X$. This completes the proof.%
{\nopagebreak\mbox{\phantom{MMMM}}\hfill $\Box$\par\addvspace{0.25cm}}

\begin{corollary}\label{HuniqueCor}
{\rm Theorem \ref{Hunique}} shows that, for a fixed multi-index $l$,
the solution of Eq.~{\rm(\ref{MainEqC})} referred to in
{\rm Theorem~\ref{Solvability}} is unique and the same
for all the $K_B$-bounded contours $\Gamma_l\subset D_l$
satisfying the inequality~{\rm(\ref{Best})}. Moreover this
solution satisfies the estimate
\begin{equation}
\label{Xr0}
\|X\|\leq r_0(B)
\end{equation}
with
\begin{equation}
\label{r0}
  r_0(B)=\Inf\limits_{\Gamma_l:\,\omega(B,\Gamma_l)>0}
  r_{\rm min}(\Gamma_l)
\end{equation}
where $r_{\rm min}(\Gamma_l)$ is given by~{\rm(\ref{rmin})} while
\begin{equation}
\label{omega}
\omega(B,\Gamma_l)=d_0^2(\Gamma_l)-4\cV_0(B,\Gamma_l).
\end{equation}
The value of $r_0(B)$ does not depend on the index $l$.
\end{corollary}

This corollary is an immediate consequence of the statement of
Theorem~{\rm\ref{Hunique}}. The only thing we want to show is
the independence of the radius $r_0(B)$ on $l$.  To prove
this, let us consider a $K_B$-bounded contour $\Gamma_l\subset
D_l$, $\Gamma_l=\mathop{\bigcup}\limits_{k=1}^\sm\Gamma_k^{l_k}$.
Denote by $\Gamma_{l'}$ a contour which is obtained from
$\Gamma_l$ by replacing a part of the curves $\Gamma_k^{l_k}$
with the conjugate ones,
\mbox{$\Gamma_k^{(-l_k)}=\{\mu:\,\overline{\mu}\in\Gamma_k^{l_k}\}$.}
Obviously, such a replacement generates, in additional to
$\Gamma_l$, \mbox{$2^\sm-1$} different contours $\Gamma_{l'}$ for
\mbox{$l'=(l'_1,l'_2,\ldots,l'_\sm)$} with $l'_k=\pm l_k$,
$k=1,2,\ldots,\sm.$ For any such contour the value of
$\cV_0(B,\Gamma_{l'})$ is the same, namely
\begin{equation}
\label{NBGprime}
    \cV_0(B,\Gamma_{l'})=\cV_0(B,\Gamma_l).
\end{equation}
Indeed,
$$
\displaystyle\int\limits_{\Gamma_k^{(-l_k)}}|d\mu|\,\|K'_B(\mu)\| =
\displaystyle\int\limits_{\Gamma_k^{l_k}}
|d\overline{\mu}|\,\|K'_B(\overline{\mu})\|,
   \quad k=1,2,\ldots,\sm.
$$
But, according to~(\ref{KBsym}),
$$
  \displaystyle\int\limits_{\Gamma_k^{l_k}}
|d\overline{\mu}|\,\|K'_B(\overline{\mu})\|=
\displaystyle\int\limits_{\Gamma_k^{l_k}}
|d{\mu}|\,\|[K'_B({\mu})]^*\|=
\displaystyle\int\limits_{\Gamma_k^{l_k}}
|d{\mu}|\,\|K'_B({\mu})\|.
$$
So that nothing happens to the value of
$
\displaystyle\int\limits_{\Gamma_l}
|d{\mu}|\,\|K'_B({\mu})\|
$
when one replaces $\Gamma_l$ by $\Gamma_{l'}$.  But this
just means that Eq.~(\ref{NBGprime}) holds true and, hence, the
infimum~(\ref{r0}) acquires the same value for any $l$.%
{\nopagebreak\mbox{\phantom{MMMM}}\hfill $\Box$\par\addvspace{0.25cm}}

\medskip

So, for a given holomorphy domain $D_l$ the solutions $X$ and
$H_1$, $H_1=A_1+X,$ do not depend on the $K_B$-bounded contours
$\Gamma_l\subset D_l$ satisfying the condition~(\ref{Best}).
But when the index $l$ changes,  $X$ and $H_1$ can also
change.  For this reason we shall supply them in the following,
when it is necessary, with the index $l$ writing,
respectively, $X^{(l)}$ and $H_1^{(l)}$, $H_1^{(l)}=A_1+X^{(l)}$.
In fact, it follows from Eq.~(\ref{NBGprime}) that if the
conditions of Theorem~\ref{Solvability} are valid for a contour
$\Gamma_l$, then they are valid for the remaining $2^{\sm-1}$ contours
$\Gamma_{l'}$ described above, too.  Therefore,
Theorem~\ref{Solvability} guarantees us, in general, the
existence of the $2^\sm$ solutions $X^{(l)}$ to the
basic equation~(\ref{MainEq}) and, hence,
the $2^\sm$ respective solutions
$H_1^{(l)}$ to the basic equation~(\ref{MainEqC}).  In the following we
shall deal only with these solutions%
\footnote{Surely, Eqs.~(\ref{MainEq}) and~(\ref{MainEqC}) are
non-linear equations and, outside the balls $\|X\|<r_{\rm
max}(\Gamma)$, they may, in principle, have other solutions,
different from the $X^{(l)}$ or $H_1^{(l)}$ the existence of
which is guaranteed by Theorem~\ref{Solvability}.}
of~(\ref{MainEq}) and~(\ref{MainEqC}).

\begin{lemma}\label{XDl}
The above solution $X^{(l)}$ of the basic
equation~{\rm(\ref{MainEqC})}, independent, for a given $D_l$, of
the contour $\Gamma_l\subset D_l$ satisfying the
condition~{\rm(\ref{Best})}, is also a solution of this equation for
any other $K_B$-bounded contour $\Gamma_l\subset D_l$ satisfying
only the condition
\begin{equation}
\label{CondMin}
\dist\{\sigma(A_1),\sigma'(A_0)\cup\Gamma_l\}>r_0(B).
\end{equation}
\end{lemma}

\noindent P~r~o~o~f~~of the lemma is reduced to an
appropriate continuous deformation of the contours, starting from a
contour satisfying~{\rm(\ref{Best})} and finishing with a
desired contour satisfying only the condition~(\ref{CondMin}).%
{\nopagebreak\mbox{\phantom{MMMM}}\hfill $\Box$\par\addvspace{0.25cm}}

Concluding the section we would like to make the following
\begin{remark}\label{MResSet}
If $\Gamma_l\subset D_l$ is a $K_B$-bounded contour satisfying
the condition~{\rm(\ref{Best})},
then the resolvent set of the transfer
function $M_1(z,\Gamma_l)$ in the domain $D(\Gamma_l)$
bounded by the curve $\Gamma_l$ and the set
$\mathop{\bigcup}_{k=1}^\sm\Delta_k^0$ is not empty.
Moreover, the curve
\begin{equation}
\label{Curve-dpopolam}
 \dist\{z,\sigma(A_1)\} = \D\frac{d_0(\Gamma_l)}{2}
\end{equation}
including its part belonging to $D(\Gamma_l)$
is entirely embedded into the resolvent set of $M_1(z,\Gamma_l)$.
In fact, there is a vicinity of the
curve~{\rm(\ref{Curve-dpopolam})} in $D(\Gamma_l)$ which is entirely
included in this set.
\end{remark}

It follows from the statement of Remark~\ref{PartA1Distances} 
that the curve~(\ref{Curve-dpopolam}) consists of $\sm$ distinct 
components surrounding respective parts 
$\sigma(A_1)\cap\Delta^0_k$ of the spectrum of $A_1$ lying 
inside the intervals $\Delta^0_k$, $k=1,2,\ldots,\sm$ and, if 
$\sigma(A_1)\setminus\sigma(A_0)\neq\emptyset$, another 
component surrounding the rest of the set $\sigma(A_1)$ lying 
outside $\sigma(A_0)$.  Every such a component is symmetric with 
respect to the real axis.

Obviously, the component of the curve~(\ref{Curve-dpopolam}) 
surrounding the set $\sigma(A_1)\setminus\sigma(A_0)$ belongs to 
the subdomain of $\C$ where $M_1(z,\Gamma_l)$ coincides with the 
initial transfer function $M_1(z)$.  Thus, at least complex 
points of this component belong to the resolvent set of 
$M_1(\cdot,\Gamma_l)$.  The component of~(\ref{Curve-dpopolam}) 
surrounding the set $\sigma(A_1)\cap\Delta^0_k$ for some 
$k=1,2,\ldots,\sm$ is entirely included in the domain 
$D(\Gamma_k^{l_k})\cup D(\Gamma_k^{(-l_k)})\cup\Delta^0_k$ where 
$D(\Gamma_k^{l_k})$ denotes domain bounded by the contour 
$\Gamma_k^{l_k}$ and interval $\Delta^0_k$ while 
$D(\Gamma_k^{(-l_k)})$ stands for domain symmetric to 
$\Gamma_k^{l_k}$ with respect to the real axis.  Since the 
function $M_1(z,\Gamma_l)$ coincides in $D(\Gamma_k^{(-l_k)})$ 
with $M_1(z)$ (see Remark~\ref{M1M1Gamma}), 
any point of (\ref{Curve-dpopolam}) lying in 
$D(\Gamma_k^{(-l_k)})$ automatically belongs to the resolvent 
set of $M_1(\cdot,\Gamma_l)$.

Further, we show that any $z\in D(\Gamma_l)$ lying 
on the curve~(\ref{Curve-dpopolam}) satisfies the inequality
\begin{equation}
\label{Estzd2}
 \dist\{z,\sigma'(A_0)\cup\Gamma_l\}\geq\D\frac{d_0(\Gamma_l)}{2}\,.
\end{equation}
We prove~(\ref{Estzd2}) by contradiction.  Let us suppose that 
there is a point $\tilde{z}\in D(\Gamma_l)$
satisfying~(\ref{Curve-dpopolam}) and such  that 
$\dist\{\tilde{z},\sigma'(A_0)\cup\Gamma_l\}<{d_0(\Gamma_l)}/{2}$.
Note that 
$\dist\{\tilde{z},\sigma'(A_0)\cup\Gamma_l\}=
\dist\{\tilde{z},\sigma'(A_0)\cup\overline{\Gamma}_l\}$ 
where as usually overlining in the notation 
$\overline{\Gamma}_l$ means closure of $\Gamma_l$
(in the case, making the closure means 
addition to $\Gamma_l$ of respective end points).
Since both sets $\sigma'(A_0)\cup\overline{\Gamma}_l$ and $\sigma(A_1)$
are closed, this implies there exists points
$z_0\in\sigma'(A_0)\cup\overline{\Gamma}_l$ and $z_1\in\sigma(A_1)$
such that
$|\tilde{z}-z_0|=
\dist\{\tilde{z},\sigma'(A_0)\cup\Gamma_l\}<{d_0(\Gamma_l)}/{2}$
and
$|\tilde{z}-z_1|=
\dist\{\tilde{z},\sigma(A_1)\}={d_0(\Gamma_l)}/{2}$.
Then it follows from the Triangle Inequality
that
$d_0(\Gamma_l)=\dist\{\sigma(A_1),\sigma'(A_0)\cup\Gamma_l\}\leq
|\tilde{z}-z_0|+|\tilde{z}-z_1|< d_0(\Gamma_l)$ and, hence,
one comes to a contradiction.
Consequently, any $z\in D(\Gamma_l)$ lying
on the curve~(\ref{Curve-dpopolam}) must satisfy~(\ref{Estzd2}).

Obviously, for $z$ satisfying~{\rm(\ref{Curve-dpopolam})} 
and~{\rm(\ref{Estzd2})} we have
$$
\|V_1(z,\Gamma_l)\| \leq \D \frac{\cV_0(B,\Gamma_l)}
{\dist\{z,\sigma'(A_0)\cup\Gamma_l\}}
$$
and
$$
\|(A_1-z)^{-1}\,V_1(z,\Gamma_l)\|\leq\D\frac{\cV_0(B,\Gamma_l)}
{\dist\{z,\sigma(A_1)\}\cdot\dist\{z,\sigma'(A_0)\cup\Gamma_l\}}
\leq\D\frac{\cV_0(B,\Gamma_l)}{d_0^2(\Gamma)/4}<1.
$$
Thus  $M_1(z,\Gamma_l)$ is invertible,
$$
M_1^{-1}(z,\Gamma_l)=\biggl(I_1+
(A_1-z)^{-1}V_1(z,\Gamma_l)\biggr)^{-1}
(A_1-z)^{-1}\,,
$$
and
$$
\|M_1^{-1}(z,\Gamma_l)\|\leq\D\frac{1}
{1-\D\frac{\cV_0(B,\Gamma_l)}{d_0^2(\Gamma)/4}}\cdot
\frac{1}{d_0(\Gamma_l)/2}\,.
$$
In the same way one can show that any real $z$ lying 
on the curve~(\ref{Curve-dpopolam}) also belongs 
to the resolvent set of $M_1(\cdot,\Gamma_l)$.  

The last statement of the remark is true due to the fact that
each regular point of $M_1(\cdot,\Gamma_l)$ is included in the
resolvent set together with a certain open neighborhood.

\section{A factorization theorem and its immediate consequences}
\label{SecFactor}

As a next step we prove the {\em factorization theorem} for
the transfer functions $M_1(z,\Gamma_l)$. This statement
will play an important role when we study the spectral
properties of the operators $H_1^{(l)}$.

\begin{theorem}\label{factorization}
Let $\Gamma_l$ be a $K_B$-bounded contour satisfying the
condition~{\rm(\ref{Best})} and $H_1^{(l)}=A_1+X^{(l)}$ with
$X^{(l)}$ the above solution of the basic
equation~{\rm(\ref{MainEq})}, $\|X^{(l)}\|\leq r_0(B)$. Then,
for $z\in\C\setminus(\sigma'(A_0)\cup\Gamma_l)$, the
transfer function $M_1(z,\Gamma_l)$ admits the factorization
\begin{equation}
\label{Mfactor}
    M_1(z,\Gamma_l)=W_1(z,\Gamma_l)\,(H_1^{(l)}-z)
\end{equation}
where $W_1(z,\Gamma_l)$ is a bounded operator in $\cH_1$,
\begin{equation}
\label{Mtild}
 W_1(z,\Gamma_l)=I_1-\D\Int_{\sigma'(A_0)\cup\Gamma_l}
 K_B(d\mu)\frac{1}{\mu-z}(H_1^{(l)}-\mu)^{-1}
\end{equation}
which is boundedly invertible for
$\dist\{z,\sigma(A_1)\}\leq{d_0(\Gamma_l)/2}$ and
\begin{equation}
\label{Mtest}
  \left\|[W_1(z,\Gamma_l)]^{-1}\right\|
 \leq \D\frac{1}{ 1-\D\frac{\cV_0(B,\Gamma_l)}
  {d_0^2(\Gamma_l)/4} }<\infty.
\end{equation}
\end{theorem}

It should be noted that the above statement recalls the known
factorization theorem by {\sc A.\,I.\,Virozub} and {\sc
V.\,I.\,Matsaev}~\cite{VirozubMatsaev} being valid for a class
of selfadjoint operator-valued functions
(see also~\cite{MarkusMatsaev}).  The results of the paper~\cite{MenShk} are
just based essentially on this theorem. However, in the case we
deal with in the present paper, the function $M_1(z,\Gamma_l)$ does
not satisfy the conditions of~\cite{VirozubMatsaev}. Moreover,
it is not even a selfadjoint operator-valued function in the
sense of~\cite{VirozubMatsaev}.

\bigskip

\noindent P~r~o~o~f~~of Theorem \ref{factorization}.  For
$z\in\C\setminus(\sigma'(A_0)\cup\Gamma_l)$,  the boundeness of
the operator $W_1(z,\Gamma_l)$ given
by~(\ref{Mtild}) is evident. To prove the
factorization~(\ref{Mfactor}) we note that
for any $z\not\in\sigma'(A_0)\cup\Gamma_l$
\begin{equation}
\label{Th2Start}
 W_1(z,\Gamma_l)\,(H_1^{(l)}-z)=H_1^{(l)}-z
-\D\Int_{\sigma'(A_0)\cup\Gamma_l}
 K_B(d\mu)\frac{1}{\mu-z}(H_1^{(l)}-\mu)^{-1}(H_1^{(l)}-z).
\end{equation}
Since
$(H_1^{(l)}-\mu)^{-1}(H_1^{(l)}-z)=I_1+(\mu-z)(H_1^{(l)}-z)^{-1}$,
one finds
$$
 \D\Int_{\sigma'(A_0)\cup\Gamma_l}
 K_B(d\mu)\frac{1}{\mu-z}(H_1^{(l)}-\mu)^{-1}(H_1^{(l)}-z) =
$$
$$
= \D\Int_{\sigma'(A_0)\cup\Gamma_l} K_B(d\mu)\frac{1}{\mu-z}
+\D\Int_{\sigma'(A_0)\cup\Gamma_l} K_B(d\mu)(H_1^{(l)}-z)^{-1}.
$$
But according to~(\ref{Mcmpl}) and (\ref{MGamma})
$$
\D\Int_{\sigma'(A_0)\cup\Gamma_l}
 K_B(d\mu)\frac{1}{\mu-z}=A_1-z-M_1(z,\Gamma_l)
$$
while according to~(\ref{MainEq})
$$
\D\Int_{\sigma'(A_0)\cup\Gamma_l}
 K_B(d\mu)(H_1^{(l)}-z)^{-1}=H_1^{(l)}-A_1.
$$
Making use of these expressions one immediately obtains
Eq.~(\ref{Mfactor}).

Further, we prove that the factor $W_1(z,\Gamma_l)$ is a
boundedly invertible operator if the condition
$\dist\{z,\sigma(A_1)\}\leq{d_0(\Gamma_l)/2}$ is valid. Under this
condition one finds
$|\mu-z|\geq\dist\{z,\sigma'(A_0)\cup\Gamma_l\}\geq
d_0(\Gamma_l)/2$, since
$\dist\{\sigma(A_1),\sigma'(A_0)\cup\Gamma_l\}=d_0(\Gamma_l)$.
On the other hand,
$H_1^{(l)}=A_1+X^{(l)}$ with $\|X^{(l)}\|<d_0(\Gamma_l)/2$
and, for $\mu\in\sigma'(A_0)\cup\Gamma_l$,
\begin{equation}
\label{HresolvEst}
\|(H_1^{(l)}-\mu)^{-1}\| <\D\frac{1}{d_0(\Gamma_l)/2}.
\end{equation}
So that
\begin{equation}
\label{forOm}
\biggl\|\, \D\Int_{\sigma'(A_0)\cup\Gamma_l}
 K_B(d\mu)\frac{1}{\mu-z}(H_1^{(l)}-\mu)^{-1}\biggr\|
< \D\frac{\cV_0(B,\Gamma_l)}{(d_0(\Gamma_l)/2)^2}<1.
\end{equation}
This means that the estimate~(\ref{Mtest}) is true.
This completes the proof.%
{\nopagebreak\mbox{\phantom{MMMM}}\hfill $\Box$\par\addvspace{0.25cm}}

\begin{corollary}\label{dmaxCor}
The statement of {\rm Theorem~\ref{factorization}}
regarding the bounded invertibility of the operator
$W_1(z,\Gamma_l)$ remains valid for $z$ in the domain
\begin{equation}
\label{Deps}
D_\varepsilon(\Gamma_l)=\left\{z:\; z\in D(\Gamma_l),
\dist\{z,\sigma(A_1)\}<\D\frac{d_{\rm max}-\varepsilon}{2}
\right\}
\end{equation}
where
\begin{equation}
\label{dmax}
d_{\rm max}=\Sup\limits_{\Gamma:\; \omega(B,\Gamma)>0}d_0(\Gamma)
\end{equation}
with $\omega(B,\Gamma)$ given by Eq.~{\rm(\ref{omega})}
while the value of $\varepsilon>0$ can be taken
arbitrarily close to zero.
\end{corollary}

Indeed, we note, first, that the value of $d_{\rm max}$ does not
depend on the index $l$ of the contour in~(\ref{dmax}), for
the same reasons that $r_0(B)$, given by~(\ref{r0}), does
not depend on $l$. Second, according to the definition of
$d_{\rm max}$, for any $\varepsilon>0$ there exists a
$K_B$-bounded contour $\Gamma_{l,\varepsilon}\subset D_l$
satisfying~(\ref{Best}) such that
$d_0(\Gamma_{l,\varepsilon})>d_{\rm max}-\varepsilon.$
Therefore, we can apply Theorem~\ref{factorization} to
$M_1(z,\Gamma_{l,\varepsilon})$. Meanwhile,
$M_1(z,\Gamma_l)=M_1(z,\Gamma_{l,\varepsilon})$ for $z\in
D(\Gamma_l)\cap D(\Gamma_{l,\varepsilon})$ (see
Sect.~\ref{Transfer_functions}) and
$W_1(z,\Gamma_l)=W_1(z,\Gamma_{l,\varepsilon})$ for such $z$,
too.  One checks the latter equality simply by deforming the
contour $\Gamma_l$ to $\Gamma_{l,\varepsilon}$ in the explicit
formula~(\ref{Mtild}).  But this just implies that
$W_1(z,\Gamma_l)$ is boundedly invertible for any $z\in
D_\varepsilon(\Gamma_l)$, since $W_1(z,\Gamma_{l,\varepsilon})$
posesses this property.%
{\nopagebreak\mbox{\phantom{MMMMM}}\hfill $\Box$\par\addvspace{0.25cm}}

\bigskip

It is easy to write some simple but useful relations
between a part of the operators $H_1^{(l)}$.  Namely, we
derive such relations between $H_1^{(l)}$ and $H_1^{(-l)}$,
$(-l)=(-l_1,-l_2,\ldots,-l_\sm)$ where $l_k$,
$k=1,2,\ldots,\sm,$ stand for the components of the multi-index
$l=(l_1,l_2,\ldots,l_\sm)$. According to our convention,
$\Gamma_{(-l)}$, $\Gamma_{(-l)}\subset D_{(-l)}$, is a contour
which is obtained from the contour $\Gamma_l$ by replacing all
the components $\Gamma_k^{l_k}$ with the conjugate ones
$\Gamma_k^{(-l_k)}$.

\begin{lemma}\label{Adjoint}
Let $\Gamma_l\subset D_l$ be a $K_B$-bounded contour
for which the conditions of {\rm Theorem~\ref{Solvability}}
are valid. Then for any
$z\in\C\setminus\biggl(\sigma'(A_0)\cup\Gamma_l\biggr)$
the following equality holds true:
\begin{equation}
\label{Hadj}
W_1(z,\Gamma_l)\,\biggl(H_1^{(l)}-z\biggr)=
\biggl(H_1^{(-l)*}-z\biggr)\,
[W_1(\overline{z},\Gamma_{(-l)})]^*\,.
\end{equation}
\end{lemma}

\noindent P~r~o~o~f.~  For $M_1(z,\Gamma_l)$ we have the
factorization formula~(\ref{Mfactor}).  The same factorization
holds as well for $M_1(\overline{z},\Gamma_{(-l)})$,
\begin{equation}
\label{Mfactors}
M_1(\overline{z},\Gamma_{(-l)})=
W_1(\overline{z},\Gamma_{(-l)})\,
\biggl(H_1^{(-l)}-\overline{z}\biggr)\,.
\end{equation}
It is easy to check that for
$\overline{z}\not\in\sigma'(A_0)\cup\Gamma_{(-l)}$
and, thus, for
$z\not\in\sigma'(A_0)\cup\Gamma_l$
\begin{equation}
\label{MsM}
\biggl[M_1(\overline{z},\Gamma_{(-l)})\biggl]^*
=M_1(z,\Gamma_l).
\end{equation}
The equality~(\ref{Hadj}), thus, follows immediately
from Eqs.~(\ref{Mfactor}),~(\ref{Mfactors})
and~(\ref{MsM}). The proof of the lemma is complete.%
{\nopagebreak\mbox{\phantom{MMMM}}\hfill $\Box$\par\addvspace{0.25cm}}

\bigskip

It is worth noting that
\begin{equation}
\label{MGadj}
\biggl[W_1(\overline{z},\Gamma_{(-l)})\biggl]^*
=I_1-\Int_{\sigma'(A_0)\cup\Gamma_l}
\biggl(H_1^{(-l)*}-\mu\biggr)^{-1}\,
\D\frac{1}{\mu-z}\, K_B(d\mu)
\end{equation}
while the $X^{(-l)*}$ determining
$H_1^{(-l)*}=A_1+X^{(-l)*}$
satisfies the equation
\begin{equation}
\label{Xadj}
X^{(-l)*}=\Int_{\sigma'(A_0)\cup\Gamma_l}
\biggl(A_1+X^{(-l)*}-\mu\biggr)^{-1}\,
K_B(d\mu)\,.
\end{equation}
One supposes here that
$\dist\{\sigma(A_1),\sigma'(A_0)\cup\Gamma_{(-l)}\}>r_0(B)$.
If, additionally, the condition~(\ref{Best}) is valid
for $\Gamma_l$, then $X^{(-l)*}$ is the only
solution of this equation. The proof of this statement repeats
literally the proof of Theorem~\ref{Solvability}.

\begin{theorem}\label{SpHalfVic}
The spectrum $\sigma(H_1^{(l)})$ of the operator
$H_1^{(l)}=A_1+X^{(l)}$ belongs to the closed $r_0(B)$-vicinity
${\cal O}_{r_0}(A_1)$ of the spectrum of $A_1$, ${\cal
O}_{r_0}(A_1)= \{z\in\C:\,\dist\{z,\sigma(A_1)\}\leq r_0(B)\}$.
If a contour $\Gamma_l\subset D_l$ satisfies~{\rm(\ref{Best})},
then the complex spectrum of $H_1^{(l)}$ belongs  to
$D_l\cap{\cal O}_{r_0}(A_1)$ while outside $D_l$ the spectrum
of $H_1^{(l)}$ is pure real.  Moreover, the spectrum
$\sigma(H_1^{(l)})$ coincides with a (subset of the) spectrum of the
transfer function $M_1(z,\Gamma_l)$. More precisely, the spectrum of
$M_1(z,\Gamma_l)$ in ${\cal O}_{d_0/2}(A_1)=\{z:\, z\in\C,
\dist\{z,\sigma(A_1)\}\leq{d_0(\Gamma_l)}/{2}\}$ is represented
only by the spectrum of $H_1^{(l)}$, i.\,e.
$\sigma\biggl(M_1(\cdot,\Gamma_l)\biggr)\cap{\cal
O}_{d_0/2}(A_1)=\sigma(H_1^{(l)})$.  Also, the following more
detailed relations hold:
\begin{eqnarray}
\label{sigp}
\sigma_p(H_1^{(l)}) & = & \sigma_p\biggl(M_1(\cdot,\Gamma_l)\biggr)
\cap{\cal O}_{d_0/2}(A_1),   \\
\label{sigc}
\sigma_c(H_1^{(l)}) &=& \sigma_c\biggl(M_1(\cdot,\Gamma_l)\biggr)
\cap{\cal O}_{d_0/2}(A_1).
\end{eqnarray}
\end{theorem}

\noindent P~r~o~o~f.~~The spectrum of $H_1^{(l)}$ belongs to
${\cal O}_{r_0}(A_1)$, since the estimate~(\ref{Xr0}) is valid
for $X^{(l)}$.  The statement regarding the spectrum of
$M_1(z,\Gamma_l)$ follows immediately from the
factorization formula~(\ref{Mfactor})
\begin{equation}
\label{invMinvH}
   [M_1(z,\Gamma_l)]^{-1}=
   (H_1^{(l)}-z)^{-1}[W_1(z,\Gamma_l)]^{-1},
\end{equation}
since $[W_1(z,\Gamma_l)]^{-1}$ exists and is bounded in ${\cal
O}_{d_0/2}(A_1)$. Since outside $\overline{D(\Gamma_l)}$ the
transfer function $M_1(z,\Gamma_l)$ coincides with the
physical-sheet transfer function $M_1(z)$ (see Remark
\ref{M1M1Gamma}), the spectrum of $M_1(\cdot,\Gamma_l)$ belongs
to $\R$ or $\overline{D(\Gamma_l)}$.  But, as we have already
established, the spectrum $\sigma(H_1^{(l)})$ represents all the
spectrum of $M_1(\cdot,\Gamma_l)$ situated in ${\cal
O}_{d_0/2}(A_1)$. Hence, the points $z\in\sigma(H_1^{(l)})$ also
belong to $\R$ or $\overline{D(\Gamma_l)}$. This just means that
for complex $z\in\sigma(H_1^{(l)})$ we have $z\in{\cal
O}_{r_0}(A_1)\cap\overline{D_l}$.

According to~(\ref{invMinvH}), not only the location of
the singularities of $[M_1(z,\Gamma_l)]^{-1}$ and
$(H_1^{(l)}-z)^{-1}$ coincide, but the properties
of these singularities are also the same and, hence,
Eqs.~(\ref{sigp}) and~(\ref{sigc}) are valid.

The proof is complete.%
{\nopagebreak\mbox{\phantom{MMMM}}\hfill $\Box$\par\addvspace{0.25cm}}

\begin{corollary}\label{HMspectr}
It follows from {\rm Corollary~\ref{dmaxCor}}
that, in fact, the complex
spectrum of the transfer function $M_1(z,\Gamma_l)$ is only represented
by the spectrum of $H_1^{(l)}$ even in a widen domain
than in the statement of
{\rm Theorem~\ref{SpHalfVic}}{\rm;} namely, in the domain
$D(\Gamma_l)\cap{\cal O}_{d_{\rm max}/2}(A_1).$
\end{corollary}

\begin{theorem}\label{HladjHl}
The spectrum of the operator $H_1^{(-l)*}$ coincides with
the spectrum of the operator
$H_1^{(l)}$. Moreover,
\begin{eqnarray}
\label{Hladjsigp}
\sigma_p(H_1^{(-l)*}) &=& \sigma_p(H_1^{(l)}) =
\sigma_p\bigl(M_1(\cdot,\Gamma_l)\bigr)
\cap{\cal O}_{d_0/2}(A_1),   \\
\label{Hladjsigc}
\sigma_c(H_1^{(-l)*}) &=& \sigma_c(H_1^{(l)}) =
\sigma_c\bigl(M_1(\cdot,\Gamma_l)\bigr)
\cap{\cal O}_{d_0/2}(A_1),
\end{eqnarray}
where $\Gamma_l$ stands for an arbitrary $K_B$-bounded contour
$\Gamma_l\subset D_l$ satisfying the condition~{\rm(\ref{Best})}.
\end{theorem}

\noindent P~r~o~o~f.~~This statement is an immediate consequence
of Lemma~\ref{Adjoint}. For a $K_B$-bounded contour
$\Gamma_l\subset D_l$ satisfying the condition~(\ref{Best})
both operators $W_1(z,\Gamma_l)$ and
$[W_1(\overline{z},\Gamma_{(-l)})]^*$ are boundedly invertible
for $z\in{\cal O}_{d_0/2}(A_1)$ and, recall,
$d_0(\Gamma_l)/2>r_0(B)$. Meanwhile, according to Eq.~(\ref{Hadj})
\begin{equation}
\label{HadjCons}
(H_1^{(l)}-z)^{-1}
[W_1(z,\Gamma_l)]^{-1}=
[W_1(\overline{z},\Gamma_{(-l)})]^{*-1}
(H_1^{(-l)*}-z)^{-1}.
\end{equation}
Therefore, the singularities of the resolvents
$(H_1^{(l)}-z)^{-1}$ and
$(H_1^{(-l)*}-z)^{-1}$
have the same location and properties
as those of $[M_1(z,\Gamma_l)]^{-1}$
in ${\cal O}_{d_0/2}(A_1)$. This assertion
implies the statement of Theorem.%
{\nopagebreak\mbox{\phantom{MMMM}}\hfill $\Box$\par\addvspace{0.25cm}}

\bigskip

\begin{theorem}\label{HlpHl2p}
Suppose that two different domains $D_{l'}$ and $D_{l''}$
include the same subdomain $D_k^{l_k}$ for some
$k=1,2,\ldots,\sm$, i.\,e., $l'_k=l''_k=l_k$.
Then the spectra of the operators $H_1^{(l')}$
and $H_1^{(l'')}$ in $D_k^{l_k}$ coincide,
\begin{equation}
\label{SpHlpHl2p}
\sigma_s(H_1^{(l')})\cap D_k^{l_k}=
\sigma_s(H_1^{(l'')})\cap D_k^{l_k}
\end{equation}
where $s=p$ or $s=c$.
\end{theorem}

\noindent P~r~o~o~f.~~The statement follows again
from Eq.~(\ref{invMinvH}) and from the identity
of $M_1(z,\Gamma_{l'})$ and $M_1(z,\Gamma_{l''})$
for
$z\in D(\Gamma_{l'})\cap D(\Gamma_{l''})\cap D_k^{l_k}$,
$\Gamma_{l'}$ and $\Gamma_{l''}$ being arbitrary
$K_B$-bounded contours satisfying~(\ref{Best}).%
{\nopagebreak\mbox{\phantom{MMMM}}\hfill $\Box$\par\addvspace{0.25cm}}

\begin{corollary}\label{SpTheSame}
{\rm(Symmetry of the resonance spectrum with respect to the real axis)}
The spectra of any two operators $H^{(l')}$ and $H^{(l'')}$ for
$l'=(l'_1,l'_2,\ldots,l'_\sm)$ and
$l''=(l''_1,l''_2,\ldots,l''_\sm)$ are related to each other as
follows
\begin{eqnarray*}
\sigma_s(H_1^{(l'')})\cap D_k^{l''_k} &=&
\sigma_s(H_1^{(l')})\cap D_k^{l'_k}
\quad\mbox{if}\quad l''_k=l'_k, \\
\sigma_s(H_1^{(l'')})\cap D_k^{l''_k} &=&
\sigma^*_s(H_1^{(l')})\cap D_k^{l''_k}
\quad\mbox{if}\quad l''_k=-l'_k
\end{eqnarray*}
where the symbol {\rm``$*$''} denotes complex conjugation
and, as previously,  $s=p$ or $s=c$.
\end{corollary}

In the following we shall use the operators
\begin{equation}
\label{Omega}
\Omega^{(l)}=\Int_{\sigma'(A_0)\cup\Gamma_l}
(H_1^{(-l)*}-\mu)^{-1}K_B(d\mu)\,
(H_1^{(l)}-\mu)^{-1}
\end{equation}
acting in $\cH_1$, where $\Gamma_l$ stands for a $K_B$-bounded
contour satisfying the condition~(\ref{Best}). The operator
$\Omega^{(l)}$ does not depend (for a fixed $l$) on the choice
of such a $\Gamma_l$.  In the same way as we came to the
estimate~(\ref{forOm}) one can obtain the following estimate for
$\Omega^{(l)}$ (see Lemma~\ref{LIntXBEBY} of
Appendix~\ref{IntOpMer}):
\begin{equation}
\label{Omest}
\|\Omega^{(l)}\|<\D\frac{\cV_0(B,\Gamma_l)}{(d_0(\Gamma_l)/2)^2}<1.
\end{equation}
Obviously, we have the equality
\begin{equation}
\label{Omadj}
\Omega^{(l)*}=\Omega^{(-l)}.
\end{equation}

It should be noted that in the case where the spectra of the
entries $A_0$ and $A_1$ do not overlap, the
operators~(\ref{Omega}) as well as the operators $H_1^{(l)}$ do
not depend on $l$.  In this case one has $H_1=A_1+B_{01}Q_{01}$
with the contraction
$Q_{01}=\Int_{\sigma(A_0)}E_0(d\mu)B_{01}(H_1-\mu)^{-1}$,
$Q_{01}:\, \cH_1\rightarrow\cH_0$ (see Theorem~5 of
Ref.~\cite{MotRem}; cf.~\cite{AdL,AdLMSr,MenShk}),  and
$\Omega=Q_{01}^*Q_{01}$. Changing the inner product in $\cH_1$
to $[\cdot,\cdot]=\lal(I_1+\Omega)\cdot,\cdot\ral$ turns
$H_1$ into a self-adjoint operator.  However, in the case we
consider in the present work this is in general not true since
$H_1^{(l)}$ can have a complex spectrum.

\begin{theorem}
\label{MHOmega}
The operators $\Omega^{(l)}$ possess the following properties%
\footnote{For the case where $A_1=\lambda I_1$ with $\lambda\in\R$
and $\mathop{\rm dim}\cH_1<\infty$ one can find formulas
similar to those in Eqs.~(\ref{MOmega}) and~(\ref{HOmega})
in the final part of Sect.~4 of Ref.~\cite{Howland}.
See also~\cite{MarkusMatsaev,MenShk,VirozubMatsaev}.}:
\begin{equation}
\label{MOmega}
-\D\frac{1}{2\pi\ri}\Int_\gamma dz\,[M_1(z,\Gamma_l)]^{-1}=
(I_1+\Omega^{(l)})^{-1}
\end{equation}
and
\begin{equation}
\label{HOmega}
-\D\frac{1}{2\pi\ri}\Int_\gamma dz\,z\,[M_1(z,\Gamma_l)]^{-1}=
(I_1+\Omega^{(l)})^{-1}H_1^{(-l)*}=
H_1^{(l)}(I_1+\Omega^{(l)})^{-1}
\end{equation}
where $\gamma$ stands for an arbitrary rectifiable closed
{\rm(}including the points at infinity if the entry $A_1$ is
unbounded{\rm)} contour going in the positive direction around the
spectrum of $H_1^{(l)}$ inside the set ${\cal
O}_{d_0(\Gamma)/2}(A_1)$.  The integration over $\gamma$ is
understood in the strong sense.
\end{theorem}
\noindent P~r~o~o~f~.~~First, we note that, using the
Closed Graph Theorem and the definition~(\ref{DefXBEBYint})
of the integral~(\ref{Omega}), one can easily check
that for any $u_1\in\cD(H_1^{(-l)*})=\cD(A_1)$
the image $\Omega^{(l)}u_1$
belongs to $\cD(H_1^{(l)})=\cD(A_1)$.
And due to~(\ref{Omest}) the operator $I_1+\Omega^{(l)}$
is a bijection of $\cD(A_1)$ on $\cD(A_1)$.

Further, we prove the validity of Eq.~(\ref{MOmega}).
At the beginning we recall that if
$z\in{\cal O}_{d_0(\Gamma_l)/2}(A_1)$, then
the factorization~(\ref{invMinvH}) and~(\ref{HadjCons})
holds for $[M_1(z,\Gamma_l)]^{-1}$
with the holomorphic functions $[W_1(z,\Gamma_l)]^{-1}$ and
$[W_1(\overline{z},\Gamma_{(-l)})]^{*-1}$
taking their values in $\bB(\cH_1,\cH_1)$.
Meanwhile the product $\Omega^{(l)}(H_1^{(l)}-z)^{-1}$
can be written as
$$
\Omega^{(l)}(H_1^{(l)}-z)^{-1}=F_1(z)+F_2(z)
$$
with
\begin{equation}
\label{F1int}
F_1(z)=\Int_{\sigma'(A_0)\cup\Gamma_l}
\D\frac{ (H_1^{(-l)*}-\mu)^{-1}\,K_B(d\mu)\,(H_1^{(l)}-\mu)^{-1}}
{\mu-z}
\end{equation}
and
$$
   F_2(z)=\biggl([W_1(\overline{z},\Gamma_{(-l)})]^*-I_1\biggr)
   (H_1^{(l)}-z)^{-1},
$$
since
$$
    (H_1^{(l)}-\mu)^{-1}(H_1^{(l)}-z)^{-1}=
    [(H_1^{(l)}-\mu)^{-1}-(H_1^{(l)}-z)^{-1}](\mu-z)^{-1}
$$
and since $[W_1(\overline{z},\Gamma_{(-l)})]^*$ is given
by~(\ref{MGadj}). Therefore,
\begin{eqnarray*}
\lefteqn{(I_1+\Omega^{(l)})\,[M_1(z,\Gamma_l)]^{-1} =}\\
&=&[M_1(z,\Gamma_l)]^{-1}+F_1(z)[W_1(z,\Gamma_l)]^{-1}
+([W_1(\overline{z},\Gamma_{(-l)})]^*-I_1)[M_1(z,\Gamma_l)]^{-1}\\
&=& F_1(z)\,[W_1(z,\Gamma_l)]^{-1}+(H_1^{(-l)*}-z)^{-1}.
\end{eqnarray*}
One should notice that the function $F_1(z)$ is holomorphic
in $z$ inside the contour $\gamma$,
$\gamma\subset{\cal O}_{d_0(\Gamma_l)/2}(A_1)$,
since the argument $\mu$ of the integrand in~(\ref{F1int})
belongs to $\sigma'(A_0)\cup\Gamma_l$
and, thereby, always \mbox{$|z-\mu|\geq d_0(\Gamma_l)/2>0$.}
Thus the term $F_1(z)[W_1(z,\Gamma_l)]^{-1}$
makes no contribution to the integral
$$
-\D\frac{1}{2\pi\ri}\Int_\gamma dz (I_1+\Omega^{(l)})
[M_1(z,\Gamma_l)]^{-1}
$$
while the resolvent $(H_1^{(-l)*}-z)^{-1}$
gives the identity $I_1$. Therefore, we have proved that
Eq.~(\ref{MOmega}) is indeed valid.

Regarding Eq.~(\ref{HOmega}) one finds
\begin{eqnarray*}
   \lefteqn{-\D\frac{1}{2\pi\ri}\Int_\gamma dz (I_1+\Omega^{(l)})
   \,z\, [M_1(z,\Gamma_l)]^{-1} =}\\
  &=&-\D\frac{1}{2\pi\ri}\Int_\gamma dz
  \,z\, F_1(z)\,[W_1(z,\Gamma_l)]^{-1}
  -\D\frac{1}{2\pi\ri}\Int_\gamma dz \,z\,(H_1^{(-l)*}-z)^{-1}
\end{eqnarray*}
where only the last integral is non-zero, giving a contribution
just equal to $H_1^{(-l)*}$. The second equation of~(\ref{HOmega})
can be checked in the same way. The proof is complete.%
{\nopagebreak\mbox{\phantom{MMMM}}\hfill $\Box$\par\addvspace{0.25cm}}

\begin{remark}\label{HlHml}
The formula~{\rm(\ref{HOmega})} implies
$$
H_1^{(l)*}=
(I_1+\Omega^{(-l)})\,H_1^{(-l)}\,(I_1+\Omega^{(-l)})^{-1}\,.
$$
\end{remark}

\begin{theorem}
\label{MPOmega}
Let $\lambda$ be an isolated eigenvalue of the operator
$H_1^{(l)}$ and, consequently, of the operator $H_1^{(-l)*}$ and
of the transfer function $M_1(z,\Gamma_l)$ taken for a
$K_B$-bounded contour $\Gamma_l$ satisfying the
condition~{\rm(\ref{Best})}. Denote by $\sP_\lambda^{(l)}$ and
$\sP_\lambda^{(-l)*}$ the respective eigenprojections of the
operators $H_1^{(l)}$ and $H_1^{(-l)*}$ and by
$P_\lambda^{(l)}$, the residue of $M_1(z,\Gamma_l)$ at
$z=\lambda$,
\begin{equation}
\label{Plambda}
\sP_\lambda^{(l)}=
-\D\frac{1}{2\pi\ri}\Int_\gamma dz\,\,(H_1^{(l)}-z)^{-1},
\end{equation}
$$
\sP_\lambda^{(-l)*}=
-\D\frac{1}{2\pi\ri}\Int_\gamma dz\,\,(H_1^{(-l)*}-z)^{-1}
$$
and
$$
P_\lambda^{(l)}=
-\D\frac{1}{2\pi\ri}\Int_\gamma dz\,\,[M_1(z,\Gamma_l)]^{-1}
$$
where $\gamma$ stands for an arbitrary rectifiable closed
contour situated in a sufficiently close vicinity
of the point $\lambda$ and going in the positive direction
around $\lambda$ so that $\gamma\cap\Gamma_l=\emptyset$
and no points of the spectrum of $M_1(\cdot,\Gamma_l)$, except the
eigenvalue $\lambda$, lie inside $\gamma$. Then the following
relations are valid:
\begin{equation}
\label{MresiduePP}
P_\lambda^{(l)}=(I_1+\Omega^{(l)})^{-1}\,\,\sP_\lambda^{(-l)*}=
\sP_\lambda^{(l)}\,\,(I_1+\Omega^{(l)})^{-1}\,.
\end{equation}
\end{theorem}
\noindent P~r~o~o~f~~is carried out in the same way as the proof
of the relation~(\ref{MOmega}) in Theorem~\ref{MHOmega},
only the path of integration is changed.%
{\nopagebreak\mbox{\phantom{MMMM}}\hfill $\Box$\par\addvspace{0.25cm}}

\section{Some properties of real eigenvalues}
\label{RealEigen}

If $\lambda$ is a real eigenvalue of $H_1^{(l)}$,
then it can not belong to the spectrum $\sigma'(A_0)$ of the
entry $A_0$ lying outside
$\displaystyle\mathop{\cup}_{k=1}^\sm\Delta_k^0$.  Indeed,
according to Theorem~\ref{SpHalfVic}, the spectrum of $H_1^{(l)}$
for arbitrary $l$ is situated in the $r_0(B)$-vicinity ${\cal
O}_{r_0}(A_1)$ of the set $\sigma(A_1)$ and in any case
$r_0(B)<\frac{1}{2}\dist\{\sigma'(A_0),\sigma(A_1)\}$  so that
automatically
\begin{equation}
\label{noninters}
\sigma'(A_0)\cap\sigma(H_1^{(l)})=\emptyset 
\qquad\mbox{and in particular} \qquad
\sigma'(A_0)\cap\sigma_p(H_1^{(l)})=\emptyset\,.
\end{equation}
Hence, such a $\lambda$ belongs either to the resolvent set
$\rho(A_0)$ of the entry $A_0$ or it is embedded into
the continuous spectrum of $A_0$ in
$\mathop{\bigcup}\limits_{k=1}^\sm\Delta^0_k$.

\begin{lemma}\label{LReal2}
If a vector $\psi^{(1)}\in\cD(A_1)$ is an eigenvector of
$H_1^{(l)}$ corresponding to a real eigenvalue
$\lambda\in\rho(A_0)$, $H_1^{(l)}\psi^{(1)}=\lambda\psi^{(1)}$,
then the vector $\Psi=(\psi^{(0)},\psi^{(1)})\in\cH$ with
\begin{equation}
\label{psi0}
\psi^{(0)}=-R_0(\lambda)B_{01}\psi^{(1)}
\end{equation}
is an eigenvector of $\bH$, $\bH\Psi=\lambda\Psi$. The converse
statement is also true: if $\lambda$, $\lambda\in\rho(A_0)$,
is a real eigenvalue of $H_1^{(l)}$ and
$\bH\Psi=\lambda\Psi$ for some $\Psi=(\psi^{(0)},\psi^{(1)})$
with $\psi^{(0)}\in\cD(A_0)$ and $\psi^{(1)}\in\cD(A_1)$,
then $\psi^{(0)}$ is related to $\psi^{(1)}$ as in~{\rm(\ref{psi0})}
{\rm(}and, therefore, $\psi^{(1)}$ can not be zero{\rm)}
and $H_1^{(l)}\psi^{(1)}=\lambda\psi^{(1)}$.
\end{lemma}

\noindent P~r~o~o~f~.~~Let us consider a $K_B$-bounded contour
$\Gamma_l\subset D_l$ satisfying the estimate~(\ref{Best}).
Since $\lambda\in\rho(A_0)\cap\R$, and therefore
$\lambda\not\in\overline{D(\Gamma_l)}$, we have
\begin{equation}
\label{MGM}
M_1(\lambda,\Gamma_l)=M_1(\lambda).
\end{equation}
So that, according to the factorization formula~(\ref{Mfactor}),
the eigenvector $\psi^{(1)}$ of $H_1^{(l)}$ is automatically
an eigenvector of the initial transfer function $M_1(\cdot)$,
$$
  (A_1-\lambda-B_{10}(A_0-\lambda)^{-1}B_{01})\psi^{(1)}=0.
$$
Introducing $\psi^{(0)}$ via~(\ref{psi0}) one immediately finds
that $\Psi=(\psi^{(0)},\psi^{(1)})$ turns out to be an
eigenvector for $\bH$, $\bH\Psi=\lambda\Psi$.

For the converse statement one first observes that if
$\bH\Psi=\lambda\Psi$ with $\Psi=(\psi^{(0)},\psi^{(1)})$ (and,
thus, $A_0\psi^{(0)}+B_{01}\psi^{(1)}=\lambda\psi^{(0)}$), then
Eq.~(\ref{psi0}) holds true. But this means that
$M_1(\lambda)\psi^{(1)}=0$ and, hence,
$M_1(\lambda,\Gamma_l)\psi^{(1)}=0$ is also true.  Then, due to
Eq.~(\ref{Mfactor}) and invertibility of $W_1(\lambda,\Gamma_l)$
(see Theorem\,\ref{factorization}),
$\psi^{(1)}$ is an eigenvector of
$H_1^{(l)}$, $H_1^{(l)}\psi^{(1)}=\lambda\psi^{(1)}$.%
{\nopagebreak\mbox{\phantom{MMMM}}\hfill $\Box$\par\addvspace{0.25cm}}

\bigskip

If an eigenvalue $\lambda$ of $H_1^{(l)}$ belongs to
$\Delta_k^0=(\mu_k^{(1)},\mu_k^{(2)})$ for some
$k=1,2,\ldots,\sm$, then
$$
|\lambda-\mu_k^{(i)}|\geq\dist\{\mu_k^{(i)},\sigma(A_1)\}-r_0(B),
\quad i=1,2,
$$
and, therefore, the $\lambda$ is situated in this case
strictly inside the interval $\Delta_k^0$.
Recall that according to our assumption the entry $A_0$
has no point spectrum inside $\Delta_k^0$. Since
$\Delta_k^0$ is a part of the continuous spectrum of $A_0$,
the resolvent $R_0(z)=(A_0-z)^{-1}$ for $z=\lambda\pm\ri0$
exists being however an unbounded operator. Nevertheless
a statement analogous to Lemma~\ref{LReal2} is valid in this
case, too.

\begin{lemma}\label{LReal3}
If a vector $\psi^{(1)}\in\cD(A_1)$ is an eigenvector of
$H_1^{(l)}$ corresponding to a real eigenvalue
$\lambda\in\Delta^0_k=(\mu^{(1)}_k,\mu^{(2)}_k)$,
$k=1,2,\ldots,\sm$, \, $H_1^{(l)}\psi^{(1)}=\lambda\psi^{(1)}$,
then either

a{\rm)} $E^0(\mu)B_{01}\psi^{(1)}=0$ for all $\mu\leq\mu_k^{(2)}$

\noindent or

b{\rm)} $E^0(\mu)B_{01}\psi^{(1)}\neq0$ for any $\mu\in\Delta^0_k$,

c{\rm)} the function $\|E^0(\mu)B_{01}\psi^{(1)}\|$ is
differentiable in $\mu$ on $\Delta^0_k$

\noindent and

d{\rm)} $\reduction{\D\frac{d}{d\mu}\|E^0(\mu)B_{01}\psi^{(1)}\|}
         {\mu=\lambda}=0.$

In both cases the vector $\psi^{(0)}$ given by~{\rm(\ref{psi0})}
exists in $\cD(A_0)$ and $\Psi=(\psi^{(0)},\psi^{(1)})$
is an eigenvector of $\bH$, $\bH\Psi=\lambda\Psi.$

The converse statement is also true. Namely, if a 
$\Psi=(\psi^{(0)},\psi^{(1)})$ with $\psi^{(0)}\in\cD(A_0)$ and 
$\psi^{(1)}\in\cD(A_1)$ is an eigenvector of $\bH$, 
$\bH\Psi=\lambda\Psi$, corresponding to an eigenvalue 
$\lambda$ of $H_1^{(l)}$, $\lambda\in\Delta^0_k$, then either 
the condition {\rm(}a{\rm)} is valid or the conditions 
{\rm(}b\,--\,d{\rm)} are valid.  In both cases the 
relation~{\rm(\ref{psi0})} is retained meaning, in particular, 
that $\psi^{(1)}\neq0$ and $\psi^{(1)}$ is an eigenvector of 
$H_1^{(l)}$ corresponding to the eigenvalue $\lambda$.
\end{lemma}

\noindent P~r~o~o~f~.~~We first prove the direct statement.
To this end we consider the equality
$$
\lal H_1^{(l)}\psi^{(1)},
\psi^{(1)}\ral=\lambda
\|\psi^{(1)}\|^2, \qquad \lambda\in\Delta_k^0,
$$
which becomes, according to Eqs.~(\ref{MainEq}) and (\ref{Return}),
\begin{equation}
\label{ZerInit}
\lal(A_1-\lambda)\psi^{(1)},
\psi^{(1)}\ral +
\Int_{\sigma'(A_0)\cup\Gamma_l}
\D\frac{\lal K_B(d\mu)\psi^{(1)},
\psi^{(1)}\ral}
{\lambda-\mu}=0.
\end{equation}
Since the denominator of the integrand is non-zero for
$\mu\in\sigma(A_0)\setminus\Delta_k^0$, we can deform
the part $\Gamma_l\setminus\Gamma_k^{l_k}$ of the
contour $\Gamma_l$ in~(\ref{ZerInit}) back into
the intervals $\Delta_i^0$, $i=1,2,\ldots,\sm,$ $i\neq k$.
As a result, Eq.~(\ref{ZerInit}) acquires the form
\begin{equation}
\label{ZerInit1}
\begin{array}{rcl}
\lal(A_1-\lambda)\psi^{(1)},
\psi^{(1)}\ral
&+&\hskip-0.5em\Int_{\sigma(A_0)\setminus\Delta^0_k}
\hskip-0.2em
\D\frac{ \lal B_{10}E_0(d\mu)B_{01}\psi^{(1)},
\psi^{(1)}\ral }{\lambda-\mu} \\
&+&\Int_{\Gamma_k^{l_k}}d\mu\,
\D\frac{ \lal K'_B(\mu)\psi^{(1)},
\psi^{(1)}\ral }{\lambda-\mu}=0.
\end{array}
\end{equation}
Obviously, the first and second terms are real, and an imaginary
component may appear in the l.\,h. side of
Eq.~(\ref{ZerInit1}) only in the third term. To find this
component one can simply transform the integration path to the
two intervals $[\mu_k^{(1)},\lambda-\varepsilon]$ and
$[\lambda+\varepsilon,\mu_k^{(2)}]$ and the semicircle
$|\mu-\lambda|=\varepsilon$, \,\, $l_k\cdot\Img\mu\geq0$,
between them.  Then taking the limit $\varepsilon\downarrow0$
one obtains
$$
\Img\Int_{\Gamma_k^{l_k}}d\mu\,
\D\frac{\lal K'_B(\mu)\psi^{(1)},
\psi^{(1)}\ral}
{\lambda-\mu}
=l_k\cdot\pi\ri\lal K'_B(\lambda)\psi^{(1)},
\psi^{(1)}\ral=0.
$$
Therefore, we have
$$
\lal K'_B(\lambda)\psi^{(1)},
\psi^{(1)}\ral=0.
$$
Meanwhile, for any $u_1\in\cH_1$ and $\mu\in\Delta^0_k$
$$
\lal K'_B(\mu)u_1,u_1\ral=\D\frac{d}{d\mu}
\lal B_{10}E^0(\mu)B_{01}u_1,u_1\ral
=\D\frac{d}{d\mu}
\|E^0(\mu)B_{01}u_1\|^2.
$$
Thus, the condition
$\lal K'_B(\lambda)u_1,u_1\ral=0$
implies that
either $\|E^0(\lambda)B_{01}u_1\|=0$ or,
if \mbox{$\|E^0(\lambda)B_{01}u_1\|\neq0$,} then
the function \mbox{$\|E^0(\mu)B_{01}u_1\|$} is differentiable
at \mbox{$\mu=\lambda$} and
\mbox{$
\reduction{\D\frac{d}{d\mu}\|E^0(\mu)B_{01}u_1\|}
{\mu=\lambda}=0.
$}

Since $\|E^0(\mu)B_{01}u_1\|$ is a non-decreasing function
of the variable $\mu$, in the first case we have to conclude that
$\|E^0(\mu)B_{01}u_1\|=0$ for all $\mu\leq\lambda$ and, hence,
$\lal K_B(\mu)u_1,u_1\ral=\|E^0(\mu)B_{01}u_1\|^2=0$ for
$\mu\leq\lambda$, too.  Since $\lal K_B(\mu)u_1,u_1\ral$ is
supposed to be a holomorphic function of $\mu\in D_k^{l_k}$ we
find $\lal K_B(\mu)u_1,u_1\ral\equiv0$ for $\mu\in D_k^{l_k}$
and, consequently, $\|E^0(\mu)B_{01}u_1\|^2=\lal
K_B(\mu)u_1,u_1\ral\equiv0$ for $\lambda<\mu\leq\mu_k^{(2)}$,
too. So that we come to the condition (a).  Applying this
condition to $u_1=\psi^{(1)}$ we find that in this case the
formula~(\ref{psi0}) makes sense and
$\Psi=(\psi^{(0)},\psi^{(1)})$ with
$
\psi^{(0)}=
-R_0(\lambda\pm\ri0)B_{01}\psi^{(1)}
\in{\cal D}(A_0)
$
is an eigenvector for $\bH$.

In the second case the (non-decreasing) function
$\|E^0(\mu)B_{01}\psi^{(1)}\|$ is non-zero (condition (b)\,),
and differentiable at any $\mu\in\Delta^0_k$ (condition (c)\,),
and
$
\reduction{\D\frac{d}{d\mu}\|E^0(\mu)B_{01}\psi^{(1)}\|}
{\mu=\lambda}=0
$
(condition (d)\,).  So that for any finite $\varepsilon,\eta>0$
such that $[\lambda-\varepsilon,\lambda+\eta]\subset\Delta^0_k$
we have the estimate
$$
\biggl\|\Int_{\lambda-\varepsilon}^{\lambda+\eta}
\D\frac{E_0(d\mu)B_{01}\psi^{(1)}}{\lambda-\mu}\biggr\|
\leq\Int_{\lambda-\varepsilon}^{\lambda+\eta}d\mu\,
\D\frac{\phantom{m}
\D\frac{d}{d\mu}\|E_0(d\mu)B_{01}\psi^{(1)}\|
\phantom{m}}
{|\lambda-\mu|\phantom{m}}
\leq C(\varepsilon,\eta)
$$
with some positive $C(\varepsilon,\eta)<\infty$.
Consequently, the integral
$$
 -R_0(\lambda\pm\ri0)B_{01}\psi^{(1)}=\Int_{\sigma(A_0)}
\D\frac{E_0(d\mu)B_{01}\psi^{(1)}}{\lambda-\mu}
$$
exists and determines an element $\psi^{(0)}\in{\cal
D}(A_0)$ such that again for
$\Psi=(\psi^{(0)},\psi^{(1)})$
one finds $\bH\Psi=\lambda\Psi$.

Let us now prove the converse statement. First, we note that
if $\bH\Psi=\lambda\Psi$ with $\Psi=(\psi^{(0)},\psi^{(1)})$,
$\psi^{(0)}\in\cD(A_0)$, $\psi^{(1)}\in\cD(A_1)$, then
\begin{equation}
\label{HfirstEq}
(A_0-\lambda)\psi^{(0)}=-B_{01}\psi^{(1)}.
\end{equation}
Let $E^0_{ac}(\mu)$ be the spectral function corresponding to
the absolutely continuous spectrum of $A_0$,
$E^0_{ac}(\mu)=E_0^{ac}\biggl((-\infty,\mu)\biggr)$.  Applying
the projection $E_0^{ac}(\delta)$ with $\delta\subset\Delta_k^0$
to both parts of Eq.~(\ref{HfirstEq})  we obtain
$$
   \Int_\delta (\mu-\lambda)dE^0_{ac}(\mu)\psi^{(0)}=
-\Int_\delta dE^0(\mu)B_{01}\psi^{(1)}
$$
(recall that we assume $E_0(\delta)B_{01}=E_0^{ac}(\delta)B_{01}$
for any Borel set $\delta\subset\mathop{\bigcup}_{k=1}^\sm\Delta_k^0$).
Thus, for the norm squares, one finds
$$
   \Int_\delta (\mu-\lambda)^2
   d\lal E^0_{ac}(\mu)\psi^{(0)},\psi^{(0)}\ral=
\Int_\delta d\mu \lal K'_B\psi^{(1)},\psi^{(1)}\ral,
$$
for an arbitrary interval $\delta\subset\Delta_k^0$.
Since the function $\lal E^0_{ac}(\mu)\psi^{(0)},\psi^{(0)}\ral$
is absolutely continuous and, hence, almost everywhere
differentiable, one further finds
$$
\lal K'_B(\mu)\psi^{(1)},\psi^{(1)}\ral=(\mu-\lambda)^2
\D\frac{d}{d\mu}\lal E^0_{ac}(\mu)\psi^{(0)},\psi^{(0)}\ral
$$
for almost all $\mu\in\Delta^0_k$. Meanwhile, the derivative
$\D\frac{d}{d\mu}\lal E^0_{ac}(\mu)\psi^{(0)},\psi^{(0)}\ral$ is
an element of $L_1\biggl(\sigma_{ac}(A_0)\biggr)$.  That is, the
function
$\lal{K'_B(\mu)}\psi^{(1)},\psi^{(1)}\ral\cdot(\mu-\lambda)^{-2}$
must also be integrable over any interval
$\delta\subset\Delta^0_k$.  Surely, this is possible only if
$\lal K'_B(\lambda)\psi^{(1)},\psi^{(1)}\ral=0$.  Now one need
only  repeat the respective consideration from the proof of the
direct part of the lemma and, as a result,  come to the conditions
(a) or (b\,--\,d). With these conditions the formula~(\ref{psi0}) is
again correct. The only thing which must be stressed in the
case of the condition (a) is the fact that
$B_{01}\psi^{(1)}\neq0$ if $\psi^{(0)}\neq0.$ Indeed, according
to Eq.~(\ref{HfirstEq}) the assumption $B_{01}\psi^{(1)}=0$
implies that $\lambda\in\sigma_p(A_0)$.  But this contradicts
our initial assumption regarding continuity of the
spectrum of $A_0$ within the intervals $\Delta^0_k$. 
Thus, $\psi^{(1)}$ can not be zero and $M_1(\lambda\pm\ri0)\psi^{(1)}=0$. 
This also means that $M_1(\lambda,\Gamma_l)\psi^{(1)}=0$ for 
arbitrary $K_B$-bounded contour $\Gamma_l\subset D_l$ satisfying the 
condition~(\ref{Best}). Then applying Theorem \ref{factorization} 
we conclude that $H_1^{(l)}\psi^{(1)}=\lambda\psi^{(1)}$. 
The proof is complete.%
{\nopagebreak\mbox{\phantom{MMMM}}\hfill $\Box$\par\addvspace{0.25cm}}

\begin{remark}\label{Remark1}
In the above proof we have also
found that if an eigenvalue $\lambda$
of $H_1^{(l)}$ is embedded into an interval $\Delta_k^0$,
$k=1,2,\ldots,\sm,$ then the function $\lal
K'_B(\mu)\psi^{(1)},\psi^{(1)}\ral\cdot(\mu-\lambda)^{-2}$ with
$\psi^{(1)}$ an eigenvector of $H_1^{(l)}$
corresponding to the $\lambda$ is
integrable over every interval $\delta\subset\Delta_k^0$.
In fact, this means that not only
$\lal{K'_B}(\lambda)\psi^{(1)},\psi^{(1)}\ral=0$ but also
\begin{equation}
\label{Der2zero}
\reduction{\D\frac{d}{d\mu}
\lal K'_B(\mu)\psi^{(1)},\psi^{(1)}\ral}{\mu=\lambda}=0.
\end{equation}
\end{remark}
The latter statement follows from the fact that the function $\lal
K'_B(\lambda)\psi^{(1)},\psi^{(1)}\ral$ is holomorphic with
respect to the variable $\mu$ in any vicinity of the point
$\lambda$ included in $\Delta^0_k\cup D_k^-\cup D_k^+$ (and,
moreover, this function is non-negative for $\mu\in\Delta^0_k$).%
{\nopagebreak\mbox{\phantom{MMMM}}\hfill $\Box$\par\addvspace{0.25cm}}

\begin{corollary}
\label{SigPointGeneral}
The statements of {\rm Lemmas \ref{LReal2}} and 
{\rm\ref{LReal3}} imply $\sigma_p(H_1^{(l)})\subset\sigma_p({\bf 
H})$.  Also, it immediately follows from these lemmas that any 
eigenvector $\psi^{(1)}$ corresponding to an eigenvalue 
$\lambda\in\sigma_p(H_1^{(l)})\cap\R$ of the operator 
$H_1^{(l)}=A+X^{(l)}$ for a certain $l=(l_1,l_2,\ldots,l_\sm)$ 
is such an eigenvector, 
$H_1^{(l')}\psi^{(1)}=\lambda\psi^{(1)}$, for the remaining 
$2^{\sm-1}$ operators $H_1^{(l')}=A_1+X^{(l')}$ for 
$l'=(l'_1,l'_2,\ldots,l'_\sm)$ with arbitrary $l'_k=\pm1$, 
$k=1,2,\ldots,\sm$.  Thus, the set $\sigma_p(H_1^{(l)})\cap\R$ 
is the same for all the $2^\sm$ operators $H_1^{(l)}$. 
\end{corollary}

\begin{lemma}
\label{labda-real}
If some $\lambda$, $\lambda\in\R$ is an isolated eigenvalue
of the operator $H_1^{(l')}=A_1+X^{(l')}$ for some
$l'=(l'_1,l'_2,\ldots,l'_\sm)$, then this $\lambda$
is also such an eigenvalue for the remaining
$2^{m-1}$ operators $H_1^{(l)}=A_1+X^{(l)}$
for $l=(l_1,l_2,\ldots,l_\sm)$ with arbitrary $l_k=\pm1$,
$k=1,2,\ldots,\sm.$ Moreover, the resolvents
for all the $2^\sm$ operators $H_1^{(l)}$ have
a pole of the first order at $z=\lambda$ allowing the decomposition
\begin{equation}
\label{Rlpole}
(H_1^{(l)}-z)^{-1}=\D\frac{{\sf P}_\lambda^{(l)}}{\lambda-z}+
\tilde{R}_\lambda^{(l)}(z)
\end{equation}
with $\tilde{R}_\lambda^{(l)}(z)$ holomorphic in a
vicinity of $\lambda$. 
Also, the factorization~{\rm(\ref{MresiduePP})} holds where
$P^{(l)}_\lambda$ does not depend on $l$, since it is the residue of
the initial inverse transfer function $R_{11}(z)=[M_1(z)]^{-1}$
at $z=\lambda$,
\begin{equation}
\label{ResidueR11}
P^{(l)}_\lambda=u-\Lim_{z\to\lambda}(\lambda-z)R_{11}(z)\,.
\end{equation}
\end{lemma}
\noindent P~r~o~o~f~.~~As the factorization 
formula~(\ref{Mfactor}) is valid for $M_1(z,\Gamma_l)$\,, any 
isolated real eigenvalue of $H_1^{(l)}$ is at the same time such 
an eigenvalue of $M_1(\cdot,\Gamma_l)$.  Since in 
$\C\setminus\overline{D(\Gamma_l)}$ the function 
$M_1(z,\Gamma_l)$ coincides with the initial (i.\,e., not 
continued yet through the continous spectrum of the entry $A_0$) 
transfer function $M_1(z)$, the point $\lambda$ must produce for 
$M_1^{-1}(z)$ the same singularity as for 
$[M_1(z,\Gamma_l)]^{-1}$.  Meanwhile, due to the 
representation~(\ref{ReprRM}) for the resolvent 
$\bR(z)=(\bH-z)^{-1}$, any singular point of the block component 
$R_{11}(z)=M_1^{-1}(z)$ of ${\bf R}(z)$ produces a singularity 
of the $\bR(z)$.  Since ${\bf H}$ is a selfadjoint operator, any 
such a point of $R_{11}(z)$ can only be a pole which is maximum 
of the first order (even if it is embedded into the continuous 
spectrum of ${\bf H}$). Since all the above is true for 
arbitrary index $l=(l_1,l_2,\ldots,l_\sm)$, $l_k=\pm1$, 
$k=1,2,\ldots,\sm,$ and since the relations~(\ref{MresiduePP}) 
hold, these considerations lead us immediately to the statements 
of the lemma.%
{\nopagebreak\mbox{\phantom{MMMM}}\hfill 
$\Box$\par\addvspace{0.25cm}}

Let $\sigma_{pri}(H_1^{(l)})$ be the set of all real isolated 
eigenvalues of the operator $H_1^{(l)}$. According to 
Lemma~\ref{labda-real} (cf. Corollary\,\ref{SigPointGeneral}) 
the set $\sigma_{pri}(H_1^{(l)})$ is the same for all 
$l=(l_1,l_2,\ldots,l_\sm)$, $l_k=\pm1$, $k=1,2,\ldots,\sm$.  
Moreover, this set coincides with the part 
$\sigma_{pri}(M_1(\cdot,\Gamma_l)$ of the set of the real 
isolated eigenvalues of the transfer function $M_1(z,\Gamma_l)$ 
belonging to ${\cal O}_{d_0/2}(A_1)$ for any  $K_B$-bounded 
contour $\Gamma_l$ satisfying the condition~(\ref{Best}),
$$
\sigma_{pri}(H_1^{(l)})=
\sigma_{pri}(M_1(\cdot,\Gamma_l))\cap{\cal O}_{d_0/2}(A_1).
$$

Since in the remainder of the Section we will consider different 
eigenvalues $\lambda\in\sigma_{pri}(H_1^{(l)})$, we will use a 
more specific notation, $\psi^{(1)}_{\lambda,j}$, 
$j=1,2,\ldots,m_\lambda$, for the respective eigenvectors of the 
$H_1^{(l)}$. The notation $m_\lambda$, $m_\lambda\leq\infty,$ 
stands for the multiplicity of the eigenvalue $\lambda$. Recall 
that every $\psi^{(1)}_{\lambda,j}$ is an eigenvector 
simultaneously for all the $H_1^{(l)}$ and 
$M_1(\lambda\pm\ri0,\Gamma_l)$, 
\mbox{$l=(l_1,l_2,\ldots,l_\sm)$} with  $l_k=\pm1$, 
$k=1,2,\ldots,\sm$ (see Lemmas \ref{LReal2} and 
{\rm\ref{LReal3}).  Since, according to Lemma~\ref{labda-real}, 
the resolvent \mbox{$(H_1^{(l)}-z)^{-1}$} has at 
$z=\lambda\in\sigma_{pri}(H_1^{(l)})$ a pole of the first order, 
the multiplicity $m_\lambda$ is, in the considered case, both 
the geometric and algebraic multiplicity of $\lambda$ (in such a 
case every element of the subspace ${\sf P}^{(l)}_\lambda\cH_1$ 
is an eigenvector of $H_1^{(l)}$ since $(H_1^{(l)}-\lambda){\sf 
P}^{(l)}_\lambda=0$).  Respective eigenvectors of the total 
matrix $\bH$ will be denoted by $\Psi_{\lambda,j}$, 
$\Psi_{\lambda,j}=(\psi^{(0)}_{\lambda,j},\psi^{(1)}_{\lambda,j})$.
It will be supposed that the $\psi^{(1)}_{\lambda,j}$ are chosen 
in such a way that the vectors $\Psi_{\lambda,j}$ are 
orthonormal, $\lal\Psi_{\lambda,j},\Psi_{\lambda',j'}\ral= 
\delta_{\lambda\lambda'}\delta_{jj'}$.  Obviously, the 
statements of Lemmas~\ref{LReal2} and~\ref{LReal3} imply that 
the eigenvectors $\Psi_{\lambda,j}$, 
$\lambda\in\sigma_{pri}(H_1^{(l)})$, $j=1,2,\ldots,m_\lambda,$ 
form an orthonormal basis in the invariant subspace of the 
operator $\bH$ corresponding to the subset 
$\sigma_{pri}(H_1^{(l)})$ of the point spectrum $\sigma_p(\bH)$ 
of $\bH$.

Let $\cH_1^{(pri)}$, $\cH_1^{(pri)}\subset\cH_1$, be the closed
span of the eigenvectors $\psi^{(1)}_{\lambda,j}$ of $H_1^{(l)}$
corresponding to the spectrum $\sigma_{pri}(H_1^{(l)})$,
$$
\cH_1^{(pri)}=\overline{{\sf V}\{\psi^{(1)}_{\lambda,j},\,
\lambda\in\sigma_{pri}(H_1^{(l)}),\, j=1,2,\ldots,m_\lambda\}}.
$$
The following statement holds.

\begin{theorem}\label{RealRieszBasis}
The system of vectors
\begin{equation}
\label{RieszBasis1}
\psi^{(1)}_{\lambda,j},\quad
\lambda\in\sigma_{pri}(H_1^{(l)}), \quad j=1,2,\ldots,m_\lambda,
\end{equation}
forms a Riesz basis of the subspace $\cH_1^{(pri)}$.
\end{theorem}

We first prove an auxiliary assertion.

\begin{lemma}\label{PosOmegaReal}
For any $l=(l_1,l_2,\ldots,l_\sm)$, $l_k=\pm1$, $k=1,2,\ldots,\sm$,
the operator $\Omega^{(l)}$ defined by Eq.~{\rm(\ref{Omega})}
is non-negative on the subspace $\cH_1^{(pri)}$.
\end{lemma}

\noindent P~r~o~o~f~.~~ It suffices to prove the assertion
for a dense subset of $\cH_1^{(pri)}$, say,
for elements
$u_1\in\cH_1^{(pri)}$ of the form
$$
   u_1=\Sum_{(\lambda,j)\in\cI} c_{\lambda,j}
   \psi^{(1)}_{\lambda,j}, \qquad c_{\lambda,j}\in\C,
$$
where $\cI$ runs through the finite subsets of the set of all possible
pairs $(\lambda,j)$ with \mbox{$\lambda\in\sigma_{pri}(H_1^{(l)})$,}\,
$j=1,2,\ldots,m_\lambda$. We have
$$
\lal\Omega^{(l)}u_1,u_1\ral=
\Sum_{(\lambda,j)\in\cI}
\Sum_{(\lambda',j')\in\cI}
c_{\lambda,j}\overline{c}_{\lambda',j'}
\Omega^{(l)}_{\lambda,j;\,\lambda',j'}
$$
with
\begin{eqnarray}
\nonumber
\Omega^{(l)}_{\lambda,j;\,\lambda',j'} &=&
\Int_{\sigma'(A_0)\cup\Gamma_l}
\lal\ K_B(d\mu)\, (H_1^{(l)}-\mu)^{-1} \psi^{(1)}_{\lambda,j},
(H_1^{(-l)}-\overline{\mu})^{-1} \psi^{(1)}_{\lambda',j'}\ral \\
\label{OmMatrEl}
&=& \Int_{\sigma'(A_0)\cup\Gamma_l}
\D\frac{\lal\ K_B(d\mu)\, \psi^{(1)}_{\lambda,j},
 \psi^{(1)}_{\lambda',j'}\ral}{(\mu-\lambda)\,(\mu-\lambda')},
\end{eqnarray}
since $\psi^{(1)}_{\lambda,j}$ and $\psi^{(1)}_{\lambda',j'}$
are eigenfunction for both $H_1^{(l)}$ and $H_1^{(-l)}$.
Due to Lemma~\ref{LReal3} (see also Remark~\ref{Remark1}
to that Lemma) one can transform the subcontours
$\Gamma_k^{l_k}\subset\Gamma_l$, $k=1,2,\ldots,\sm,$
back to the respective intervals $\Delta_k^0$, on which
$K_B(d\mu)=B_{10}E_0(d\mu)B_{01}$, even in case
$\lambda\in\Delta_k^0$ and/or $\lambda'\in\Delta_k^0$.
After such a transformation one can use Eq.~(\ref{psi0})
to express $\psi^{(0)}_{\lambda,j}$ in terms of
$\psi^{(1)}_{\lambda,j}$ and $\psi^{(0)}_{\lambda',j'}$
in terms of $\psi^{(1)}_{\lambda',j'}$. As a result, one finds
\begin{equation}
\label{OmMatrElRes}
\Omega^{(l)}_{\lambda,j;\,\lambda',j'}=
\lal\psi^{(0)}_{\lambda,j},\psi^{(0)}_{\lambda',j'}\ral
\quad \mbox{(independent of $l$)}
\end{equation}
and, hence,
$$
\lal\Omega^{(l)}u_1,u_1\ral=\|u_0\|^2\geq 0
$$
with
$u_0=\Sum_{(\lambda,j)\in\cI}c_{\lambda,j}\psi^{(0)}_{\lambda,j}$.
Thus, the operator $\Omega^{(l)}$ is non-negative on a dense
subset of $\cH_1^{(pri)}$ and, consequently, it is
non-negative on the whole subspace $\cH_1^{(pri)}$, too.
The proof of the lemma is complete.%
{\nopagebreak\mbox{\phantom{MMMM}}\hfill $\Box$\par\addvspace{0.25cm}}

\bigskip

Thus, one can introduce a new inner product in $\cH_1^{(pri)}$,
\begin{equation}
\label{NewInnerProduct}
[u_1,v_1]_{\cH_1^{(pri)}}=\lal(I_1+\Omega^{(l)})u_1,v_1\ral,
\qquad u_1,v_1\in\cH_1^{(pri)},
\end{equation}
topologically equivalent to the initial inner product
$\lal\cdot,\cdot\ral$, since $I_1+\Omega^{(l)}\geq I_1$ on
$\cH_1^{(pri)}$ and since, in view of the estimate~(\ref{Omest}), the
operator $I_1+\Omega^{(l)}$ is boundedly invertible.  (One
can even check that the restriction of $H_1^{(l)}$ on
$\cD(A_1)\cap\cH_1^{(pri)}$ does not depend on $l$
and is an operator in $\cH_1^{(pri)}$ which is self-adjoint with
respect to the inner product $[\cdot,\cdot]_{\cH_1^{(pri)}}$.)

\bigskip

P~r~o~o~f~~of Theorem \ref{RealRieszBasis}.
We prove that the vector system~(\ref{RieszBasis1})
is an orthonormal system with respect to the inner
product $[\cdot,\cdot]_{\cH_1^{(pri)}}$. Indeed, according
to Eqs.~(\ref{OmMatrEl}) and~(\ref{OmMatrElRes}) we have
$$
\lal\Omega^{(l)}\psi^{(1)}_{\lambda,j},\psi^{(1)}_{\lambda',j'}\ral=
\Omega^{(l)}_{\lambda,j;\,\lambda',j'}=
\lal\psi^{(0)}_{\lambda,j},\psi^{(0)}_{\lambda',j'}\ral.
$$
Thus,
$$
  [\psi^{(1)}_{\lambda,j},\psi^{(1)}_{\lambda',j'}]=
\lal\psi^{(1)}_{\lambda,j},\psi^{(1)}_{\lambda',j'}\ral+
\lal\psi^{(0)}_{\lambda,j},\psi^{(0)}_{\lambda',j'}\ral=
\lal\Psi_{\lambda,j},\Psi_{\lambda',j'}\ral=
\delta_{\lambda\lambda'}\delta_{jj'}.
$$
In addition, the system~(\ref{RieszBasis1}) is complete
in $\cH_1^{(pri)}$ and the inner product $[\cdot,\cdot]_{\cH_1^{(pri)}}$
is topologically equivalent to the initial inner product
$\lal\cdot,\cdot\ral$. According to a theorem of {\sc N.\,K.\,Bari}
(Theorem~VI.2.1 of~\cite{GK}) this means that
the system~(\ref{RieszBasis1}) constitutes a basis of
$\cH_1^{(pri)}$ equivalent to an orthonormal one, i.\,e.,
it is a Riesz basis. The proof is complete.%
{\nopagebreak\mbox{\phantom{MMMM}}\hfill $\Box$\par\addvspace{0.25cm}}

\bigskip


\section{\mbox{The operators $H_1^{(\lowercase{l})}$ in the case of}
\mbox{a finite-dimensional space $\cH_1$}}
\label{Sfinite-dim}

If $n_{\cH_1}=\mathop{\rm dim}\cH_1<\infty$, then the operators $A_1$
and $H_1^{(l)}$ are simply $n_{\cH_1}\times n_{\cH_1}$ scalar
matrices. In this case the resolvent of $H_1^{(l)}$ admits the
representation (see e.\,g.,~\cite{Kato}, pp.~39--44)
\begin{equation}
\label{HResExpan}
(H_1^{(l)}-z)^{-1}=\Sum_{i=1}^s \left(
-\D\frac{\sP_i^{(l)}}{z-\lambda_i^{(l)}}-
\Sum_{1\leq k\leq n_i-1}\,\,
\D\frac{[\sN_i^{(l)}]^k}{(z-\lambda_i^{(l)})^{k+1}}
\right).
\end{equation}
Here $\lambda_i^{(l)}$, $i=1,2,\ldots,s$, $s\leq n_{\cH_1}$,
stand for the different eigenvalues of $H_1^{(l)}$,
$\sP_i^{(l)}$ for the eigenprojections and $\sN_i^{(l)}$,
$\sN_i^{(l)}=(H_1^{(l)}-\lambda_i^{(l)})\sP_i^{(l)}=
\sP_i^{(l)}(H_1^{(l)}-\lambda_i^{(l)})$ for the eigennilpotents
corresponding to the $\lambda_i^{(l)}$.
The $n_i$, $n_i\geq 1$ denote the pole orders of the
resolvent $(H_1^{(l)}-z)^{-1}$ and, consequently, of the
inverse transfer function $[M_1(z,\Gamma_l)]^{-1}$ at
$z=\lambda_i^{(l)}$\,; \, if $n_i=1$, then $\sN_i^{(l)}=0$
and the eigenvalue $\lambda_i^{(l)}$ is said to be {\em
semisimple}.

Recall that all the eigenvalues $\lambda_i^{(l)}$,
$i=1,2,\ldots,s,$ belong to the set $\cO_{r_0(B)}(A_1)$, see
Corollary~\ref{HuniqueCor} to Theorem~\ref{Hunique}.

We assume that the enumeration of the $\lambda_i^{(l)}$ for
$H_1^{(l)}$ with different $l$ is co-ordinated in accordance
with the statement of Corollary~\ref{SpTheSame} to
Theorem~\ref{HlpHl2p}: If $l=(l_1,l_2,\ldots,l_\sm)$ and
$l'=(l'_1,l'_2,\ldots,l'_\sm)$ with $l'_k=l_k$ for a certain
$k=1,2,\ldots,\sm$, and $\lambda_i^{(l)}\in D_k^{l_k}$ then
$\lambda_i^{(l')}=\lambda_i^{(l)}$. Also,
$\lambda_i^{(-l)}={\overline\lambda}_i^{(l)}$ is accepted.  The
indication of $l$ in the notation $n_i$ of the pole orders
in~(\ref{HResExpan}) is omitted, since for a given $i$ the pole
order does not depend on $l$ according to the factorization
formulas~(\ref{invMinvH}) and~(\ref{HadjCons}).

The kernel $\cG_i^{(l)}=\mathop{\rm
Ker}(H_1^{(l)}-\lambda_i^{(l)})$ is called the geometric
eigenspace for the eigenvalue $\lambda_i^{(l)}$,
$i=1,2,\ldots,s$. For any $u_1\in\cG_i^{(l)}$ one has
$H_1^{(l)}u_1=\lambda_i^{(l)}u_1$.  The subspace
$\cM_i^{(l)}=\sP_i^{(l)}\cH_1$ is called the algebraic
eigenspace for $\lambda_i^{(l)}$;
$m_i\equiv\mathop{\rm dim}\cM_i^{(l)}\geq n_i$
and $\cG_i^{(l)}\subset\cM_i^{(l)}$, so that
$g_i\equiv\mathop{\rm dim}\cG_i^{(l)}\leq m_i$.
Similarly to the pole order $n_i$, the {\em algebraic
and geometric multiplicities $m_i$ and $g_i$ of the eigenvalue
$\lambda_i^{(l)}$ for a given $i$ do not depend on $l$.}
The system of algebraic eigenspaces $\cM_i^{(l)}$ is linearly
independent and complete,
\begin{equation}
\label{HMdecompos}
\cM_1^{(l)}\plusp\cM_2^{(l)}\plusp\ldots\plusp\cM_s^{(l)}=\cH_1,
\end{equation}
and, thus,
$$
%
m_1+m_2+\ldots+m_s=n_{\cH_1}
$$
and
\begin{equation}
\label{ProjSum}
   \sP_1^{(l)}+\sP_2^{(l)}+\ldots+\sP_s^{(l)}=I_1.
\end{equation}

Any vector of $\cM_i^{(l)}$ is called a {\em root vector} of the
operator $H_1^{(l)}$ corresponding to the eigenvalue
$\lambda_i^{(l)}$.

Recall some properties of the eigenprojections
and eigennilpotents:
\begin{equation}
\label{PNproperties}
\begin{array}{c}
\sP_i^{(l)}\sP_j^{(l)}=\delta_{ij}\sP_i^{(l)},
\qquad
\sP_i^{(l)}\sN_j^{(l)}=
\sN_j^{(l)}\sP_i^{(l)}=\delta_{ij}\sP_j^{(l)},\\
\sN_i^{(l)}\sN_j^{(l)}=\delta_{ij}[\sN_i^{(l)}]^{2},
\qquad [\sN_i^{(l)}]^{n_i}=0\quad\mbox{but}\quad
[\sN_i^{(l)}]^{n_i-1}\neq 0.
\end{array}
\end{equation}

The spectral representation for the $H_1^{(l)}$, in terms of
$\sP_i^{(l)}$ and $\sN_i^{(l)}$, is
\begin{equation}
\label{HSpExpansion}
   H_1^{(l)}=\Sum_{i=1}^{s}(\lambda_i^{(l)}\sP_i^{(l)}+\sN_i^{(l)}),
\end{equation}
the decomposition being unique.

As we already established in Sect. \ref{RealEigen}, if
$\lambda_i^{(l)}\in\R$, then $n_i=1$ (see Lemma~\ref{labda-real})
and thus $\sN_i^{(l)}=0$.  Therefore, in the case of a real
eigenvalue $\lambda_i^{(l)}$, any (root) vector of the subspace
$\cM_i^{(l)}$ is an eigenvector of the operator $H_1^{(l)}$
corresponding to $\lambda_i^{(l)}$, and
$\cG_i^{(l)}=\cM_i^{(l)}$.

Let $\Gamma_l\subset D_l$ be a $K_B$-bounded contour
satisfying the condition~(\ref{Best}). According to
Theorem~\ref{SpHalfVic}, the spectrum of the transfer function
$M_1(\cdot,\Gamma_l)$ is represented in the set
$\cO_{d_0(\Gamma_l)/2}(A_1)$ [and even in the set
$\biggl(\R\cup D(\Gamma_l)\biggr)\cap\cO_{d_{\rm max}/2}(A_1)$,
according to Corollary~\ref{HMspectr} to the theorem] just by
the spectrum of the operator $H_1^{(l)}$.  Thus, the transfer
function $M_1(\cdot,\Gamma_l)$ has in
$\cO_{d_0(\Gamma_l)/2}(A_1)$ only discrete spectrum
consisting of the eigenvalues $\lambda_i^{(l)}$,
$i=1,2,\ldots,s$. Due to Eq.~(\ref{invMinvH}) the inverse
function $[M_1(\cdot,\Gamma_l)]^{-1}$ has
poles at $z=\lambda_i^{(l)}$
of the same orders $n_i$ as the resolvent
$(H_1^{(l)}-z)^{-1}$.

\begin{lemma}\label{MeqPN}
The eigenprojections $\sP_i^{(l)}$ and eigennilpotents
$\sN_i^{(l)}$, $i=1,2,\ldots,s,$ of the operator $H_1^{(l)}$
satisfy the equations
\begin{equation}
\label{MPeq}
M_1(\lambda_i^{(l)},\Gamma_l)\sP_i^{(l)}=
\sN_i^{(l)}-\Sum_{1\leq k\leq n_i-1} \,\,\D\frac{1}{k!}\,
V_1^{(k)}(\lambda_i^{(l)},\Gamma_l)[\sN_i^{(l)}]^k
\end{equation}
where $V_1^{(k)}(\lambda,\Gamma_l)$
stands for the $k$-th derivative of the function
$V_1(z,\Gamma_l)$, defined by Eq.~{\rm(\ref{Mcmpl})},
at $z=\lambda$,
\begin{equation}
\label{Vderiv}
V_1^{(k)}(\lambda,\Gamma_l)=(-1)^k\,k!
\Int_{\sigma'(A_0)\cup\Gamma_l} K_B(d\mu)
\D\frac{1}{(\lambda-\mu)^{k+1}},
\qquad k=0,1,2,\ldots\,.
\end{equation}
\end{lemma}

\noindent P~r~o~o~f~.~~Write the basic equation~(\ref{MainEq})
for $H_1^{(l)}$ as follows
\begin{equation}
\label{BasicM}
A_1=H_1^{(l)}-\Int_{\sigma'(A_0)\cup\Gamma_l} K_B(d\mu)\,
(H_1^{(l)}-\mu)^{-1}\,.
\end{equation}
Multiplying both parts of Eq.~(\ref{BasicM}) by $\sP_i^{(l)}$
from the right and taking into account equalities
\begin{equation}
\label{H1Pl}
H_1^{(l)}\sP_i^{(l)}=\lambda_i^{(l)}\sP_i^{(l)}+\sN_i^{(l)}
\end{equation}
and
\begin{equation}
\label{H1ResPl}
(H_1^{(l)}-z)^{-1}\sP_i^{(l)}=
-\D\frac{\sP_i^{(l)}}{z-\lambda_i^{(l)}}-
\Sum_{1\leq k\leq n_i-1}\,\,
\D\frac{[\sN_i^{(l)}]^k}{(z-\lambda_i^{(l)})^{k+1}}
\end{equation}
which follow from Eqs.~(\ref{HSpExpansion}),  (\ref{HResExpan})
and~(\ref{PNproperties}), one finds
$$
A_1\sP_i^{(l)}=\lambda_i^{(l)}\sP_i^{(l)}+\sN_i^{(l)}-
\Int_{\sigma'(A_0)\cup\Gamma_l} K_B(d\mu) \left(
-\D\frac{\sP_i^{(l)}}{\mu-\lambda_i^{(l)}}-
\Sum_{1\leq k\leq n_i-1}\,\,
\D\frac{[\sN_i^{(l)}]^k}{(\mu-\lambda_i^{(l)})^{k+1}}
\right)\,.
$$
But according to Eqs.~(\ref{Mcmpl}) and~(\ref{Vderiv})
this is just the equation~(\ref{MPeq}) which we wanted
to prove.%
{\nopagebreak\mbox{\phantom{MMMM}}\hfill $\Box$\par\addvspace{0.25cm}}

\begin{remark}\label{MNilpEq}
The eigennilpotents $\sN_i^{(l)}$, $i=1,2,\ldots,s$
also satisfy the equations
\begin{equation}
\label{MNeq}
\begin{array}{rcl}
M_1(\lambda_i^{(l)},\Gamma_l)[\sN_i^{(l)}]^{n_i-p} & = &
[\sN_i^{(l)}]^{n_i-p+1}-\Sum_{1\leq k\leq p-1}\,\,\D\frac{1}{k!}\,
V_1^{(k)}(\lambda_i^{(l)},\Gamma_l)[\sN_i^{(l)}]^{n_i-p+k}, \\
 && p=1,2,\ldots,n_i-1.
\end{array}
\end{equation}
\end{remark}

One obtains the Eqs.~(\ref{MNeq}) simply by multiplying both
parts of Eqs.~(\ref{MPeq}) from the right by
$[\sN_i^{(l)}]^{n_i-p}$, $p=1,2,\ldots,\mbox{$n_i-1$}$.

\begin{remark}\label{MNPrepr}
Writing $A_1$ as $A_1=H_1^{(l)}-V_1(H_1^{(l)},\Gamma_l)$
and then using Eqs.~{\rm(\ref{HResExpan}), (\ref{HSpExpansion})}
and~{\rm(\ref{Vderiv})} one can represent the transfer function
$M_1(z,\Gamma_l)$ as
\begin{eqnarray*}
M_1(z,\Gamma_l)  &=& \Sum_{i=1}^s \biggl\{
(\lambda_i^{(l)}-z)\sP_i^{(l)}+\sN_i^{(l)}+V_1(z,\Gamma_l)\sP_i^{(l)} \\
 & & -V_1(\lambda_i^{(l)},\Gamma_l)\sP_i^{(l)}
 -\Sum_{1\leq k\leq n_i-1} \,\,\D\frac{1}{k!}\,
V_1^{(k)}(\lambda_i^{(l)},\Gamma_l)[\sN_i^{(l)}]^k
\biggr\}.
\end{eqnarray*}
\end{remark}

\bigskip

In fact, the Eqs.~(\ref{MPeq}) considered together with the
conditions~(\ref{PNproperties}) determine the eigenprojections
$\sP_i^{(l)}$ and eigennilpotents $\sN_i^{(l)}$,
$i=1,2,\ldots,s$ uniquely, at least under the additional
condition
\begin{equation}
\label{AdditCond}
\biggl\|\Sum_{i=1}^s(\lambda_i^{(l)}\sP_i^{(l)}+\sN_i^{(l)})-A_1\biggr\|
< r_{\rm max}(\Gamma_l),
\end{equation}
where $r_{\rm max}(\Gamma_l)$ is given by Eq.~(\ref{rmax}).
Namely, the following assertion holds.

\begin{theorem}\label{NPinverse}
Let $\Gamma_l\subset D_l$ be a $K_B$-bounded contour satisfying
the condition~{\rm(\ref{Best})}. Also, let $\lambda_i^{(l)}$,
$i=1,2,\ldots,s$ be the eigenvalues and $n_i$, $i=1,2,\ldots,s$,
the respective pole orders of the transfer function
$M_1(z,\Gamma_l)$ in the domain $\cO_{d_0(\Gamma_l)/2}(A_1)$.
Then the system of Eqs.~{\rm(\ref{ProjSum})},
{\rm(\ref{PNproperties})} and~{\rm(\ref{MPeq})} for \,
$i,j=1,2,\ldots,s$ under the condition~{\rm(\ref{AdditCond})}
determines uniquely the complete system of eigenprojections and
eigennilpotents for the operator $H_1^{(l)}$.
\end{theorem}

\noindent P~r~o~o~f~.~~Obviously, one has to prove only
the uniqueness of the solution for the
system~(\ref{PNproperties}), (\ref{MPeq}), since the existence of
such a solution is already guaranteed by Lemma~\ref{MeqPN}.

Let $\tsP_i^{(l)}$, $\tsN_i^{(l)}$, $i=1,2,\ldots,s$
be a solution of the system~{\rm(\ref{ProjSum})}, 
(\ref{PNproperties}) and (\ref{MPeq})
under the condition~(\ref{AdditCond}). Consider the operator
\begin{equation}
\label{tildeHSpExpansion}
   \tH_1^{(l)}=\Sum_{i=1}^{s}
(\lambda_i^{(l)}\tsP_i^{(l)}+\tsN_i^{(l)}).
\end{equation}
Since the $\tsP_i^{(l)}$, $\tsN_i^{(l)}$
satisfy Eqs.~(\ref{ProjSum}), (\ref{PNproperties}),
the equation~(\ref{tildeHSpExpansion}) is at the same time
the spectral representation for $\tH_1^{(l)}$
and, consequently,
\begin{equation}
\label{tHResExpan}
(\tH_1^{(l)}-z)^{-1}=\Sum_{i=1}^s \left(
-\D\frac{\tsP_i^{(l)}}{z-\lambda_i^{(l)}}-
\Sum_{1\leq k\leq n_i-1}\,\,
\D\frac{[\tsN_i^{(l)}]^k}{(z-\lambda_i^{(l)})^k}
\right)\,.
\end{equation}
Rewrite Eqs.~(\ref{MPeq}) for $\tsP_i^{(l)}$ and $\tsN_i^{(l)}$
as
\begin{equation}
\label{tMPeq}
\lambda_i^{(l)}\tsP_i^{(l)}+\tsN_i^{(l)}=A_1\tsP_i^{(l)}
+V_1(\lambda_i^{(l)},\Gamma_l)\tsP_i^{(l)}
+\Sum_{1\leq k\leq n_i-1} \,\,\D\frac{1}{k!}\,
V_1^{(k)}(\lambda_i^{(l)},\Gamma_l)[\tsN_i^{(l)}]^k
\end{equation}
and represent then the derivatives
$V_1^{(k)}(\lambda_i^{(l)},\Gamma_l)$ by~(\ref{Vderiv}). Summing
over $i$ in~(\ref{tMPeq}) and taking into
account~(\ref{tHResExpan}) one finds
\begin{equation}
\label{tHBasic}
\tH_1^{(l)}=A_1+\Int_{\sigma'(A_0)\cup\Gamma_l}
K_B(d\mu)\,(\tH_1^{(l)}-\mu)^{-1},
\end{equation}
that is, $\tH_1^{(l)}$ satisfies the basic
equation~(\ref{MainEq}) and the difference
$\tX_1^{(l)}=\tH_1^{(l)}-A_1$ the basic equation~(\ref{MainEqC}).
Meanwhile the condition~(\ref{AdditCond})
means $\|\tX_1^{(l)}\|<r_{\rm max}(\Gamma_l)$.
Then it follows from Theorem~\ref{Solvability}
that $\tX_1^{(l)}=X_1^{(l)}$ and, hence, $\tH_1^{(l)}=H_1^{(l)}$.
Since the spectral representation~(\ref{HSpExpansion})
for $H_1^{(l)}$ is unique, one must conclude that
$\tsP_i^{(l)}=\sP_i^{(l)}$ and $\tsN_i^{(l)}=\sN_i^{(l)}$,
$i=1,2,\ldots,s$ and this completes the proof.%
{\nopagebreak\mbox{\phantom{MMMM}}\hfill $\Box$\par\addvspace{0.25cm}}

Thus, the eigenprojections $\sP_i^{(l)}$ and eigennilpotents
$\sN_i^{(l)}$, $i=1,2,\ldots,s$ of the operator $H_1^{(l)}$ can
be called the eigenprojections and eigennilpotents of the
transfer function $M_1(z,\Gamma_l)$ in the domain
$\cO_{d_0(\Gamma_l)/2}(A_1)$.  By Lemma~\ref{MeqPN} and
Theorem~\ref{NPinverse} this definition is correct.

\bigskip

Recall that a basis $\{e_j\}_{j=1}^{n_{\cH_1}}$ of the space $\cH_1$
is said to be {\em adapted} to the
decomposition~(\ref{HMdecompos}) if the first several elements of
$\{e_j\}_{j=1}^{n_{\cH_1}}$ belong to $\cM_1$, the following
several elements belong to $\cM_2$ and so on. In the case considered here
such an adapted basis consists of root vectors corresponding to
the eigenvalues $\lambda_i^{(l)}$, $i=1,2,\ldots,s,$ including
the resonances.  One can choose, in particular, a basis consisting
of the eigenvectors and associated vectors
which reduces every eigennilpotent
$\sN_i^{(l)}$ to Jordan canonical form (see e.\,g.~\cite{Kato},
pp.~22, 43).

To conclude the section we consider the case
where all the eigenvalues $\lambda_i^{(l)}$ are semisimple,
i.\,e., $\sN_i^{(l)}=0$, $i=1,2,\ldots,s.$

Let $\{\psi_{ij}^{(l)},\,j=1,2,\ldots,m_i\}$ be a basis of the
eigenspace $\cG_i^{(l)}$ ($\cG_i^{(l)}=\cM_i^{(l)}$ in this case).
The union of these bases for $i=1,2,\ldots,s$ is a basis of the
space $\cH_1$. Denote by
$\{\varphi_{ij}^{(l)},\,i=1,2,\ldots,s;\,j=1,2,\ldots,m_i\}$
the biorthogonal basis,
$\lal\psi_{ij}^{(l)},\varphi_{i'j'}^{(l)}\ral=\delta_{ii'}\delta_{jj'}$.
Then the vectors \mbox{$\varphi_{ij}^{(l)}$,
$j=1,2,\ldots,m_i$} for a given $i=1,2,\ldots,s$
are automatically eigenvectors of the operator $H_1^{(l)*}$ corresponding
to the eigenvalue $\overline{\lambda}_i^{(l)}$
while
$\sP_i^{(l)}=\Sum_{j=1}^{m_i}\psi_{ij}^{(l)}
\lal\cdot,\varphi_{ij}^{(l)}\ral$. Eq.~(\ref{HOmega}) implies that
$$
\psi_{ij}^{(-l)}=(I_1+\Omega^{(-l)})^{-1}\varphi_{ij}^{(l)},
\quad j=1,2,\ldots,m_i
$$
are the linearly independent eigenvectors of the operator
$H_1^{(-l)}$ corresponding to the eigenvalue
$\lambda_i^{(-l)}=\overline{\lambda}_i^{(l)}$ and
\begin{equation}
\label{OmBiorth}
\lal\psi_{ij}^{(l)},(I_1+\Omega^{(-l)})\psi_{i'j'}^{(-l)}\ral=
\lal(I_1+\Omega^{(l)})\psi_{ij}^{(l)},\psi_{i'j'}^{(-l)}\ral=
\delta_{ii'}\delta_{jj'}\,.
\end{equation}
Therefore, one comes to the following assertion.

\begin{lemma}\label{SemisimpleProj}
If all the eigenvalues $\lambda_i^{(l)}$, $i=1,2,\ldots,s$
of the operator $H_1^{(l)}$ are semisimple, then the spectral
projections $\sP_i^{(l)}$ can be written as
\begin{equation}
\label{PpsiOm}
  \sP_i^{(l)}=\Sum_{j=1}^{m_i}\psi_{ij}^{(l)}
  \lal\cdot,\psi_{ij}^{(-l)}\ral(I_1+\Omega^{(l)})
\end{equation}
where the eigenvectors $\psi_{ij}^{(l)}$ and $\psi_{ij}^{(-l)}$
of the operators $H_1^{(l)}$ and $H_1^{(-l)}$
{\rm(}$H_1^{(l)}\psi_{ij}^{(l)}=\lambda_i^{(l)}\psi_{ij}^{(l)}$,
$H_1^{(-l)}\psi_{ij}^{(-l)}=
\overline{\lambda}_i^{(l)}\psi_{ij}^{(-l)}${\rm)}
are normalized according to Eqs.~{\rm(\ref{OmBiorth})}.
At the same time
$$
H_1^{(l)}=\Sum_{i=1}^s\lambda_i^{(l)}
\Sum_{j=1}^{m_i}\psi_{ij}^{(l)}
  \lal\cdot,\psi_{ij}^{(-l)}\ral(I_1+\Omega^{(l)}).
$$
\end{lemma}
\begin{remark}\label{SemisMResidue}
It follows from the relations~{\rm(\ref{MresiduePP})}
and~{\rm(\ref{PpsiOm})} that, in the case considered here, the residues
$P_i^{(l)}$ of the transfer function $M_1(z,\Gamma)$ at
$z=\lambda_i^{(l)}$, $i=1,2,\ldots,s$ read as follows:
$$
    P_i^{(l)}=\Sum_{j=1}^{m_i}\psi_{ij}^{(l)}
  \lal\cdot,\psi_{ij}^{(-l)}\ral\,.
$$
The total sum of these residues represents an invertible operator
and
$$
  \left(\Sum_{i=1}^{s}P_i^{(l)}\right)^{-1}=(I_1+\Omega^{(l)})^{-1}\,.
$$
\end{remark}

\section{Completeness and basis properties of the
$H_1^{(\lowercase{l})}$
root vectors in the case of an infinite-dimensional space $\cH_1$}
\label{BasisnessGeneral}

In the present Section we restrict ourselves to the case where
the entry $A_1$ has pure discrete spectrum only, i.\,e., the resolvent
$R_1(z)=(A_1-z)^{-1}$ is a compact operator in $\cH_1$ for any
$z\in\rho(A_1)$.

\begin{lemma}\label{HlCompact}
If the entry $A_1$ has compact resolvent, then the operators
$H_1^{(l)}$ have compact resolvents, too.
\end{lemma}
\noindent This statement is a simple consequence of
Theorem~IV.3.17 of~\cite{Kato}, since $H_1^{(l)}$ since the
difference \mbox{$H_1^{(l)}-A_1=X^{(l)}$} is a bounded operator
(see Theorem~\ref{Solvability}).%
{\nopagebreak\mbox{\phantom{MMMM}}\hfill $\Box$\par\addvspace{0.25cm}}

\begin{lemma}\label{XlCompact}
If the entry $A_1$ has compact resolvent, then the solutions
$X^{(l)}$ of the basic equation~{\rm(\ref{MainEqC})} are compact
operators.
\end{lemma}

\noindent P~r~o~o~f~.~~According to Lemma \ref{HlCompact},
the resolvent \mbox{$(H_1^{(l)}-\mu)^{-1}$} is a compact operator
for any $\mu$ belonging
to an arbitrary $K_B$-bounded contour $\Gamma_l$
satisfying the condition~(\ref{Best}), since for such a contour
\mbox{$\dist\{\sigma(H_1^{(l)},\Gamma_l)\}>d_0(\Gamma_l)/2>0.$}
Thus, any finite integral sum for the integral defining $X^{(l)}$,
\mbox{$
X^{(l)}=\Int_{\sigma'(A_0)\cup\Gamma_l}
K_B(d\mu)\,(H_1^{(l)}-\mu)^{-1},
$}
is a compact operator. But under the $K_B$-boundedness
condition~(\ref{Bbound}) the integral sums converge to $X^{(l)}$
with respect to the operator norm topology
(see Appendix~\ref{IntOpMer}). Thus, $X^{(l)}$ must be
a compact operator.%
{\nopagebreak\mbox{\phantom{MMMM}}\hfill $\Box$\par\addvspace{0.25cm}}

\bigskip

Denote by $\cH_{1,\lambda}^{(l)}$ the algebraic eigenspace of
$H_1^{(l)}$ corresponding to an eigenvalue $\lambda$,
$\cH_{1,\lambda}^{(l)}=\sP_\lambda^{(l)}\cH_1$ where the
eigenprojection $\sP_\lambda^{(l)}$ is given by
Eq.~(\ref{Plambda}).  Let $m_\lambda$ be the algebraic
multiplicity, $m_\lambda=\mathop{\rm dim}\cH_\lambda^{(l)}$,
$m_\lambda<\infty$, and $\sN_\lambda^{(l)}$ respective
eigennilpotent,
$\sN_\lambda^{(l)}=(H_1^{(l)}-\lambda)\sP_\lambda^{(l)}$.  The
eigenprojections $\sP_\lambda^{(l)}$ and eigennilpotents
$\sN_\lambda^{(l)}$ for different $\lambda$ again satisfy
Eqs.~(\ref{PNproperties}) as well as Eqs.~(\ref{H1Pl}) and
(\ref{H1ResPl}) are valid (see~\cite{Kato}, \S\,6.5 of
Chapter~III). Repeating literally the proof of Lemma~\ref{MeqPN}
one can check easily that for the case considered now the 
statement this lemma is still valid.

Let $\psi^{(l)}_{\lambda,i},$ $i=1,2,\ldots,m_\lambda,$ be the root
vectors of $H_1^{(l)}$ forming a basis of the algebraic
eigenspace $\cH_{1,\lambda}^{(l)}$. In the following we will try to give
an answer on the question when the union of such bases in
$\lambda$ forms a basis of the total space $\cH_1$. But, in
any case, we already have  an assertion regarding
completeness of the system
\begin{equation}
\label{RootVecSystem}
\{\psi^{(l)}_{\lambda,i},\,\,\lambda\in\sigma(H_1^{(l)}),\,\,
                              i=1,2,\ldots,m_\lambda\}.
\end{equation}

\begin{theorem}\label{Completeness}
The closure of the linear span of the system~{\rm(\ref{RootVecSystem})}
coincides with $\cH_1$,
$$
\overline{{\sf V}\{\psi^{(l)}_{\lambda,i},\,\,
          \lambda\in\sigma(H_1^{(l)}),\,\,
          i=1,2,\ldots,m_\lambda\}}=\cH_1\,.
$$
\end{theorem}
\noindent This assertion is a particular case
of Theorem~V.10.1 from~\cite{GK}.%
{\nopagebreak\mbox{\phantom{MMMM}}\hfill $\Box$\par\addvspace{0.25cm}}

\bigskip

The following statement concerns the case where the basis
property of the system~(\ref{RootVecSystem}) follows immediately
from the general basis property (Theorem~\ref{RealRieszBasis})
of the eigenvectors corresponding to the real isolated point
spectrum eigenvalues of the operator $H_1^{(l)}$.

\begin{theorem}\label{CSpectrumAbove}
Let the entry $A_1$ have compact resolvent and be semibounded
from below. Suppose the set
$\mathop{{\bigcup}}\limits_{k=1}^{\sm}\Delta^0_k$ is bounded
from above, i.\,e., $\mu_\sm^{(2)}<\infty$.  Then the operator
$H_1^{(l)}$ has only a finite number of complex eigenvalues
{\rm(}resonances{\rm)}. It can be represented as
$$
   H_1^{(l)}=H_{1,R}^{(l)}+H_{1,C}^{(l)}
$$
with $H_{1,R}^{(l)}=H_1^{(l)}\sP_R^{(l)}$ and
$H_{1,C}^{(l)}=H_1^{(l)}\sP_C^{(l)}$ where $\sP_R^{(l)}$
and $\sP_C^{(l)}$ stand for the projections on the invariant
subspaces $\cH_{1,R}^{(l)}$,
$\cH_{1,R}^{(l)}=\sP_R^{(l)}\cH_1$, and $\cH_{1,C}^{(l)}$,
$\cH_{1,C}^{(l)}=\sP_C^{(l)}\cH_1$, corresponding respectively
to the real and complex spectrum of $H_1^{(l)}$.
The restriction of $H_{1,R}^{(l)}$ to $\cD(A_1)\cap\cH_{1,R}^{(l)}$
represents an operator which is similar to a selfadjoint
one while for the finite-dimensional component
$\reduction{H_{1,C}^{(l)}}{\cH_{1,C}^{(l)}}$ one can find the
eigenprojections and eigennilpotents using the statements
of {\rm Lemma \ref{MeqPN}} and {\rm Theorem~\ref{NPinverse}}.
Combining a basis of the subspace $\cH_{1,C}^{(l)}$
consisting of the root vectors for
$\reduction{H_{1,C}^{(l)}}{\cH_{1,C}^{(l)}}$ with a Riesz basis
of the subspace $\cH_{1,R}^{(l)}$ constructed from the
eigenvectors of $\reduction{H_{1,R}^{(l)}}{\cH_{1,R}^{(l)}}$
one gets a Riesz basis of the space $\cH_1$.
\end{theorem}

\noindent P~r~o~o~f~.~~According to Theorem \ref{SpHalfVic}
the complex spectrum of the operator $H_1^{(l)}$
belongs to the set $D_l\cap\cO_{r_0(B)}(A_1)$ and even
to the domains $D(\Gamma_l)$ restricted by
$\mathop{{\bigcup}}\limits_{k=1}^{\sm}\Delta^0_k$ and
arbitrary $K_B$-bounded contours $\Gamma_l\subset D_l$
satisfying the condition~(\ref{Best}).
The rest of the spectrum of $H_1^{(l)}$ is real.
Obviously, in the case concerned,
the set $D(\Gamma_l)\cap\cO_{r_0(B)}(A_1)$ is bounded
even if the domain $D_l$ is unbounded. Since the spectrum of
$H_1^{(l)}$ is discrete (see Lemma~\ref{HlCompact}),
only a finite number of the $H_1^{(l)}$ eigenvalues
can be situated in $D(\Gamma_l)\cap\cO_{r_0(B)}(A_1)$ and
these eigenvalues
generate the finite-dimensional eigenprojections.
Thus, the projection $\sP^{(l)}_C$, being a sum of the
individual eigenprojections, is finite-dimensional, too.
Multiplying both parts of the
basic equation for $H_1^{(l)}$, written as Eq.~(\ref{BasicM}),
by $\sP^{(l)}_C$ from the right and separating
the eigenprojections and eigennilpotents corresponding
to the individual resonances, one further comes to
the statements of Lemma~{\rm\ref{MeqPN}} and
Theorem~{\rm\ref{NPinverse}} restricted to the subspace
$\cH_{1,C}^{(l)}$. Noting that
\mbox{$\cH_1=\cH_{1,R}^{(l)}\plusp\cH_{1,C}^{(l)}$},
and then applying Theorem~\ref{RealRieszBasis}
one gets the remaining statements.%
{\nopagebreak\mbox{\phantom{MMMM}}\hfill $\Box$\par\addvspace{0.25cm}}

\begin{remark}\label{SemibAbove}
The statement of {\rm Theorem~\ref{CSpectrumAbove}}
remains valid if the entry $A_1$ has compact resolvent
and is semibounded from above while the set
$\mathop{\cup}_{k=1}^\sm\Delta_k^0$ is semibounded from
below, i.\,e., $\mu_1^{(1)}>-\infty.$
\end{remark}

In what follows we need a few definitions and statements
from Chapter~VI of the book~\cite{GK}.

Let $\{e_k\}_{k=1}^\infty$  be a basis of a Hilbert space
$\cN$. If there exists an orthonormal basis
$\{e'_k\}_{k=1}^\infty$ of $\cN$ such that
$$
    \Sum_{k=1}^\infty\|e_k-e'_k\|^2<\infty,
$$
then the basis $\{e_k\}_{k=1}^\infty$ is said to be quadratically
close to an orthonormal basis. Also, such a basis is called
a Bari basis. Any Bari basis is at the same time a
Riesz basis (see Theorem~VI.2.3 of~\cite{GK}).

A sequence $\{\cN_k\}_{k=1}^\infty$ of non-zero subspaces
$\cN_k\subset\cN$ is said to be a basis (of subspaces) of the
Hilbert space $\cN$ if any vector $x\in\cN$
can be expanded in a unique way in a series of the form
$$
        x=\Sum_{k=1}^\infty x_k
$$
where $x_k\in\cN_k$, $k=1,2,\ldots$\,\,.

A sequence $\{\cN_k\}_{k=1}^\infty$ of non-zero subspaces
$\cN_k\subset\cN$ is said to be $\omega$-linearly independent
if the equality
$$
        \Sum_{k=1}^\infty x_k=0, \quad
        x_k\in\cN_k,\,\, k=1,2,\ldots\,,
$$
is not possible for
$$
        0<\Sum_{k=1}^\infty \|x_k\|^2<\infty\,.
$$

A sequence $\{\cN_k\}_{k=1}^\infty$ of subspaces
$\cN_k\subset\cN$ is said to be quadratically close to an
orthogonal basis (of subspaces) of the the space $\cN$ if there
exists a sequence of pairwise orthogonal subspaces
$\cN'_k\subset\cN$ such that
$\opla_{k=1}^\infty\cN'_k=\cN$ and
$$
    \Sum_{k=1}^\infty \|\sP_{\cN_k}-\sP_{\cN'_k}\|^2<\infty
$$
where $\sP_{\cN_k}$ and $\sP_{\cN'_k}$, $k=1,2,\ldots,$
stand for the orthogonal projections of $\cN$ onto
$\cN_k$ and $\cN'_k$, respectively.

The {\em minimal angle}\,\, $\phi(\cN',\cN'')$, $0\leq\phi\leq\pi/2,$
between two subspaces $\cN'$ and $\cN''$ is defined as
$$
\cos\phi(\cN',\cN'')=\Sup\limits_{\mbox{\scriptsize$
\begin{array}{c}
                x'\in\cN',\,\,x''\in\cN''\\
                  \|x'\|=\|x''\|=1
\end{array}$}}
|\lal x',x''\ral|\,.
$$

\begin{theorem}\label{Markus1}
{\rm\cite{Mar}} Let $\{\cN_k\}_{k=1}^\infty$
be a complete, $\omega$-linearly independent sequence
of finite-dimensional subspaces in $\cN$ such that
\begin{equation}
\label{Cosines}
\Sum^\infty_{\mbox{\scriptsize$
\begin{array}{c}
                i,j=1\\
                i\neq j
\end{array}$}}
\cos^2\phi(\cN_i,\cN_j)<\infty\,.
\end{equation}
Then $\{\cN_k\}_{k=1}^\infty$ is a basis of the space
$\cN$, quadratically close to an orthogonal one.
\end{theorem}

\begin{theorem}\label{GKMarkus1}
{\rm(Proposition VI.5.6 of~\cite{GK})}
If the condition~{\rm (\ref{Cosines})} in
{\rm Theorem~\ref{Markus1}} can be replaced with
\begin{equation}
\label{Cosines1}
\Sum^\infty_{\mbox{\scriptsize$
\begin{array}{c}
                i,j=1\\
                i\neq j
\end{array}$}}
\Min\{\nu_i,\nu_j\}\,\cos^2\phi(\cN_i,\cN_j)<\infty
\end{equation}
where $\nu_i=\mathop{\rm dim}\cN_i$, then the union of orthonormal
vector bases of the subspaces $\cN_k,$ $k=1,2,\ldots,$
forms a Bari basis of the space $\cN$.
\end{theorem}

The following statement is a particular case of a more general 
proposition from the final part of \S\,III.7.3 of 
Ref.~\rm\cite{GK} regarding criteria for a linear operator to 
belong to a certain class of compact operators. 
It also represents Theorem~5 of \S\,11.5 of 
Ref.~\cite{BS}.

\begin{theorem}
\label{CriterionTrace}
If for a bounded linear operator $T$ acting in a Hilbert space
$\cN$ the condition
$$
   \sum\limits_{k=1}^\infty \|Te_k\|<\infty
$$
is valid for some orthonormal basis $\{e_k\}_{k=1}^\infty$
of the space $\cN$ then $T$ is an operator of the trace class.
\end{theorem}

Let us return to the operators $H_1^{(l)}$, now in the case
where the intersection
$\biggl(\mathop{{\bigcup}}\limits_{k=1}^{\sm}\Delta^0_k\biggr)
\cap\sigma(A_1)$ includes infinitely many points and, thus,
\mbox{$\mathop{\rm dim}\cH_1=\infty$}.

For the sake of simplicity we assume that the entry $A_1$ as
above, has compact resolvent and is semibounded from below.
Then the previous assumption means that at least the interval
$\Delta_\sm^0$ is infinite,
$\Delta_\sm^0=(\mu_\sm^{(1)},+\infty)$.  The eigenvalues
$\lambda_i^{(A_1)}$, $i=1,2,\ldots\,$, of the operator $A_1$
will be enumerated in increasing order,
$\lambda_1^{(A_1)}<\ldots<\lambda_i^{(A_1)}
<\lambda_{i+1}^{(A_1)}<\ldots\,$\, and
\mbox{$\Lim_{i\to\infty}\lambda_i^{(A_1)}=+\infty$} exists.

Suppose further that there is a number $i_0$ such that
for any $i\geq i_0$
\begin{equation}
\label{DifLambda}
  \lambda_i^{(A_1)}-\lambda_{i-1}^{(A_1)}>2r>2r_0(B),
\end{equation}
with a fixed value $r$ while $r_0(B)$ is given by~(\ref{r0}).
Let $\gamma_0$ be a circle centered at
\mbox{$z=(\lambda_1^{(A_1)}+\lambda_{i_0-1}^{(A_1)})/2$}
and having the radius
\mbox{$(\lambda_{i_0-1}^{(A_1)}-\lambda_1^{(A_1)})/2+r$}
while the $\gamma_i$ for $i\geq i_0$ are the circles with
centers $\lambda_i^{(A_1)}$ and the radius $r$. Obviously,
the union
$\Inter\gamma_0\mathop{\bigcup}\limits_{i\geq i_0}\Inter\gamma_i$
of the interiors of the circles $\gamma_i$, $i=0,i_0,i_0+1,\ldots,$
covers all the spectrum of $H_1^{(l)}$, since
$\sigma(H_1^{(l)})\subset\cO_{r_0(B)}(A_1)$. At the same time
\begin{equation}
\label{InterInter}
\overline{\Inter\gamma_i}\cap\overline{\Inter\gamma_j}
=\emptyset\quad \mbox{if \, $i\neq j$}.
\end{equation}
Thus, one can introduce the projections
\begin{equation}
\label{sQDef}
\sQ_i^{(l)}=
-\D\frac{1}{2\pi\ri}\Int_{\gamma_i} dz\,\,(H_1^{(l)}-z)^{-1},
\qquad i=0,i_0,i_0+1,\ldots,
\end{equation}
and, then, the subspaces
$
    \cN_i^{(l)}=\sQ_i^{(l)}\cH_1
$
which are invariant under $H_1^{(l)}$. Due to Eqs.~(\ref{InterInter})
one has
\begin{equation}
\label{QQ}
\sQ^{(l)}_i\sQ^{(l)}_j=\delta_{ij}\sQ^{(l)}_i.
\end{equation}
Each projection $\sQ_i^{(l)}$ represents a sum of the
eigenprojections~(\ref{Plambda}) corresponding to the eigenvalues
$\lambda^{(l)}$ of $H_1^{(l)}$ belonging to $\Inter\gamma_i$.
Since the algebraic eigenspaces for different eigenvalues
are linearly independent, the dimension
$\mathop{\rm dim}\cN_i^{(l)}$ coincides with sum of the algebraic
multiplicities for the $\lambda^{(l)}$ lying inside $\gamma_i$.
We introduce also the (orthogonal) projections
\begin{equation}
\label{sPADef}
\sP_i^{(A_1)}=
-\D\frac{1}{2\pi\ri}\Int_{\gamma_i} dz\,\,(A_1-z)^{-1},
\qquad i=0,i_0,i_0+1,\ldots\,\,\,.
\end{equation}
Obviously, for $i\geq i_0$ the projections $\sP_i^{(A_1)}$
are simply the eigenprojections of the entry $A_1$ corresponding
to the eigenvalues $\lambda_i^{(A_1)}$ while $\sP_0^{(A_1)}$ is
the sum of the eigenprojections for $A_1$ corresponding to the
eigenvalues $\lambda_1^{(A_1)},\lambda_2^{(A_1)},
\ldots,\lambda_{i_0-1}^{(A_1)}.$
In the following by $\varphi_{ij}^{(A_1)}$,
$j=1,2,\ldots,n_i^{(A_1)}$, $n_i^{(A_1)}=\mathop{\rm
dim}\sP_i^{(A_1)}\cH_1<\infty$ we understand an orthonormal
basis of the subspace $\sP_i^{(A_1)}\cH_1$.  For $i\geq i_0$ the
vectors $\varphi_{ij}^{(A_1)}$ are automatically eigenvectors of
$A_1$,
$A_1\varphi_{ij}^{(A_1)}=\lambda_i^{(A_1)}\varphi_{ij}^{(A_1)}$.
The sequence
$\{\varphi_{ij}^{(A_1)},\,\,i=0,i_0,i_0+1,\ldots,\,\,j=1,\ldots,
n_i^{(A_1)}\}$ forms an orthonormal basis of $\cH_1$.

\begin{lemma}\label{NomegaLInd}
Under the condition~{\rm(\ref{DifLambda})} the sequence
\begin{equation}
\label{Nsequence}
  \{\cN_i^{(l)},\,\,i=0,i_0,i_0+1,\ldots\}
\end{equation}
of subspaces $\cN_i^{(l)}=\sQ_i^{(l)}\cH_1$
is $\omega$-linearly independent
and complete in $\cH_1$.
\end{lemma}

\noindent P~r~o~o~f~.~~The completeness of
the sequence~(\ref{Nsequence}) follows immediately from
Theorem~\ref{Completeness}.
Regarding the $\omega$-linear independence of this
sequence, it suffices to prove
$\omega$-independence
for the subsequence $\{\cN_i^{(l)}\}_{i=i_0}^\infty$.
Suppose there is a sequence $\{x_i\}_{i=i_0}^\infty$,
\mbox{$x_i\in\sN_i^{(l)}$,} such that
\begin{equation}
\label{SumNormX}
  0<\Sum_{i=i_0}^\infty\|x_i\|^2<\infty
\end{equation}
but
\begin{equation}
\label{SumX}
\Lim_{n\to\infty}\Sum_{i=i_0}^n x_i=0.
\end{equation}
The condition~(\ref{SumNormX}) implies that there are nonzero
elements among the $x_i$, say, an element $x_k$, $k\geq i_0$. Since
the projection $\sQ^{(l)}_k$ is a continuous operator,
the equality
\begin{equation}
\label{QSumLimX}
\sQ^{(l)}_k\Lim_{n\to\infty}\Sum_{i=i_0}^n x_i=
\Lim_{n\to\infty}\Sum_{i=i_0}^n \sQ^{(l)}_k x_i
\end{equation}
holds. But, due to Eqs.~(\ref{QQ}),
$\sQ^{(l)}_k{x_i}=\delta_{ik}x_i$
and the r.\,h. side of~(\ref{QSumLimX}) gives $x_k$ while
the l.\,h. side gives zero, because of~(\ref{SumX}). Thus,
$x_k$ must be zero, too, and one comes to a contradiction
which means that the sequence~(\ref{Nsequence})
is $\omega$-linearly independent.%
{\nopagebreak\mbox{\phantom{MMMM}}\hfill $\Box$\par\addvspace{0.25cm}}

\begin{lemma}\label{NPHdimensions}
If, instead of~{\rm(\ref{DifLambda})}, the condition
\begin{equation}
\label{DifLambda4}
  \lambda_i^{(A_1)}-\lambda_{i-1}^{(A_1)}>2r>4r_0(B)\qquad
\forall i\geq i_0
\end{equation}
holds, then  \mbox{$\mathop{\rm dim}\cN_i^{(l)}=
\mathop{\rm dim}\sP_i^{(A_1)}\cH_1$,}
$i=0,i_0,i_0+1,\ldots$\,\,.
\end{lemma}

\noindent P~r~o~o~f~.~ The proof is based on ideas from the proof of
Theorem~V.4.15 in~\cite{Kato}.  Our goal is to show that the
differences $\sQ_i^{(l)}-\sP_i^{(A_1)}$, $i=0,i_0,i_0+1,\ldots$,
have norms smaller than unity. Obviously,
\begin{eqnarray}
\label{QminP0}
\|\sQ_i^{(l)}-\sP_i^{(A_1)}\| &=& \frac{1}{2\pi}
\biggl\|\Int_{\gamma_i}dz\,
(A_1-z)^{-1}X^{(l)}(H_1^{(l)}-z)^{-1}\biggr\|          \\
\nonumber
&&\leq\frac{1}{2\pi}\Int_{\gamma_i}|dz|\,
\|(A_1-z)^{-1}\|\,\|X^{(l)}\|\,\|(H_1^{(l)}-z)^{-1}\|.
\end{eqnarray}
First, we deal with $i=0$. Obviously, one can deform
the circle $\gamma_0$ in~(\ref{QminP0}) into the line
\mbox{$\Real z=\lambda_{i_0-1}+r$}, since, in the half-plane
\mbox{$\Real z\leq\lambda_{i_0-1}+r$}, the integrand behaves
like $1/z^2$ as \mbox{$z\to\infty$}. On this line,
\mbox{$\|(A_1-z)^{-1}\|\leq(r^2+\eta^2)^{-1/2}$} where
\mbox{$\eta=\Img z$}. At the same time
\mbox{$\|X^{(l)}\|\leq r_0(B)$}. Then it follows from
the identity~(\ref{ResIdentity}) that
\mbox{$\|(H_1^{(l)}-z)^{-1}\|\leq1/(\sqrt{r^2+\eta^2}-r_0)$},
$r_0\equiv r_0(B)$, and one obtains the estimate
\begin{eqnarray*}
\|\sQ_0^{(l)}-\sP_0^{(A_1)}\| &\leq&
\D\frac{r_0}{2\pi}\Int_{-\infty}^\infty d\eta\,
\D\frac{1}{\sqrt{r^2+\eta^2}(\sqrt{r^2+\eta^2}-r_0)} \\
 && =\D\frac{r_0}{\pi}\Int_0^\infty d\eta\,
\left(1+\D\frac{r_0}{\sqrt{r^2+\eta^2}}\right)\cdot
\D\frac{1}{r^2+\eta^2-r_0^2}.
\end{eqnarray*}
Estimating the fraction $r_0/\sqrt{r^2+\eta^2}$ by
$r_0/r$ and calculating the integral, one obtains
$$
\|\sQ_0^{(l)}-\sP_0^{(A_1)}\|\leq
\D\frac{r_0}{2}\left(1+\D\frac{r_0}{r}\right)
\D\frac{1}{\sqrt{r^2-r_0^2}}.
$$
Since $r>2r_0$, we find finally
$$
 \|\sQ_0^{(l)}-\sP_0^{(A_1)}\|\leq\D\frac{\sqrt{3}}{4}
<\D\frac{1}{2}.
$$

We can show, further, that
\begin{equation}
\label{QPdifsmall}
\|\sQ_i^{(l)}-\sP_i^{(A_1)}\|<1
\end{equation}
for $i\geq i_0$, too.
Indeed, for $z\in\gamma_i$, $i\geq i_0,$ we have
$$
 \|(A_1-z)^{-1}\|\leq \D\frac{1}{r},\qquad
\|(H_1^{(l)}-z)^{-1}\|\leq \D\frac{1}{r-r_0}\,.
$$
Substitution of these estimates into~(\ref{QminP0}) gives
$$
\|\sQ_i^{(l)}-\sP_i^{(A_1)}\|\leq\D\frac{r_0}{r-r_0}
$$
and, under the condition~(\ref{DifLambda4}),
the inequalities~(\ref{QPdifsmall}) hold true.

Thus, we have proved that for any $i=0,i_0,i_0+1,\ldots$
the estimate~(\ref{QPdifsmall}) is valid. But such an estimate implies
that the subspaces $\cN_i^{(l)}$ and $\sP_i^{(A_1)}\cH_1$
are isomorphic to each other (see, e.\,g.,~\cite{Kato}, \S\,4.6 of
Chapter~I) and, consequently,
\mbox{$\mathop{\rm dim}\cN_i^{(l)}=
\mathop{\rm dim}\sP_i^{(A_1)}\cH_1$.} The proof is complete.%
{\nopagebreak\mbox{\phantom{MMMM}}\hfill $\Box$\par\addvspace{0.25cm}}

\begin{theorem}\label{KatoV4-15-16}
Assume $\lambda_{i+1}^{(A_1)}-\lambda_i^{(A_1)}\to\infty$
as $i\to\infty.$ Let $i_0$ be a number such
that~{\rm(\ref{DifLambda4})} holds. Then the following limit exists
\begin{equation}
\label{sLimQ}
\mathop{s-\rm lim}\limits_{n\to\infty}
\sum\limits_{i=0,i\geq i_0}^n \sQ_i^{(l)}=I_1.
\end{equation}
Additionally, assume that
\begin{equation}
\label{InvSquareConvergent}
\sum\limits_{i=1}^\infty
{(\lambda_{i+1}^{(A_1)}-\lambda_i^{(A_1)})^{-2}}<\infty\,.
\end{equation}
Then~{\rm(\ref{sLimQ})} is true for any renumbering of
$\sQ_i^{(l)}$. Moreover, there exists a constant $C$ such that
$\biggl\|\Sum_{i\in\cI}\sQ_i^{(l)}\biggr\|\leq C$ for any finite
set $\cI$ of integers $i=0$, $i\geq i_0$.
\end{theorem}
This theorem represents a slightly extended statement of
Theorems~V.4.15 and~V.4.16 of~\cite{Kato} (which treated
only the case where all the eigenvalues $\lambda_i^{(A_1)}$
of the operator $A_1$ were simple). The proof of
Theorem~\ref{KatoV4-15-16} is realized in exactly the same
way as the proof of the mentioned theorems in~\cite{Kato}
and, thus, we omit it.

\begin{remark}\label{KatoV4-15-16-Rem1}
Eq.~{\rm(\ref{sLimQ})} implies that
\begin{equation}
\label{sLimP}
\mathop{s-\rm lim}\limits_{n\to\infty}
\Sum_{i=0,i\geq i_0}^n\,\,\,
\Sum_{\lambda\in\Inter\gamma_i} \sP_\lambda^{(l)}=I_1
\end{equation}
where $\lambda$ stand for the eigenvalues of the
operator $H_1^{(l)}$ and $\sP_\lambda^{(l)}$ for the respective
eigenprojections.
If, additionally, the
inequality~{\rm(\ref{InvSquareConvergent})} holds and all the
eigenvalues $\lambda_i^{(A_1)}$ are simple, then one can
renumber the eigenprojections $\sP_\lambda^{(l)}$ in
Eq.~{\rm(\ref{sLimP})} in any way {\rm(}see
{\rm Theorem~V.4.16 of~\cite{Kato}).}
\end{remark}

\begin{lemma}\label{XPdecrease}
As before, assume $\Delta_\sm=(\mu_\sm^{(1)},+\infty)$.
Also, suppose that there is a $K_B$-bounded contour
$\Gamma_l\subset D_l$ satisfying~{\rm(\ref{Best})} and such that a part
of its component $\Gamma_\sm^{l_\sm}$ coincides with the ray
$\tilde{\Delta}_\sm^{0}=[\mu_0,\ri b_0+\infty)$
where $\mu_0\in D_\sm^{l_\sm}$,
$\mu_0=a_0+\ri b_0$ with $a_0,b_0\in\R$. Additionally, suppose
that the remaining part
\mbox{$\tilde{\Gamma}_l=\Gamma_l\setminus\tilde{\Delta}_\sm^{0}$}
of the contour $\Gamma_l$ belongs to the half-plane
$\Real\mu<a_0$, and for $\mu\in\tilde{\Delta}_\sm^{0}$
\begin{equation}
\label{KBderEst}
\|K_B'(\mu)\|\leq\tilde{C}(1+|\Real\mu|)^{-\theta}
\end{equation}
with $\tilde{C}>0$ and $\theta>1$. Then the estimate
\begin{equation}
\label{XPEst}
\|X^{(l)}\sP_i^{(A_1)}\|\leq{C}{(1+|\lambda_i^{(A_1)}|)^{-1}},
\quad i=i_0,i_0+1,\ldots,
\end{equation}
is valid with some $C>0.$
\end{lemma}

\noindent P~r~o~o~f~.~~Since the basic equation~(\ref{MainEqC})
for $X^{(l)}$ can be written as
$$
X^{(l)}=\Int_{\sigma'(A_0)\cup\Gamma_l}
K_B(d\mu)\,[(A_1-\mu)^{-1}-
(H_1^{(l)}-\mu)^{-1}X_1^{(l)}(A_1-\mu)^{-1}],
$$
one finds
\begin{equation}
\label{XPfirst}
X^{(l)}\sP_i^{(l)}=\Int_{\sigma'(A_0)\cup\Gamma_l}
\D\frac{K_B(d\mu)\,T_i(\mu)}{\lambda_i^{(A_1)}-\mu}, \qquad i\geq i_0,
\end{equation}
with
$$
   T_i(\mu)=[I_1-(H_1^{(l)}-\mu)^{-1}X^{(l)}]\sP_i^{(A_1)}.
$$
Due to the estimates~(\ref{Xr0}) and~(\ref{HresolvEst}) the functions
$T_i(\mu)$ on $\sigma'(A_0)\cup\Gamma_l$ are bounded,
$\|T_i(\mu)\|\leq c$ where the constant $c$ is the same
for all $i\geq i_0$ being determined only by $r_0(B)$ and
$d_0(\Gamma_l)$. According to one of our assumptions,
the condition $\mu\in\tilde{\Gamma}_l$ implies that
$\Real\mu<a_0$ and, thus,
\mbox{$|\lambda_i^{(A_1)}-\mu|^{-1}<|\lambda_i^{(A_1)}-a_0|^{-1}$}
for sufficiently large $\lambda_i^{(A_1)}$,
$\lambda_i^{(A_1)}>a_0$. Since the condition~(\ref{Bbound})
is valid, this assertion immediately implies that the
estimate~(\ref{XPEst}) holds at least for the contribution
to $X^{(l)}\sP_i^{(A_1)}$ from the set
$\sigma'(A_0)\cup\tilde{\Gamma}_l$.

As to the contribution from the ray $\tilde{\Delta}_\sm^{0}$,
one writes $K_B(d\mu)$ on this ray as
$K_B(d\mu)=K'_B(\mu)d\mu$ and applies the
inequality~(\ref{KBderEst}). Here, the elementary estimate
\begin{equation}
\label{ElemEst}
\Int_0^{+\infty} dx\,\D\frac{1}{(1+x)^\theta(1+|x-\lambda|)}
\leq c(\theta)\,\D\frac{1}{1+|\lambda|},
\qquad \theta>1,\quad\lambda\in\R,
\end{equation}
is useful with a certain $c(\theta)>0$ depending only on $\theta$.
Using this inequality one easily finds that the
estimate~(\ref{XPEst}) holds for the contribution to
$X^{(l)}\sP_i^{(A_1)}$ from the ray $\tilde{\Delta}_\sm^{0}$, too,
and this completes the proof.%
{\nopagebreak\mbox{\phantom{MMMM}}\hfill $\Box$\par\addvspace{0.25cm}}

\begin{theorem}\label{X-TrClass}
Let, in addition to the conditions of {\rm Lemma \ref{XPdecrease}},
the condition~{\rm(\ref{InvSquareConvergent})}
be valid and the sequence of the differences
$\lambda_{i+1}^{(A_1)}-\lambda_{i}^{(A_1)}$
be monotone starting from some $i=k$, $k\geq1$, i.\,e.,
\begin{equation}
\label{DifLamMonot}
\lambda_{i+2}^{(A_1)}-\lambda_{i+1}^{(A_1)}\geq
\lambda_{i+1}^{(A_1)}-\lambda_{i}^{(A_1)}, \qquad i\geq k.
\end{equation}
Also, let
\begin{equation}
\label{dimPHbound}
n_i^{(A_1)}=
\mathop{\rm dim}\sP_i^{(A_1)}\cH_1\leq n_{\rm max}^{(A_1)}
\end{equation}
where $n_{\rm max}^{(A_1)}$ is a finite number,
the same for all $i=0,i_0,i_0+1,\ldots$\,\,.
Then $X^{(l)}$ is an operator of the trace class.
\end{theorem}

To prove this theorem we need the following simple auxiliary
statement.

\begin{lemma}
\label{SeriesAux}
Let a sequence $\{a_n\}_{n=1}^\infty$ have positive elements,
$a_n>0$,  $n\in\N$ and starting from some number $N$ be
monotone, i.\,e.  $a_{n+1}\geq a_n$ for $n\geq N$. Also, let
\begin{equation}
\label{Inva2Convergent}
\sum_{n=1}^\infty \frac{1}{a_n^2}<\infty\,.
\end{equation}
Then the series $\sum_{n=1}^\infty b_n$ with
$b_n=(a_1+a_2+\ldots+a_n)^{-1}$ is convergent.
\end{lemma}

\noindent P~r~o~o~f~.~~First, one notes that for $2k\geq 2N$
$$
b_{2k}=\frac{1}{a_1+\ldots+a_k+a_{k+1}+\ldots+a_{2k}}<
\frac{1}{a_{k+1}+\cdots+a_{2k}}\leq\frac{1}{k\,a_{k+1}}
$$
and, similarly,
$$
b_{2k+1}<\frac{1}{a_{k+1}+\cdots+a_{2k}+a_{2k+1}}
\leq\frac{1}{(k+1)\,a_{k+1}}\,.
$$
This means that for $m>N$
\begin{eqnarray}
\label{EstSeq2N2m}
\sum_{n=2N}^{2m}\,b_n &=&\sum_{k=N}^m \,b_{2k}+\sum_{k=N}^{m-1}b_{2k+1}\,
\leq \sum_{k=N}^{m}\,\frac{1}{k\,a_{k+1}}+\sum_{k=N}^{m-1}\,
\frac{1}{(k+1)\,a_{k+1}}\\
\nonumber
&\leq &\left(\sum_{k=N}^{m}\,\frac{1}{k^2}\right)^{1/2}
\left(\sum_{k=N}^{m}\,\frac{1}{a_{k+1}^2}\right)^{1/2}+
\left(\sum_{k=N}^{m-1}\,\frac{1}{(k+1)^2}\right)^{1/2}
\left(\sum_{k=N}^{m-1}\,\frac{1}{a_{k+1}^2}\right)^{1/2}\,.
\end{eqnarray}
Since the condition~(\ref{Inva2Convergent}) is assumed
and the series $\sum_{n=1}^\infty\,n^{-2}$ is convergent,
it immediately follows from~(\ref{EstSeq2N2m})
that the series $\sum_{n=1}^\infty b_n$ considered
is convergent, too, and this completes the proof.%
{\nopagebreak\mbox{\phantom{MMMM}}\hfill $\Box$\par\addvspace{0.25cm}}

\noindent P~r~o~o~f~~o~f~~Theorem \ref{X-TrClass}.
Under the conditions~(\ref{InvSquareConvergent})
and~(\ref{DifLamMonot})
the series of the inverse eigenvalues
of the entry $A_1$ is convergent
\begin{equation}
\label{LamTrace}
\sum\limits^{\infty}_{\mbox{\scriptsize$
\begin{array}{c}
                i=1\\
              \lambda_i^{(A_1)}\neq 0
\end{array}$}}
|\lambda_i^{(A_1)}|^{-1}<\infty\,.
\end{equation}
Indeed if one takes
$a_i=\lambda^{(A_1)}_{i+1}-\lambda^{(A_1)}_i$ then the sequence
$\{b_i\}_{i=1}^\infty$ of Lemma~\ref{SeriesAux} with
$b_i=(a_1+\ldots+a_i)^{-1}$ is represented just by
$b_i=1/\lambda^{(A_1)}_{i+1}$ (except the case
$\lambda^{(A_1)}_{i+1}=0$).  If all the eigenvalues
$\lambda^{(A_1)}_i$ and differences
$\lambda^{(A_1)}_{i+1}-\lambda^{(A_1)}_i$ are positive, then to
prove~(\ref{LamTrace}) one can immediately use the statement of
Lemma~\ref{SeriesAux}.  In the case of presence of a (finite)
number of negative $\lambda^{(A_1)}_i$ and/or of a (finite)
number of negative differences
$\lambda^{(A_1)}_{i+1}-\lambda^{(A_1)}_i$ one has to omit in the
sum in the l.\,h.  side of~(\ref{LamTrace}) all the negative
eigevalues $\lambda^{(A_1)}_i$ and/or all the eigenvalues
$\lambda^{(A_1)}_i$ generating negative differences
$\lambda^{(A_1)}_i-\lambda^{(A_1)}_j$ with $j<i$. Then, after
appropriate shift in numbering of the remaining eigenvalues
$\lambda^{(A_1)}_i$, Lemma~\ref{SeriesAux} can be applied and
thus, the inequality~(\ref{LamTrace}) will again hold true.

Further, consider the quantity
$$
\sum\limits_{i=0,i\geq i_0}^\infty\sum\limits_{j=1}^{n_i^{(A_1)}}
\|X^{(l)}\varphi_{ij}^{(A_1)}\|=
\sum\limits_{i=0,i\geq i_0}^\infty\sum\limits_{j=1}^{n_i^{(A_1)}}
\|X^{(l)}\sP_i^{(A_1)}\varphi_{ij}^{(A_1)}\|\leq
\sum\limits_{i=0,i\geq i_0}^\infty
\|X^{(l)}\sP_i^{(A_1)}\|
\sum\limits_{j=1}^{n_i^{(A_1)}}
\|\varphi_{ij}^{(A_1)}\|\,.
$$
Since the estimate~(\ref{XPEst})
as well as the condition~(\ref{dimPHbound}) hold true,
one finds
$$
\sum\limits_{i=0,i\geq i_0}^\infty\sum\limits_{j=1}^{n_i^{(A_1)}}
\|X^{(l)}\varphi_{ij}^{(A_1)}\|\leq
n_0^{(A_1)}\|X^{(l)}\sP_0^{(A_1)}\|+
C\,n_{\rm max}^{(A_1)}\,\sum_{i=i_0}^\infty
(1+|\lambda_i^{(A_1)}|)^{-1}.
$$
Due to the inequality (\ref{LamTrace}) the
above quantity is finite and, thus, according to
Theorem~\ref{CriterionTrace} the operator $X^{(l)}$ is indeed of
the trace class. The proof is complete.%
{\nopagebreak\mbox{\phantom{MMMM}}\hfill $\Box$\par\addvspace{0.25cm}}
\begin{remark}
\label{Rezolv1TrClass}
Under the condition~{\rm(\ref{dimPHbound})}, the
inequality~{\rm(\ref{LamTrace})} implies that the resolvent
$(A_1-z)^{-1}$ is of the trace class for any $z\in\rho(A_1)$
since the sum
$\sum\limits_{j=1}^{n_0^{(A_1)}}\|(A_1-z)^{-1}\varphi_{0j}^{(A_1)}\|$ is
finite while the series
$\sum\limits_{i=i_0}^\infty\sum
\limits_{j=1}^{n_i^{(A_1)}}\|(A_1-z)^{-1}\varphi_{ij}^{(A_1)}\|\leq
\sum\limits_{i=i_0}^\infty\,n_i^{(A_1)}{|\lambda^{(A_1)}_i-z|}^{-1}$
is convergent {\rm(see Theorem~\ref{CriterionTrace})}.
\end{remark}
\begin{theorem}\label{BasisPropertyInf}
Let the condition~{\rm(\ref{InvSquareConvergent})} and
the conditions of {\rm Lemma \ref{XPdecrease}} be valid.
Then the inequality~{\rm(\ref{Cosines})} holds
for the subspaces $\cN_i^{(l)}$, $i=0,i_0,i_0+1,\ldots$\,\,.
This implies that the sequence~{\rm(\ref{Nsequence})} forms
a basis of the space $\cH_1$, quadratically close to an orthogonal one.

If, additionally, the condition~{\rm(\ref{dimPHbound})} holds, then
the union of orthonormal vector bases of the subspaces
$\cN_i^{(l)}$, $i=0,i_0,i_0+1,\ldots,$ forms a Bari basis of
$\cH_1$.
\end{theorem}

\noindent P~r~o~o~f~.~~Let
$r_i=\frac{1}{2}\Min\{\lambda_{i+1}^{(A_1)}-\lambda_{i}^{(A_1)},
\lambda_{i}^{(A_1)}-\lambda_{i-1}^{(A_1)}\}$, $i\geq i_0.$
Under the condition~(\ref{InvSquareConvergent}),
$r_i\to\infty$ as $i\to\infty$ and, moreover,
\begin{equation}
\label{InvRsquare}
\Sum_{i=i_0}^\infty \D\frac{1}{r_i^2}<\infty\,.
\end{equation}
Denote by $\tilde{\gamma}_i$, $i\geq i_0$, the circle centered
at $z=\lambda_i^{(A_1)}$ and having the radius $r_i$. Consider
the difference $\sQ_i^{(l)}-\sP_i^{(l)}$ replacing $\gamma_i$ in
the definitions~(\ref{sQDef}) and~(\ref{sPADef}) for $i\geq i_0$
with $\tilde{\gamma}_i$. Applying the resolvent identities
$$
(H_1^{(l)}-z)^{-1}-(A_1-z)^{-1}=
-(H_1^{(l)}-z)^{-1}X^{(l)}(A_1-z)^{-1}=
-(A_1-z)^{-1}X^{(l)}(H_1^{(l)}-z)^{-1}
$$
twice gives
\begin{equation}
\label{QPDD}
\sQ_i^{(l)}=\sP_i^{(A_1)}+\sD'_i+\sD''_i
\end{equation}
where
$$
\sD'_i=\D\frac{1}{2\pi\ri}\Int_{\tilde{\gamma}_i} dz\,
(A_1-z)^{-1}X^{(l)}(A_1-z)^{-1}
$$
and
$$
\sD''_i=-\D\frac{1}{2\pi\ri}\Int_{\tilde{\gamma}_i} dz\,
(A_1-z)^{-1}X^{(l)}(H_1^{(l)}-z)^{-1}X^{(l)}(A_1-z)^{-1}\,.
$$
It is easy to obtain the following estimates
(see also the proof of Theorem~V.4.16 in~\cite{Kato}):
\begin{equation}
\label{DDestimates}
\|\sD'_i\|\leq c\,\D\frac{1}{1+r_i}, \qquad
\|\sD''_i\|\leq c\,\D\frac{1}{(1+r_i)^2}\,
\end{equation}
with some $c>0$, the same for all $i\geq i_0$.

Further, consider the minimal angle
$\phi(\cN_i^{(l)},\cN_j^{(l)})$ between the subspaces
$\cN_i^{(l)}$ and $\cN_j^{(l)}$, $i,j\geq i_0$, $i\neq j$.
To estimate this angle it suffices to evaluate the inner
product
$\lal\sQ_i^{(l)}x,\sQ_j^{(l)}y\ral$
for $x\in\cN_i^{(l)}$, $y\in\cN_j^{(l)}$, $\|x\|=\|y\|=1$.
Substituting~(\ref{QPDD}) one obtains
\begin{eqnarray}
\nonumber
|\lal\sQ_i^{(l)}x,\sQ_j^{(l)}y\ral| &\leq &
|\lal\sP_i^{(A_1)}x,\sD'_j y\ral+\lal\sD'_i x,\sP_j^{(A_1)}y\ral|
+|\lal\sP_i^{(A_1)}x,\sD''_j y\ral|+|\lal\sD''_i x,\sP_j^{(A_1)}y\ral| \\
\label{QxQy}
&&+|\lal\sD'_i x,\sD'_j y\ral|+|\lal\sD'_i x,\sD''_j y\ral|
+|\lal\sD''_i x,\sD'_j y\ral|+|\lal\sD''_i x,\sD''_j y\ral|\,.
\end{eqnarray}
The term $|\lal\sP_i^{(A_1)}x,\sP_j^{(A_1)}y\ral|$ is
absent in the r.\,h. side of Eq.~(\ref{QxQy}) since
\mbox{$\sP_i^{(A_1)}\sP_j^{(A_1)}=0$} for $i\neq j$.
Meanwhile, according to~(\ref{DDestimates}),
the last four terms
can be estimated together by
\begin{equation}
\label{crirj}
c\,\,\,\D\frac{1}{1+r_i}\cdot\D\frac{1}{1+r_j}
\end{equation}
with another constant $c$ which does not depend on $i,j$.
The estimation of the terms
\mbox{$|\lal\sP_i^{(A_1)}x,\sD''_j y\ral|$} and
\mbox{$|\lal\sD''_i x,\sP_j^{(A_1)}y\ral|$}
is simple, too. Consider, for example, the term
\mbox{$|\lal\sD''_i x,\sP_j^{(A_1)}y\ral|=
|\lal\sP_j^{(A_1)}\sD''_i x,y\ral|$}. For $j\geq i_0$
we find
\begin{equation}
\label{PjR1}
\sP_j^{(A_1)}(A_1-z)^{-1}=
\sP_j^{(A_1)}(\lambda_j^{(A_1)}-z)^{-1}
\end{equation}
and, thus,
$$
\lal\sD''_i x,\sP_j^{(A_1)}y\ral=-\D\frac{1}{2\pi\ri}\,\,
\,\Int_{\tilde{\gamma}_i} dz\,\,
\D\frac{\lal\sP_j^{(A_1)}X^{(l)}(H_1^{(l)}-z)^{-1}
X^{(l)}(A_1-z)^{-1}x,y\ral}{\lambda_j^{(A_1)}-z}\,.
$$
Since $i\neq j$, one observes that
\mbox{$|\lambda_j^{(A_1)}-z|\geq
|\lambda_j^{(A_1)}-\lambda_i^{(A_1)}|-r_i$} if
$z\in\tilde{\gamma}_i$.  Consequently,
$$
|\lal\sD''_i x,\sP_j^{(A_1)}y\ral|\leq
\D\frac{\|X^{(l)}\|^2}
{(|\lambda_j^{(A_1)}-\lambda_i^{(A_1)}|-r_i)\,(r_i-\|X^{(l)}\|)}\,.
$$
Since
\begin{equation}
\label{DifLambdaSumR}
|\lambda_j^{(A_1)}-\lambda_i^{(A_1)}|\geq r_i+r_j,
\end{equation}
the term $|\lal\sD''_i x,\sP_j^{(A_1)}y\ral|$ can be estimated
by~(\ref{crirj}), too, and the same estimate
holds for \mbox{$|\lal\sP_i^{(A_1)}x,\sD''_j y\ral|$}.

Regarding the first term on the r.\,h. side of~(\ref{QxQy}),
it can be greatly
simplified by using the identity~(\ref{PjR1}) and then
the Residue Theorem. As a result one finds
\begin{equation}
\label{QxQyFirstTerm}
\lal\sP_i^{(A_1)}x,\sD'_j y\ral+\lal\sD'_i x,\sP_j^{(A_1)}y\ral=
\D\frac{\lal(X^{(l)}-X^{(l)*})\sP_i^{(l)}x,\sP_j^{(l)}y\ral}
{\lambda_j^{(A_1)}-\lambda_i^{(A_1)}}\,.
\end{equation}
Applying the inequalities~(\ref{XPEst})
and~(\ref{DifLambdaSumR}) one concludes that the
term~(\ref{QxQyFirstTerm}) can be easily estimated
again by~(\ref{crirj}).  However, the
series~(\ref{InvRsquare}) is convergent. This just implies
that
\begin{equation}
\label{CosinesLJ0}
\Sum^\infty_{\mbox{\scriptsize$
\begin{array}{c}
                i,j=i_0\\
                i\neq j
\end{array}$}}
\cos^2\phi(\cN_i^{(l)},\cN_j^{(l)})<\infty\,.
\end{equation}
An almost literal repetition of the previous consideration shows
that
\begin{equation}
\label{CosinesL0I}
\Sum^\infty_{i=i_0}
\cos^2\phi(\cN_0^{(l)},\cN_i^{(l)})<\infty\,,
\end{equation}
too. The inequalities~(\ref{CosinesLJ0}) and~(\ref{CosinesL0I})
imply that the condition of Theorem~\ref{Markus1} holds and,
thus, the sequence~{\rm(\ref{Nsequence})} indeed forms a basis
of the space $\cH_1$, quadratically close to an orthogonal one.
The second statement of the theorem is now a trivial
consequence of Theorem~\ref{GKMarkus1}.
The proof is complete.%
{\nopagebreak\mbox{\phantom{MMMM}}\hfill $\Box$\par\addvspace{0.25cm}}

\bigskip

\section{The simplest example}\label{Simple-Example}
In the present Section we consider the operator
matrix~(\ref{twochannel}) with the entry $A_0$ being
the multiplication operator,
\begin{equation}
\label{A0Multiplication}
 (A_0 u_0)(\mu)=\mu\, u_0(\mu),
\end{equation}
considered in $\cH_0=L_2(0,a)$, $0<a\leq+\infty$.  The domain of the
operator $A_0$ is \mbox{$\cD(A_0)=\left\{u_0\in L_2(0,a):\,
\displaystyle\int_0^a d\mu\,\mu^2|u_0(\mu)|^2<\infty\right\}$}\,.
Surely, if $a<\infty$, then $\cD(A_0)=\cH_0$.  The spectrum
of $A_0$ only consists of absolutely
continuous spectrum coinciding with the interval $[0,a]$.

As $A_1$ we take a diagonal numerical matrix,
\begin{equation}
\label{A1Lambda}
A_1=\Lambda=\diag\{\lambda_1,\lambda_2,\ldots,\lambda_k,\ldots\},
\qquad \lambda_k\in\R,
\end{equation}
and $\cH_1=\C^n$, $1\leq n\leq\infty$ (by $\C^\infty$ we
understand the Hilbert space $l_2$). The domain of the entry
$A_1$ is given by
\mbox{$\cD(A_1)=\left\{u_1=(u_1^{(1)},u_1^{(2)},\ldots,u_1^{(k)},\ldots)
\in\C^n:\,{\displaystyle\sum}_{k=1}^n\,\lambda_k^2\,
u_1^{(k)2}<\infty\right\}$.} Of course, if $n<\infty$,
then $\cD(A_1)=\C^n$.

The coupling operator $B_{01}$ acts on
$u_1\in\C^n$, $u_1=(u_1^{(1)},u_1^{(2)},\ldots,u_1^{(k)},\ldots)^t$,
as
\begin{equation}
\label{B01Example}
  (B_{01}u_1)(\mu)\equiv B(\mu)u_1\equiv
  \Sum_{k=1}^n b_k(\mu)u_1^{(k)}
\end{equation}
where $b_k\in L_2(0,a)$, $k=1,2,\ldots$, while $B(\mu)$
stands for the matrix-row of the values $b_k(\mu)$
of the functions $b_k$ for a fixed $\mu\in(0,a)$.
Boundedness of the entry $B_{01}$ means
\begin{equation}
\label{B01bound}
   \|B_{01}\|^2=\Sup\limits_{\|u_1\|=1} \Int_0^a d\mu\,
|B(\mu)u_1|^2<\infty.
\end{equation}
Obviously, the adjoint operator $B_{10}=B_{01}^*$ is given by
$$
B_{10}=(\lal\cdot,b_1\ral,\lal\cdot,b_2\ral,\ldots,
\lal\cdot,b_k\ral,\ldots)^t.
$$
Under the condition~(\ref{B01bound}) the operator $\bH$ is
selfadjoint on the domain $\cD(\bH)=\cD(A_0)\oplus\cD(A_1)$.

Note that this example is sufficiently universal. In particular
the Hamiltonians for the quantum-mechanical two-body systems
with internal structure used in
Refs.~\cite{Dashen,JaffeLow,YaF88,KMMMP,MotJMPh91,
Pavlov1984,PavlovShushkov,Simonov} can be reduced
to just the present example.  Note also that for
$n=1$ the operator matrix $\bH$ described represents one of the
well known Friedrichs models~\cite{Fried}.

The definition~(\ref{A0Multiplication}) of the entry $A_0$
represents at the same time its spectral decomposition,
that is, the spectral function $E^0(\mu)$ is given by (see,
e.\,g.,~\cite{BS})
\begin{equation}
\label{MeasureMult}
   \biggl(E^0(\mu)u_0\biggr)(\nu)=\left\{
\begin{array}{cl}
   u_0(\nu)\,, & \nu\leq\mu   \\
     0\,       &  \nu>\mu
\end{array}
\right.
\end{equation}
for $0<\mu\leq a$ and $E^0(\mu)=0$ for $\mu\leq 0$
while $E^0(\mu)=I_0$ for $\mu>a$. Consequently,
the product~(\ref{KBproduct}) for $0\leq\mu\leq a$
reads
$$
    K_B(\mu)=\Int_0^\mu d\nu\,[B(\nu)]^* B(\nu)
$$
and formally
\begin{equation}
\label{KBderModel}
    K'_B(\mu)=[B(\mu)]^* B(\mu)\,.
\end{equation}
In particular, if $n=1$, then $K'_B(\mu)=|b_1(\mu)|^2$.

One of our central assumptions of Sect.~\ref{Transfer_functions}
was the assumption regarding the holomorphy of the function
$K'_B(\mu)$ in a vicinity of (a part of) $\sigma_c(A_0)$.
Thus, in the case here one has to assume that the function $K'_B$
given by~(\ref{KBderModel}) takes values
in $\bB(\C^n,\C^n)$ for any $\mu\in[0,a]$ and admits
analytic continuation as
\mbox{$K'_B:\,D\to\bB(\C^n,\C^n)$} on
a domain $D$, $D\supset(0,a)$, symmetric with respect
to the real axis, that is
$D=D^-\cup D^+\cup(0,a)$ with $D^\pm\subset\C^\pm$,
$D^\pm=\{z:\,z\in D^\mp\}$, and
\begin{equation}
\label{N0BExample}
 \cV_0(B)=\Int_0^a d\mu\,\biggl\|[B(\mu)]^* B(\mu)\biggr\|<\infty\,.
\end{equation}
Obviously, the condition~(\ref{N0BExample})  implies
the inequality~(\ref{B01bound}).

The transfer function $M_1(z)$ for the model concerned reads
$$
M_1(z)=\Lambda-z+V_1(z)\qquad\mbox{with}\qquad
V_1(z)=\Int_0^a d\mu\,\D\frac{K'_B(\mu)}{z-\mu}.
$$
Consider the quantity
$$
V_1^{(0)}\mathop{\mbox{\large$=$}}\limits^{\rm def}
\Sup\limits_{\|u_1\|=1}(-\lal V_1(0)u_1,u_1\ral)=
\Sup\limits_{\|u_1\|=1} \Int_0^a d\mu\, \D\frac{|B(\mu)u_1|^2}{\mu}
$$
and, in the case of a finite $a$, the quantity
$$
V_1^{(a)}\mathop{\mbox{\large$=$}}\limits^{\rm def}
\Sup\limits_{\|u_1\|=1}\lal V_1(a)u_1,u_1\ral=
\Sup\limits_{\|u_1\|=1} \Int_0^a d\mu\, \D\frac{|B(\mu)u_1|^2}{a-\mu}.
$$
Denote by $\lambda_{\rm min}$ and $\lambda_{\rm max}$,
respectively, the lower and upper bounds for $\Lambda$,
$
\lambda_{\rm min}=\Inf\limits_{\|u_1\|=1}
\lal\Lambda u_1,u_1\ral
$
and
$
\lambda_{\rm max}=\Sup\limits_{\|u_1\|=1}
\lal\Lambda u_1,u_1\ral.
$
In the case of a finite $n$, $\lambda_{\rm min}$ and
$\lambda_{\rm max}$ coincide, respectively, with the minimal and
maximal eigenvalues of $\Lambda$.

Considering the quadratic form $\lal M_1(z)u_1,u_1\ral$
for $z<0$ and then, if $a$ is finite, for $z>a$
one can easily check that the following assertion holds true.

\begin{lemma}\label{Lem1-Example}
If $V_1^{(0)}$ is finite and
\begin{equation}
\label{CondLam-0}
V_1^{(0)}<\lambda_{\rm min},
\end{equation}
then the operator $\bH$ has no spectrum below $z=0$.
If, in the case of a finite $a$, $V_1^{(a)}$ is finite and
\begin{equation}
\label{CondLam-a}
V_1^{(a)}<a-\lambda_{\rm max},
\end{equation}
then the operator $\bH$ has no spectrum above $z=a$.
\end{lemma}
\begin{remark}\label{Rem1-Example}
The inequalities $V_1^{(0)}<\infty$  and $V_1^{(a)}<\infty$ imply,
respectively, \mbox{$K'_B(0)=0$} and \mbox{$K'_B(a)=0$}.
\end{remark}
\begin{remark}\label{Rem2-Example}
Suppose $n=1$ and, thus, $\lambda_{\rm min}=\lambda_{\rm
max}=\lambda_1$.  Considering the graphs of the functions
$y=z-\lambda$ and $y=V_1(z)$ at $z<0$ one can easily check that
if $\lambda_1\in(0,a)$, then for $V_1^{(0)}>\lambda_1$ the
transfer function $M_1(z)$ has a single negative eigenvalue.
Respectively, considering the graphs of the same functions at $z>a$
one observes that for $V_1^{(a)}>a-\lambda_1$ there exists a
single eigenvalue situated to the right from $a$.
\end{remark}

Therefore if the conditions~(\ref{CondLam-0}) and~(\ref{CondLam-a})
are valid the entire spectrum of the operator $\bH$ (and, hence,
the spectrum of the transfer function $M_1(z)$) must belong
to the interval $[0,a]$. Meanwhile, the eigenvalues of $\Lambda$
which are embedded initially into the continuous spectrum of $A_0$
can survive in this interval only in exceptional cases. Obviously,
if $\lambda\in(0,a)$ and $M_1(\lambda\pm\ri 0)u_1=0$, $u_1\neq 0$,
then the following conditions must hold (cf.~Lemma~\ref{LReal3},
condition (d), and Eq.~(\ref{MeasureMult})):
\begin{equation}
\label{BLamZero}
\lal K'_B(\lambda)u_1,u_1\ral=\|B(\lambda)u_1\|^2=0,
\end{equation}
$$
\lal(\Lambda-\lambda)u_1,u_1\ral+\mathop{\rm V.p.}
\Int_0^a d\mu\,\D\frac{\lal K'_B(\mu)u_1,u_1\ral}{\lambda-\mu}=0\,.
$$
These conditions may be hard to satisfy. In particular,
for $n=1$ Eq.~(\ref{BLamZero}) implies $b_1(\lambda)=0$. And if
one knows that $b_1(\mu)\neq 0$ for any $\mu\in(0,a)$, then no
point spectrum of $\bH$ can be situated in the interval $(0,a)$.
One understands that the embedded eigenvalue does not disappear.
It simply shifts into the unphysical sheet(s) and turns into a
pair of conjugate resonances which are eigenvalues of the
continued transfer function $M_1(z)$.

Suppose that there exist $K_B$-bounded contours
$\Gamma^\pm\subset D^\pm$ (see Sect.~\ref{Transfer_functions})
such that the condition~(\ref{Best}) holds. In this case these
are contours for which
\begin{equation}
\label{SolvCondEx}
\cV_0(B,\Gamma^\pm)=\Int_{\Gamma^\pm}
|d\mu|\,\|K'_B(\mu)\|<\D\frac{1}{4}\,d_0^2(\Gamma^\pm)
\end{equation}
where $d_0(\Gamma^\pm)=\dist\biggl\{\Gamma^\pm,
\{\lambda_k\}_{k=1}^n\biggr\}$\,.
Then, according to Theorem~\ref{Solvability},
one can construct two operators $H_1^{(+)}$ and $H_1^{(-)}$
the spectrum of which exhausts the spectrum of the respective
(continued) transfer functions $M_1(z,\Gamma^+)$ and
$M_1(z,\Gamma^-)$ in the set ${\cal O}_{d_{\rm max}/2}(A_1)$
where $d_{\rm max}$ is given by~(\ref{dmax}).

Finally, we give an illustration for the assertion of
Theorem~\ref{Solvability} for the simplest case of the model
(\ref{A0Multiplication})--(\ref{B01Example}) with $n=1$, $a=2R$,
$\lambda_1=R$ and $b_1(\mu)\equiv\beta$ where $R,\beta$ are some
positive numbers, $R,\beta\in{\R}^+$.  In this case the
basic equation~(\ref{MainEq}) coincides with the equation
$M_1(z,\Gamma^\pm)=0$ and the solutions $H_1^{(\pm)}$ if they
exist are operators in $\C$ defined by multiplication by
respective resonances.  Obviously
$$
\cV_0(B,\Gamma^\pm)=\beta^2\Int_{\Gamma^\pm}
|d\mu|=\beta^2\,\ell_{\Gamma^\pm}.
$$
Let ${\Gamma^\pm}$ be, say, semicircles,
${\Gamma^\pm}=\{z:\, |z-R|=R,\, z\in{\C}^\pm\}$. Then
$\cV_0(B,\Gamma^\pm)=\pi\,\beta^2 R$ and $d_0=d_0(\Gamma^\pm)=R$.
Thus, the solvability condition~(\ref{SolvCondEx}) reads now
as $\beta^2<\D\frac{R}{4\pi}$.
This means that Theorem~\ref{Solvability}
guarantees the unique solvability of the basic equation~(\ref{MainEq})
in any semidisc ${\cal S}_r^\pm\subset{\C}^\pm$ of radius $r$
centered at the point $z=R$, with $r$ satisfying the inequalities
$r_{\rm min}\leq r < r_{\rm max}$ where
\begin{eqnarray*}
r_{\rm min}&=&\D\frac{d_0}{2}-
\sqrt{\D\frac{d_0^2}{4}-\cV_0(B,\Gamma^\pm)}=
\D\frac{R}{2}-\sqrt{\D\frac{R^2}{4}-\pi\beta^2R}<\D\frac{R}{2} \\
r_{\rm max} &=& d_0-\sqrt{\cV_0(B,\Gamma^\pm)}=
R-\sqrt{\pi\beta^2R}>\D\frac{R}{2}.
\end{eqnarray*}
For instance, if $\beta^2=\D\frac{3}{16}\D\frac{R}{\pi}$ then
$r_{\rm min}=\D\frac{R}{4}$ and $r_{\rm
max}=R\left(1-\D\frac{\sqrt{3}}{4}\right)\approx0.6R$.  In this
case the solution (number) $H_1^{\pm}$ belongs to the semidisc
$|z-R|\leq\D\frac{R}{4}$, $z\in{\C}^\pm$, and no other solutions
exist in the semidisc $|z-R|\leq
R\left(1-\D\frac{\sqrt{3}}{4}\right)$, $z\in{\C}^\pm$.

In fact the model (\ref{A0Multiplication})--(\ref{B01Example})
with $n=1$, $0<a<\infty$ and $b_1(\mu)\equiv\beta$,
$\beta\in{\R}^+$, allows to calculate
the function $V_1(z)$ in
explicit form:
\begin{equation}
\label{V1Explic}
    V_1(z)=\beta^2\ln\D\frac{z}{z-a}
\end{equation}
where the physical-sheet logarithm branch is chosen in
such a way that
$$
\left(\ln\D\frac{z}{z-a}\right)_{\rm phys}
=\ln|z|-\ln|z-a| \qquad\mbox{for $z>a$.}
$$
The expression~(\ref{V1Explic}) gives an opportunity to treat
the physical as well as unphysical sheets of the transfer function
$M_1(z)=\lambda_1-z+V_1(z)$ immediately.

In the case considered both values $V_1^{(0)}$ and $V_1^{(a)}$ are
infinite. Thus the equation $M_1(z)=0$ has two roots in the
physical sheet (see Remark~\ref{Rem2-Example}), say $z_0$, $z_0<0$,
and $z_a$, $z_a>a$, representing eigenvalues of the operator $\bf H$.
One can even calculate the main terms of their asymptotics
as $\beta\to 0$:
$$
   z_0\sim -a\exp(-\lambda_1/\beta^2), \qquad
   z_a\sim a\,\{1+\exp[-(a-\lambda_1)/\beta^2]\}.
$$

Riemann surface of the function $M_1$ coincides with that of
$V_1$.  We denote unphysical sheets of this surface by
$\Pi_\nu$, $\nu=\pm1,\pm2,\ldots,$ assuming
$$
    \reduction{V_1(z)}{\Pi_\nu}=\beta^2
 \left[\left(\ln\D\frac{z}{z-a}\right)_{\rm phys}+2\pi\ri\,\nu\right]
$$
and that $\nu=0$ in this equation corresponds to the physical sheet
$\Pi_0$.

Put as previously $a=2R$, $\lambda_1=R$ and take
$\beta=\D\sqrt{\frac{R}{2}}\,\tilde{\beta}$.  We want to show
that the equation $\reduction{M_1(z)}{\Pi_\nu}=0$ has at least one
solution (resonance) in each unphysical sheet $\Pi_\nu$,
$\nu=\pm1,\pm2,\ldots$, with no restrictions on $\tilde{\beta}$
and $R$ such that $0<\tilde{\beta}<+\infty$, $0<R<+\infty$.

First, we note that for $z$ lying in the line $\Real z=R$,
$z=R\,(1+\ri\,\tan\varphi)$ where $0\leq\varphi<\D\frac{\pi}{2}$
or $\D\frac{3\pi}{2}<\varphi\leq 2\pi$, this equation can be
transformed into the following equations for the argument
$\varphi$:
\begin{equation}
\label{EqUpPhi}
\tan\varphi=\tilde{\beta}^2\left(\varphi-\D\frac{\pi}{2}+\pi\nu\right),
\quad \nu\in\Z,
\qquad\varphi\in\left[0,\D\frac{\pi}{2}\right)
\end{equation}
and
\begin{equation}
\label{EqDnPhi}
\tan\varphi=\tilde{\beta}^2\left(\varphi-\D\frac{3\pi}{2}+\pi\nu\right),
\quad \nu\in\Z,
\qquad\varphi\in\left(\D\frac{3\pi}{2},2\pi\right].
\end{equation}
Considering the graphs $y=\tan\varphi$ and
$y=\tilde{\beta}^2\left(\varphi-\D\frac{\pi}{2}+\pi\nu\right)$
for $0\leq\varphi<\D\frac{\pi}{2}$  one immediately checks that
Eq.~(\ref{EqUpPhi}) has no solutions for entire $\nu\leq0$ while
it necessarily gets a single root $\varphi_\nu$ for any entire
positive $\nu$, corresponding to a resonance
$z_\nu=R\,(1+\ri\,\tan\varphi_\nu)$ belonging to the upper
halfplane of the unphysical sheet $\Pi_\nu$ with
$\nu=1,2,3,\ldots$\,\,. At the same time, considering the graphs
$y=\tan\varphi$ and
$y=\tilde{\beta}^2\left(\varphi-\D\frac{3\pi}{2}+\pi\nu\right)$
for $\D\frac{3\pi}{2}<\varphi\leq2\pi$  one finds that
Eq.~(\ref{EqDnPhi}) has no solutions for entire $\nu\geq0$ and
necessarily gets a single root $\varphi_\nu$ for any entire
negative $\nu$, corresponding to a resonance
$z_\nu=R\,(1+\ri\,\tan\varphi_\nu)$
belonging to the lower halfplane of $\Pi_\nu$ with
$\nu=-1,-2,-3,\ldots$\,\,. (Also one observes that the
resonances $z_\nu$ and $z_{-\nu}$ are situated symmetrically
with respect to the real axis.)  Therefore we have proved that
the resonance set of the transfer function $M_1(z)$ in the model
considered is indeed nonempty in every unphysical sheet.


\appendix
\section{The norm of an operator with respect to a spectral measure}
\label{app1}
\small

Let $\cH',\cH''$ be separable Hilbert spaces, not necessarily
distinct, and $T\in\bB(\cH',\cH'')$.  Let $E$ be the spectral
measure of a self-adjoint operator in $\cH'$ with the support
$\sigma=\Supp E$, $\sigma\subset\R$. By the $E$-norm of the operator
$T$ we understand a number $\|T\|_E$ defined as
\begin{equation}
\label{Om}
\|T\|_E^2=\Sup_{\{\delta_k\}}\Sum_k \|TE(\delta_k)T^*\|
\end{equation}
where  $\{\delta_k\}$ stands for a finite or countable complete
system of pairwise nonintersecting subsets
of the set $\sigma$ measurable with respect to
$E$, i.\,e., $\delta_k$ are Borel subsets of $\sigma$, with
$\delta_k\cap\delta_l=\emptyset$ if $k\neq l$ and
$\bigcup_k\delta_k=\sigma$. For $S\in\bB(\cH'',\cH')$ we define
$\|S\|_E\mathop{\mbox{\large$=$}}\limits^{\rm def} \|S^*\|_E$.

One can easily check that
\begin{equation}
\label{TTE}
 \|T\|\leq\|T\|_E\,.
\end{equation}
Indeed,
$\|T\|^2=\mbox{\large$\|$}|T|\mbox{\large$\|$}^2=\|TT^*\|.$
Since $\Sum_k E(\delta_k)=I'$ with $I'$ being the identity operator in
$\cH'$, we conclude that
$$
\|TT^*\|=\|T\Sum_k E(\delta_k)T^*\|\leq\Sum_k\|TE(\delta_k)T^*\|\,.
$$
From this we immediately obtain~(\ref{TTE}). The equality $\|T\|=\|T\|_E$
is attained if the support $\sigma$ of the measure $E$ consists
of a single point.

\begin{lemma}\label{TEnorm}
The following equalities are valid:
\begin{equation}
\label{TEnormEq}
\|T\|_E^2=\Sup_{\{\delta_k\}}\Sum_k \|TE(\delta_k)\|^2
=\Sup_{\{\delta_k\}}\Sum_k \|E(\delta_k)T^*\|^2\,.
\end{equation}
\end{lemma}

\noindent P~r~o~o~f.
To begin with we note that for any $E$-measurable set
$\delta$
$$
\|TE(\delta)T^*\|=\|TE(\delta)\cdot E(\delta)T^*\|
\leq\|TE(\delta)\|\cdot\|E(\delta)T^*\|.
$$
Since $[TE(\delta)]^*=E(\delta)T^*$, we have
$\|[TE(\delta)]^*\|=\|E(\delta)T^*\|$. Therefore,
\begin{equation}
\label{Eq1}
\|TE(\delta)T^*\|\leq\|TE(\delta)\|^2=\|E(\delta)T^*\|^2.
\end{equation}
On the other hand, for any $f\in\cH''$
$$
  \|E(\delta)T^*f\|^2=\lal E(\delta)T^*f,E(\delta)T^*f\ral=
\lal TE(\delta)T^*f,f\ral\leq\|TE(\delta)T^*\|\,\|f\|^2
$$
and this means
\begin{equation}
\label{NEq2}
\|TE(\delta)\|^2=\|E(\delta)T^*\|^2\leq\|TE(\delta)T^*\|.
\end{equation}
It follows from~(\ref{Eq1}) and~(\ref{NEq2}) that, in fact,
\begin{equation}
\label{TETeq}
\|TE(\delta)T^*\|=\|TE(\delta)\|^2=\|E(\delta)T^*\|^2\,.
\end{equation}
Further, we can take $\delta=\delta_k$ and sum in~(\ref{TETeq})
over $k$. Then we can take for resulting sums the
exact upper bounds. Finally, one finds that Eqs.~(\ref{TEnormEq})
are indeed valid. The proof is complete.%
{\nopagebreak\hfill $\Box$\par\addvspace{0.25cm}}

Obviously, $\|\alpha T\|_E=|\alpha|\,\|T\|_E\,.$

At the same time, if the operators $T_1,T_2:\cH'\rightarrow\cH''$
have finite $E$-norms then their sum $T_1+T_2$ has a finite
$E$-norm and
\begin{equation}
\label{T12sum}
\|T_1+T_2\|_E\leq\|T_1\|_E+\|T_2\|_E\,.
\end{equation}
Indeed,
$$
\begin{array}{rcl}
\left(\Sum_k\|(T_1+T_2)E(\delta_k)\|^2\right)^{1/2} &\leq &
\left(\Sum_k(\|T_1 E(\delta_k)\|+\|T_2 E(\delta_k)\|)^2\right)^{1/2} \\
 &\leq& \left(\Sum_k\|T_1 E(\delta_k)\|^2\right)^{1/2} +
 \left(\Sum_k\|T_2 E(\delta_k)\|^2\right)^{1/2}\,.
\end{array}
$$
Using the statement of Lemma~\ref{TEnorm},
we come immediately to~(\ref{T12sum}).

The condition $T=0$, if $\|T\|_E=0$, follows from the
inequality~(\ref{TTE}).

Hence, the function \mbox{$\|\cdot\|_E$} is indeed a norm.  Due 
to~(\ref{TTE}), the limit of each Cauchy sequence of operators 
from $\bB(\cH',\cH'')$ having finite $E$-norms and converging 
with respect to the norm \mbox{$\|\cdot\|_E$} is automatically 
an element of $\bB(\cH',\cH'')$ having finite $E$-norm, too. 
Therefore, the operators from $\bB(\cH',\cH'')$, having finite 
$E$-norms, constitute a Banach space.

If the spectral measure $E$ corresponds to a self-adjoint
operator having a simple pure discrete spectrum, then, evidently,
$\|T\|_E$ coincides with the Hilbert-Schmidt norm $\|T\|_2$,
$\|T\|_2=\Sum_{n=1}^\infty\|Te_n\|^2$, where $\{e_n\}$ is an
arbitrary orthonormal basis of $\cH'$.  In general we
have only the inequalities
\begin{equation}
\label{TTET2}
\|T\|\leq\|T\|_E\leq\|T\|_2\,.
\end{equation}

Along with the $E$-norm~(\ref{Om}) one may consider as well a whole
family of operator norms defined with respect to a spectral
measure $E$:
\begin{equation}
\label{Epnorm}
\|T\|_{p,E}=\left(\Sup_{\{\delta_k\}}
\Sum_k\|TE(\delta_k)\|^p\right)^{1/p}=
\left(\Sup_{\{\delta_k\}}
\Sum_k\|E(\delta_k)T^*\|^p\right)^{1/p}, \quad p\geq 1.
\end{equation}
The norm~(\ref{Om}) is a particular case of these norms for
$p=2$.  Many properties of the norms~(\ref{Epnorm})
are similar to those for the respective norms $\|\,\cdot\,\|_p$ on
classes of compact operators.  Note, in particular, that if, for
$T_1:\cH'\rightarrow\cH''$, $T_2:\cH''\rightarrow\cH'$ and
${1}/{p}+{1}/{q}=1$, the norms $\|T_1\|_{p,E}$ and
$\|T_2\|_{q,E}$ are finite, then
$$
\Sup_{\{\delta_k\}}\Sum_k\|T_1 E(\delta_k)T_2\| \leq
\|T_1\|_{p,E}\cdot\|T_2\|_{q,E}\,.
$$
Note also that for any $1\leq p<\infty$
$$
    \|T\|\leq\|T\|_{p,E}\leq\|T\|_p.
$$
We do not describe properties of $\|\cdot\|_{p,E}$ for
arbitarary $p$, since, in this work, we use only the norm
$\|\cdot\|_E\equiv\|\cdot\|_{2,E}$.

\section{The integral of an operator-valued function over a
         spectral measure}
\label{IntOpMer}
\small

To avoid confusion with the measure $E_j$, we shall denote
the spectral function of the self-adjoint operator
$A_j$, $j=0,1$, by $E^j(\mu)$
(i.\,e., with superscript):
$E^j(\mu)=E_j\biggl((-\infty,\mu)\biggr)$, $\mu\in{\R}$.
Recall that $E^j(\mu)$ is a projection-valued function
satisfying the conditions of monotonicity, $E^j(\mu_1)\leq
E^j(\mu_2)$ for $\mu_1<\mu_2$, and completeness, $\mathop{s{-}\rm
lim}\limits_{\mu\rightarrow-\infty}E^j(\mu)=0$, $\mathop{s{-}\rm
lim}\limits_{\mu\rightarrow+\infty}E^j(\mu)=I_j$.  In addition,
this function is continuous from the left, $\mathop{s{-}\rm
lim}\limits_{\mu'\uparrow\mu}E^j(\mu')=E^j(\mu).$

Let $F(\mu)$ be a function defined on an interval $[a,b]$,
$-\infty<a<b<+\infty$, whose values are bounded operators acting
from $\cH_j$ to $\cH_i$, that is,
$F:[a,b]\rightarrow\bB(\cH_j,\cH_i)$, $i,j=0,1$ and where it is not
necessary that $i\neq j$. Following~\cite{AdLMSr}, we say the function
$F$ is uniformly (strongly, weakly) integrable from the right
over the spectral measure $E_j$ on $[a,b)$ if the limit
\begin{equation}
\label{defIntR}
\Int_a^b F(\mu)\,dE^j(\mu)\,\,\,
\mathop{\mbox{\large$=$}}\limits^{\rm def}
\Lim_{\Max_{k=1}^n |\delta_k^{(n)}|\rightarrow0}\,\,
\sum\limits_{k=1}^n  F(\xi_k)\, E_j(\delta_k^{(n)})
\end{equation}
exists considered in the sense of the uniform (strong, weak) operator
topology. Here, $\delta_k^{(n)}=[\mu_{k-1},\mu_k)$ and
$|\delta_k^{(n)}|=\mu_k-\mu_{k-1}$,  $k=1,2,\ldots,n$, where
$\mu_0,\mu_1,\ldots,\mu_n$ is any subsequence of numbers from the
interval $[a,b]$ satisfying the conditions
$a=\mu_0<\mu_1<\ldots<\mu_n=b$. By $\xi_k$ we understand an
arbitrary point of $\delta_k^{(n)}$. The limit
value~(\ref{defIntR}), if it exists, is called the right integral
of the function $F$ over the measure $E_j$ on $[a,b)$ in the sense
of Riemann-Stieltjes.

Similarly, we say the function
$G:[a,b]\rightarrow\bB(\cH_i,\cH_j)$, $i,j=0,1$  is uniformly
(strongly, weakly) integrable from the left over the spectral
measure $E_j$ on $[a,b)$, $-\infty<a<b<+\infty$, if it exists the
limit
\begin{equation}
\label{defIntL}
\Int_a^b dE^j(\mu)\,G(\mu)\,\,\,
\mathop{\mbox{\large$=$}}\limits^{\rm def}
\Lim_{\Max_{k=1}^n |\delta_k^{(n)}|\rightarrow0}\,\,
\sum\limits_{k=1}^n E_j(\delta_k^{(n)}) G(\xi_k)
\end{equation}
considered in the sense of the uniform (strong, weak) operator
topology.  The limit value~(\ref{defIntL}), if it exists, is called
the left integral of the function $G$ over measure $E_j$ on
$[a,b)$ in the sense of Riemann-Stieltjes.

Since the Banach spaces $\bB(\cH_i,\cH_j)$ and
$\bB(\cH_j,\cH_i)$ are closed with respect not only to the uniform
$(u)$ but also to the strong $(s)$ and weak $(w)$ convergence of
operators, the integrals~(\ref{defIntR}) and~(\ref{defIntL}), if
they exist in some sense, determine certain bounded operators
belonging, respectively, to $\bB(\cH_i,\cH_j)$ and $\bB(\cH_j,\cH_i)$.

Evidently, if the integrals~(\ref{defIntR}) and~(\ref{defIntL})
exist in the sense of the strong operator topology, they exist
as well in the sense of the weak operator topology. In turn,
existence of these integrals in the sense of the uniform operator
topology implies their existence in the sense of the strong as
well as weak operator topology. The following
simple statement holds.

\begin{lemma}\label{Integr1}
The function $F(\mu)$, $F:[a,b]\rightarrow\bB(\cH_j,\cH_i)$, is
integrable in the sense of the uniform operator topology over
the measure $E_j$ in $[a,b)$ from the left iff the function
$[F(\mu)]^*$ is integrable over the same measure from the right,
and also in the sense of the uniform operator topology.  The
same statement is true with respect to the simultaneous
integrability of these functions with respect to the weak
operator topology.  In general it only follows from the
existence of one of the
integrals $\Int_a^b F(\mu)\,dE^j(\mu)$ and $\Int_a^b
dE^j(\mu)\,[F(\mu)]^*$ in the sense of the strong operator topology
that the other one exists with respect to the weak
operator topology.  In all cases the integrability of $F(\mu)$
and $[F(\mu)]^*$ implies the equality
\begin{equation}
\label{JJadj}
\left[\Int_a^b F(\mu)\,dE^j(\mu)\right]^*
     =\Int_a^b dE^j(\mu)\,[F(\mu)]^*.
\end{equation}
\end{lemma}

\noindent P~r~o~o~f. Note that
$$
\left(\sum\limits_{k=1}^n  F(\xi_k)\, E_j(\delta_k^{(n)})\right)^*
=\sum\limits_{k=1}^n   E_j(\delta_k^{(n)})\,[F(\xi_k)]^*\,.
$$
Therefore, the validity of the first and second statements
of the lemma follows from continuity of the involution
$T\rightarrow T^*$ with respect to
$u$- and $w$-convergence in $\bB(\cH_j,\cH_i)$.
The last statement follows from the fact
that strong convergence of a sequence of operators in
$\bB(\cH_j,\cH_i)$ implies also weak convergence of this
sequence. However, one can not claim that the sequence
of respective adjoint operators converges strongly,
since the involution $T\rightarrow T^*$ is not continuous
with respect to $s$-convergence (see, e.\,g.,~\cite{BS},
\S\,5 of Chapter~2).

The proof of the lemma is complete.%
{\nopagebreak\mbox{\phantom{MMMM}}\hfill $\Box$\par\addvspace{0.25cm}}

\medskip

Some sufficient conditions for the integrability of an
operator-valued function $F(\mu)$ in the sense of the uniform
operator topology are given in the following statement.
\begin{lemma}\label{Integr2}
Any operator function $F$, $F:[a,b]\rightarrow\bB(\cH_i,\cH_j)$,
which satisfies the Lipschitz condition
\begin{equation}
\label{Lipschitz}
\|F(\mu_2)-F(\mu_1)\|\leq C_F\,|\mu_2-\mu_1| \qquad
\forall\mu_1,\mu_2\in[a,b]
\end{equation}
with some constant $C_F>0$, is
right-integrable with respect to $E_j$ in the sense of the
operator norm topology.
\end{lemma}

A proof of this statement can be found in Ref.~\cite{AdLMSr}
(see in~\cite{AdLMSr} Lemma~7.2 and Remark~7.3).

\medskip

The integrals $\Int_a^b F(\mu)\,dE^j(\mu)$ and
$\Int_a^b dE^j(\mu)\,G(\mu)\,$ with $a=-\infty$
or $b=+\infty$ are understood as respective limits,
if they exist, of integrals with finite bounds;
for example,
$$
\Int_a^b dE^j(\mu)\,G(\mu)\mathop{\mbox{\large$=$}}\limits^{\rm def}
\Lim_{a'\downarrow a,\, b'\uparrow b}
\Int_{a'}^{b'} dE^j(\mu)\,G(\mu).
$$

Also, we define
\begin{equation}
\label{DefIntSp}
\Int_{\sigma(A_j)} dE_j(\mu)\,G(\mu)
\mathop{\mbox{\large$=$}}\limits^{\rm def}
\Int_a^b dE^j(\mu)\,G(\mu), \qquad
\Int_{\sigma(A_j)} F(\mu)\,dE_j(\mu)
\mathop{\mbox{\large$=$}}\limits^{\rm def}
\Int_a^b F(\mu)\,dE^j(\mu)
\end{equation}
where $(a,b)$ is an arbitrary open interval entirely containing
the set $\sigma(A_j)$.  This definition is correct, since the
support of the spectral measure $E_j$ is just the
spectrum $\sigma(A_j)$.

\begin{lemma}\label{IntegrN}
Let a function $X:\sigma(A_j)\rightarrow\bB(\cH_i,\cH_i)$ be
bounded,
$\|X\|_\infty=\Sup\limits_{\mu\in\sigma(A_j)}\|X(\mu)\|<\infty$,
and satisfy the Lipschitz condition~{\rm(\ref{Lipschitz})}.
Then, if the $E_j$-norm $\|B_{ij}\|_{E_j}$ of the operator
$B_{ij}$ is finite, then the integrals
$$
\Int_{\sigma(A_j)} dE_j(\mu)\,B_{ji}\,X(\mu)
\quad
\mbox{and}
\quad
\Int_{\sigma(A_j)} X(\mu)\,B_{ij}\,dE_j(\mu)
$$
exist%
\footnote{\small The function $X$ can be extended outside the set
$\sigma(A_j)$ in an arbitrary way when the
definitions~(\ref{defIntR}), (\ref{defIntL})
and~(\ref{DefIntSp}) are used, retaining only the Lipschitz
condition~(\ref{Lipschitz}).}
in the sense of the operator norm topology, and
the following estimates are valid for their norms:
$$
\biggl\|\Int_{\sigma(A_j)} dE_j(\mu)\,B_{ji}\,X(\mu)\biggr\|
\leq\|B_{ji}\|_{E_j}\cdot\|X\|_\infty
\quad
\mbox{and}
\quad
\biggl\|\Int_{\sigma(A_j)} X(\mu)\,B_{ij}\,dE_j(\mu)\biggr\|
\leq\|B_{ji}\|_{E_j}\cdot\|X\|_\infty.
$$
\end{lemma}

\noindent P~r~o~o~f~.~ The proof
will be given for the case of the integral
$\Int_{\sigma(A_j)} dE_j(\mu)\,B_{ji}\,X(\mu)$. To this end let
us consider a partition
$\left\{\delta_k^{(n)}\right\}_{k=1}^{n}$ of an interval $[a',b')$
for finite $a'$, $b'$ and the respective integral sum. We have
for this sum:
$$
\left\|\Sum_{k=1}^n E_j(\delta_k^{(n)})B_{ji}X(\xi_k)f\right\|^2=
\lal\Sum_{k=1}^n E_j(\delta_k^{(n)})B_{ji}X(\xi_k)f,
\Sum_{m=1}^n E_j(\delta_m^{(n)})B_{ji}X(\xi_m)f\ral
$$
where $\left\{\xi_k\in\delta_k^{(n)}\right\}_{k=1}^{n}$
is an arbitrary set of points belonging to the intervals
$\delta_k^{(n)}$. Since $E_j(\delta_k^{(n)})E_j(\delta_m^{(n)})=0$
for $\delta_k^{(n)}\cap\delta_m^{(n)}=\emptyset$, we find
\begin{eqnarray}
\nonumber
\left\|\Sum_{k=1}^n E_j(\delta_k^{(n)})B_{ji}X(\xi_k)f\right\|^2
 &=& \lal\Sum_{k=1}^n [X(\xi_k)]^* B_{ij}
E_j(\delta_k^{(n)})B_{ji}X(\xi_k)f,f\ral  \\
\nonumber
&\leq& \Sum_{k=1}^n \|B_{ij}E_j(\delta_k^{(n)})B_{ji}\|
\cdot\|X(\xi_k)\|^2\cdot\|f\|^2  \\
\nonumber
&\leq& \Reduction{\cV_j(B)}{[a',b')}\|X\|^2_\infty\|f\|^2.
\end{eqnarray}
This means
$$
\left\|\Sum_{k=1}^n E_j(\delta_k^{(n)})B_{ji}X(\xi_k)\right\|
\leq\sqrt{\Reduction{\cV_j(B)}{[a',b')}}\cdot\|X\|_\infty.
$$
Since the total variation $\cV_j(B)=\|B_{ji}\|_{E_j}^2$ is
supposed to be finite, we have
$$
\biggr\|\Int_{a'}^{b'} dE^j(\mu)\,B_{ji}\,X(\mu)\biggr\|
\mathop{\longrightarrow}\limits_{
\begin{array}{c}
a'\rightarrow+\infty \\
(b'>a')
\end{array}}
0,
$$
This means that the integral
$\Int_{a}^{+\infty} dE^j(\mu)\,B_{ji}\,X(\mu)$
with a finite lower bound $a$ converges with respect to the
uniform operator topology. Existence of the integral
$\Int_{-\infty}^{b} dE^j(\mu)\,B_{ji}\,X(\mu)$
with a finite $b$ can be proved in the same way.
The existence of the integral
$\Int_{\sigma'(A_j)}dE_j(\mu)\,B_{ji}\,X(\mu)$
with respect to the operator norm topology as well as its
norm estimate follow immediately from these results.%
{\nopagebreak\mbox{\phantom{MMMM}}\hfill $\Box$\par\addvspace{0.25cm}}

\bigskip

In the same way as the for the previous integrals one can define
and treat the integral
\begin{equation}
\label{IntXBEBY}
\Int_{\sigma(A_j)}X(\mu)B_{ij}dE_j(\mu)B_{ji}Y(\mu)
\end{equation}
with $X,Y$ like the function $X$ in Lemma \ref{IntegrN}
and the same $B_{ij}=B_{ji}^*$. First, we extend $X$
and $Y$ outside the set $\sigma(A_j)$ retaining the
Lipschitz condition~(\ref{Lipschitz}) and introduce
for $-\infty<a<b<+\infty$ the value
\begin{equation}
\label{DefXBEBYint}
\Int_a^b X(\mu)B_{ij}dE^j(\mu)B_{ji}Y(\mu)\,\,\,
\mathop{\mbox{\large$=$}}\limits^{\rm def}
\Lim_{\Max_{k=1}^n |\delta_k^{(n)}|\rightarrow0}\,\,
\sum\limits_{k=1}^n  X(\xi_k)B_{ij}E_j(\delta_k^{(n)})B_{ji}Y(\xi_k)
\end{equation}
with $\delta_k^{(n)}$ and $\xi_k^{(n)}$ taken
as in~(\ref{defIntR}). Then we consider the limits
$a\to-\infty$ and/or $b\to+\infty$ if necessary.
As a result one has the following

\begin{lemma}\label{LIntXBEBY}
Let the functions $X,Y:\,\sigma(A_j)\to\bB(\cH_i,\cH_i)$
be bounded, $\|X\|_\infty<\infty,$ $\|Y\|_\infty<\infty,$ and satisfy
the Lipschitz condition~{\rm(\ref{Lipschitz})}
and let $\|B_{ij}\|_{E_j}<\infty$, also. Then the
integral~{\rm(\ref{IntXBEBY})}
exists in the sense of the operator norm topology and
$$
\biggl\|\Int_{\sigma(A_j)}X(\mu)B_{ij}dE_j(\mu)B_{ji}Y(\mu)\biggr\|
\leq \|B_{ij}\|_{E_j}^2\cdot\|X\|_\infty\cdot\|Y\|_\infty\,.
$$
\end{lemma}


\end{document}